\documentclass[preprint,3p,12pt]{elsarticle}
\usepackage{mathrsfs}
\usepackage{amsmath}
\usepackage{stmaryrd}
\usepackage{bbding}
\usepackage{dcolumn}
\usepackage{graphicx}
\usepackage{amsfonts}
\usepackage{amssymb}
\usepackage{wrapfig}
\usepackage{subfigure}
\usepackage{makeidx}
\usepackage{bm}
\usepackage{epsf}
\usepackage{color}
\usepackage{epsfig}
\usepackage{setspace}
\usepackage{graphicx}
\usepackage{epstopdf}
\usepackage{psfrag}
\usepackage{subfigure}
\usepackage{multirow}
\usepackage{diagbox}
\usepackage{verbatim}

\usepackage{makecell}
\usepackage{float}
\usepackage{esint}
\usepackage{comment}
\usepackage{movie15}
\usepackage{hyperref}
\usepackage[ruled,linesnumbered]{algorithm2e}
\usepackage{algorithmicx}
\usepackage{multirow}
\usepackage[noend]{algpseudocode}
\usepackage{amsthm}
\theoremstyle{remark}
\newtheorem{remark}{Remark}

\newcommand{\mathsym}[1]{{}}
\newcommand{\unicode}[1]{{}}
\begin{document}

    \color{red}
\title{A Geometric Multigrid-Accelerated Compact Gas-Kinetic Scheme for Fast Convergence in High-Speed Flows on GPUs}
\color{black}

	\author[XJTU,HKUST1]{Hongyu Liu}
    \ead{hliudv@connect.ust.hk}
	
	\author[XJTU]{Xing Ji\corref{cor1}}
	\ead{jixing@xjtu.edu.cn}

        \author[XJTU]{Yuan Fu}
	\ead{yuanfu@stu.xjtu.edu.cn}
	
	\author[HKUST1,HKUST2,HKUST3]{Kun Xu}
	\ead{makxu@ust.hk}

	\address[XJTU]{Shaanxi Key Laboratory of Environment and Control for Flight Vehicle, Xi'an Jiaotong University, Xi'an, China}
	\address[HKUST1]{Department of Mathematics, Hong Kong University of Science and Technology, Clear Water Bay, Kowloon, Hong Kong}
	\address[HKUST2]{Department of Mechanical and Aerospace Engineering, Hong Kong University of Science and Technology, Clear Water Bay, Kowloon, Hong Kong}
	\address[HKUST3]{Shenzhen Research Institute, Hong Kong University of Science and Technology, Shenzhen, China}
	\cortext[cor1]{Corresponding author}

\begin{abstract}

Implicit methods and GPU parallelization are two distinct yet powerful strategies for accelerating high-order CFD algorithms. However, few studies have successfully integrated both approaches within high-speed flow solvers. The core challenge lies in preserving the robustness of implicit algorithms in the presence of strong discontinuities, while simultaneously enabling massive thread parallelism under the constraints of limited GPU memory.
To address this, we propose a GPU-optimized, geometric multigrid-accelerated, high-order compact gas kinetic scheme (CGKS) that incorporates three key innovations:
 (1) a multi-color lower-upper symmetric Gauss-Seidel scheme that eliminates thread conflicts and preserves memory efficiency, serving as an implicit smoother on coarse grids; (2) a discontinuity-adaptive relaxation technique and a multigrid prolongation process, based on a discontinuous feedback factor, which dynamically stabilize shock regions without compromising convergence in smooth zones; and (3) a three-layer V-cycle geometric parallel multigrid strategy specifically tailored for unstructured meshes.
Extensive tests on multi-dimensional subsonic to hypersonic flows demonstrate that our GPU-based high-performance solver achieves one to two orders of magnitude faster convergence compared to previous explicit solvers. More importantly, it preserves the shock-capturing robustness of the explicit CGKS and exhibits strong scalability on GPU architectures. This work presents a unified framework that synergistically leverages implicit acceleration and GPU optimization for high-speed flow simulations, effectively overcoming traditional trade-offs between parallelism, memory constraints, and numerical stability in high-order methods.

\end{abstract}


\maketitle

\section{Introduction}

High-order methods with enhanced spatial resolution have become widely adopted in computational fluid dynamics, particularly for the accurate design of complex geometrical configurations. The extension of these methods to unstructured meshes further enhances their applicability to real-world engineering problems. However, traditional high-order finite volume and finite difference methods are challenging to implement on unstructured meshes due to their reliance on large stencils.
To address this limitation, compact high-order schemes—such as the discontinuous Galerkin (DG) method \cite{shu2003fd-fv-weno-dg-review}, flux reconstruction (FR) \cite{Huynh2007FR}, and the correction procedure via reconstruction (CPR) \cite{CPR2011}—have seen significant advancements in recent research. These methods offer high computational efficiency and spatial resolution in smooth regions. Nevertheless, their stability and robustness can be compromised in the presence of strong discontinuities \cite{cheng2009high}, which restricts their effectiveness in high-speed flow simulations.

The high-order compact gas-kinetic scheme (CGKS) \cite{ji2018compact,zhao2023high,zhao2019compact}, recently developed based on gas-kinetic theory, offers significant advancements in computational fluid dynamics. CGKS utilizes the time-dependent distribution function, which provides a time accurate evolution solution at the cell interface from the kinetic model eqaution \cite{BGK}. By leveraging this precise time-dependent distribution function, the Navier–Stokes flux functions can be accurately computed, enabling the evaluation of time-accurate macroscopic flow variables.
Within the finite volume framework, CGKS updates cell-averaged flow variables and, importantly, also updates cell-averaged slopes as additional degrees of freedom within each cell. This approach maximizes the amount of information obtained while minimizing the stencil size. Utilizing both the cell-averaged flow variables and their slopes, a Hermite Weighted Essentially Non-Oscillatory (HWENO) method \cite{li2021multi} can be employed for reconstruction, as demonstrated in our previous work on third-order CGKS \cite{zhang2023high}.
CGKS demonstrates excellent performance in both smooth and discontinuous flow regimes, benefiting from a reliable evolution process grounded in mesoscopic gas-kinetic theory. It also shows strong capabilities in unsteady compressible flow applications, such as computational aeroacoustics \cite{zhao2019acoustic} and implicit large eddy simulation \cite{ji2021compact}. In CGKS, the discontinuity feedback factor (DF) extends beyond the traditional finite volume framework by identifying the presence of discontinuities within a cell for the upcoming time step. This mechanism further enhances the robustness of CGKS in handling strong discontinuities, such as shocks and rarefaction waves.

However, simulating steady-state wall-bounded flows often requires millions of iterations to achieve convergence, resulting in substantial computational costs. This slow convergence is primarily due to the stringent time-step restrictions imposed by explicit time-marching schemes, which are dictated by the global minimum cell size—often as small as one millionth of the largest cell in the domain.
To address this, acceleration techniques such as local time stepping \cite{jeanmasson2019localtime} and LU-SGS \cite{li2017performance} are widely used in various CFD codes for their simplicity and robustness. These methods have also been applied to the gas-kinetic scheme (GKS) and the unified gas kinetic scheme (UGKS) \cite{zhu2016implicit}. More recently, the generalized minimal residual method (GMRES) has gained considerable attention for its rapid convergence rate, and CGKS has demonstrated its effectiveness in accelerating steady-state solutions \cite{yang2024Implicit}. However, the convergence rate of GMRES is highly sensitive to the choice of initial values and preconditioners; poor choices can lead to instability or even divergence.
Meanwhile, the rapid advancement of GPUs has led to significant breakthroughs in accelerating CFD solvers. Notable examples include HiFiLES \cite{hifi}, AMREX \cite{zhang2019amrex}, and ZEFR \cite{romero2020zefr}, all of which have leveraged GPU architectures to achieve remarkable gains in simulation speed and efficiency. Despite these successes, adapting implicit solvers to GPUs presents significant challenges. For instance, GMRES requires storing large, sparse matrices, which is difficult to achieve efficiently on GPUs with limited memory. Similarly, the widely used LU-SGS method may suffer from divergence due to massive thread contention. While multigrid methods are inherently parallelizable, their robustness in the presence of strong shocks can be limited. To date, there has been little success in adapting fully implicit all-flow solvers in GPUs.

Building on previous work \cite{liu2024compact}, we have developed a GPU-accelerated geometric multigrid method to achieve faster convergence and enhanced robustness. This development incorporates four main improvements and innovations.
(1). Efficient Smoothing via LU-SGS Filtering:
The choice of an appropriate filter is critical for accelerating convergence. In this work, we employ the LU-SGS method as our filter, which effectively damps high-frequency errors and serves as a matrix-free implicit solver, conserving precious GPU memory. To resolve the data dependency and thread-racing issues inherent in the standard LU-SGS method—which can destabilize computations—a universal rainbow coloring algorithm is introduced. This technique partitions computational cells into distinct groups by assigning different colors to neighboring cells, thereby eliminating thread conflicts and ensuring numerical stability \cite{bocharov2020implicit}.
(2). Discontinuity Feedback (DF) for Adaptive Multigrid Prolongation:
Unlike traditional geometric multigrid prolongation processes, which often require computationally expensive upwind stencil selection, we introduce the discontinuity feedback (DF) factor as a generalized discontinuity prediction tool within the multigrid cycle. Because DF adapts to the local flow field, the prolongation process becomes flow-adaptive, enhancing both convergence efficiency and robustness. DF predicts the strength of internal discontinuities for the next time step based on current cell-boundary characteristics, eliminating the need for stencil selection and improving computational efficiency.
(3). Adaptive Relaxation Using DF:
Traditional relaxation methods employ a fixed relaxation factor, which limits both convergence rate and robustness in the presence of strong discontinuities. Here, we utilize the DF factor for adaptive relaxation, dynamically adjusting the relaxation process according to the intensity of discontinuities. This enables optimal acceleration in smooth regions, while reverting to explicit time-stepping in areas with pronounced discontinuities, thereby maintaining stability and efficiency.
(4). Full GPU-Based Computation:
All computational processes—including spatial discretization, flux calculation, and multigrid iteration—are executed entirely on the GPU, dramatically accelerating solution speed. By leveraging both software and hardware acceleration, the CGKS achieves unprecedented solving speed for steady-state problems.
The proposed GMG-CGKS has been tested on inviscid and viscous flows across a range of Mach numbers on hybrid unstructured meshes. The results demonstrate significant improvements in accuracy, robustness, and computational efficiency compared to previous methods.

The remainder of this paper is organized as follows.
Section 2 introduces the three-dimensional BGK equation, the finite volume framework, and the construction of CGKS on hybrid unstructured meshes.
Section 3 details the third-order linear spatial reconstruction and the nonlinear limiting procedure.
Section 4 describes the construction of multi-color groups on unstructured meshes for the multi-color LU-SGS scheme, as well as its GPU implementation.
Section 5 presents the geometric multigrid iteration coupled with the discontinuity feedback (DF) mechanism.
Section 6 provides numerical examples, including both inviscid and viscous flow computations.
Finally, Section 7 concludes the paper.

\section{Gas kinetic scheme in finite volume framework}
\subsection{3-D BGK equation}
The Boltzmann equation \cite{cercignani1988boltzmann} describes the evolution of molecules at the mesoscopic scale. It indicates that each particle will be transported at a constant velocity or encounter a two-body collision. The BGK \cite{BGK} model simplifies the Boltzmann equation by replacing the full collision term with a relaxation model.
The 3-D gas-kinetic BGK equation \cite{BGK} is
\begin{equation}\label{bgk}
	f_t+\textbf{u}\cdot\nabla f=\frac{g-f}{\tau},
\end{equation}
where $f=f(\textbf{x},t,\textbf{u},\xi)$ is the gas distribution function, which is a function of space $\textbf{x}$, time $t$, phase space velocity $\textbf{u}$, and internal variable $\xi$.
$g$ is the equilibrium state
and $\tau$ is the collision time, which means an averaged time interval between two sequential collisions.  $g$ is expressed as a Maxwellian distribution function.
Meanwhile, the collision term on the right-hand side (RHS) of Eq.~\eqref{bgk} should satisfy the compatibility condition
\begin{equation*}\label{compatibility}
	\int \frac{g-f}{\tau} \pmb{\psi} \text{d}\Xi=0,
\end{equation*}
where $\pmb{\psi}=(1,\textbf{u},\displaystyle \frac{1}{2}(\textbf{u}^2+\xi^2))^T$,
$\text{d}\Xi=\text{d}u_1\text{d}u_2\text{d}u_3\text{d}\xi_1...\text{d}\xi_{K}$,
$K$ is the number of internal degrees of freedom, i.e.
$K=(5-3\gamma)/(\gamma-1)$ in the 3-D case, and $\gamma$
is the specific heat ratio. The details of the BGK equation can be found in \cite{xu2014directchapter2}.

In the continuous flow regime, distribution function $f$ can be taken as a small-scale expansion of Maxwellian distribution. Based on the Chapman-Enskog expansion \cite{CE-expansion}, the gas distribution function can be expressed as \cite{xu2014directchapter2},
\begin{align*}
	f=g-\tau D_{\textbf{u}}g+\tau D_{\textbf{u}}(\tau
	D_{\textbf{u}})g-\tau D_{\textbf{u}}[\tau D_{\textbf{u}}(\tau
	D_{\textbf{u}})g]+...,
\end{align*}
where $D_{\textbf{u}}={\partial}/{\partial t}+\textbf{u}\cdot \nabla$.
Through zeroth-order truncation $f=g$, the Euler equation can be obtained. The Navier-Stokes (NS) equations,
\begin{equation*}\label{ns-conservation}
	\begin{split}
		\textbf{W}_t+ \nabla \cdot \textbf{F}(\textbf{W},\nabla \textbf{W} )=0,
	\end{split}
\end{equation*}
 can be obtained by taking first-order truncation, i.e.,
\begin{align} \label{ce-ns}
	f=g-\tau (\textbf{u} \cdot \nabla g + g_t),
\end{align}
with $\tau = \mu / p$ and $Pr=1$.

Benefiting from the time-accurate gas distribution function, a time evolution solution at the cell interface is provided by the gas kinetic solver, which is distinguishable from the Riemann solver with a constant solution \cite{yang2022comparison}. This is a crucial point to construct the compact high-order gas kinetic scheme.
\begin{align}\label{point}
	\textbf{W}(\textbf{x},t)=\int \pmb{\psi} f(\textbf{x},t,\textbf{u},\xi)\text{d}\Xi,
\end{align}
and the flux at the cell interface can also be obtained
\begin{equation}\label{f-to-flux}
	\textbf{F}(\textbf{x},t)=
	\int \textbf{u} \pmb{\psi} f(\textbf{x},t,\textbf{u},\xi)\text{d}\Xi.
\end{equation}

\subsection{Finite volume framework}
The boundary of a three-dimensional arbitrary polyhedral cell $\Omega_i$ can be expressed as
\begin{equation*}
	\partial \Omega_i=\bigcup_{p=1}^{N_f}\Gamma_{ip},
\end{equation*}
where $N_f$ is the number of cell interfaces for cell $\Omega_i$.
$N_f=4$ for tetrahedron, $N_f=5$ for prism and pyramid, $N_f=6$ for hexahedron.
The semi-discretized form of the finite volume method for conservation laws can be written as
\begin{equation}\label{semidiscrete}
	\frac{\text{d} \textbf{W}_{i}}{\text{d}t}=\mathcal{L}(\textbf{W}_i)=-\frac{1}{\left| \Omega_i \right|} \sum_{p=1}^{N_f} \int_{\Gamma_{ip}}
	\textbf{F}(\textbf{W}(\textbf{x},t))\cdot\textbf{n}_p \text{d}s,
\end{equation}
with
\begin{equation*}\label{f-to-flux-in-normal-direction}
	\textbf{F}(\textbf{W}(\textbf{x},t))\cdot \textbf{n}_p=\int\pmb{\psi}  f(\textbf{x},t,\textbf{u},\xi) \textbf{u}\cdot \textbf{n}_p \text{d}\Xi,
\end{equation*}
and for ease of subsequent exposition, the residual form of the equation is presented here as Eq.~\eqref{high-order-flux-res}
\begin{equation}\label{high-order-flux-res}
\frac{\text{d} \textbf{W}_{i}}{\text{d}t}=-\textbf{R}^{n}=-\frac{1}{\left| \Omega_i \right|\Delta t} \sum_{p=1}^{N_f} \int_{\Gamma_{ip}}\int_{0}^{\Delta t}
	\textbf{F}(\textbf{W}(\textbf{x},t))\cdot\textbf{n}_p \text{d}t\text{d}s,
\end{equation}
where $\textbf{W}_{i}$ is the cell averaged values over cell $\Omega_i$, $\left|
\Omega_i \right|$ is the volume of $\Omega_i$, $\textbf{F}$ is the interface fluxes, and $\textbf{n}_p=(n_1,n_2,n_3)^T$ is the unit vector representing the outer normal direction of $\Gamma_{ip}$.
Through the iso-parametric transformation,
the Gaussian quadrature points can be determined and $\textbf{F}_{ip}(t)$ can be approximated by the numerical quadrature
\begin{equation*}\label{fv-3d-general-quadrature}
	\int_{\Gamma_{ip}}
	\textbf{F}(\textbf{W}(\textbf{x},t))\cdot\textbf{n}_p \text{d}s =  S_{i,p} \sum_{k=1}^{M} \omega_k
	\textbf{F}(\textbf{x}_{p,k},t)\cdot\textbf{n}_p,
\end{equation*}
where $S_{i,p}$ is the area of $\Gamma_{ip}$. Through the iso-parametric transformation,
in the current study, the linear element is considered.
To meet the requirement of a third-order spatial accuracy,
three Gaussian points are used for a triangular face, and four Gaussian points are used for a quadrilateral face.
In the computation, the fluxes are obtained under the local coordinates.
The details are in \cite{ji2021gradient,pan2020high,JI2024112590}.

\subsection{Gas kinetic solver}
To obtain the numerical flux at the cell interface, the integration solution based on the BGK equation is used
\begin{equation}\label{integral1}
	f(\textbf{x},t,\textbf{u},\xi)=\frac{1}{\tau}\int_0^t g(\textbf{x}',t',\textbf{u},\xi)e^{-(t-t')/\tau}\text{d}t'
	+e^{-t/\tau}f_0(\textbf{x}-\textbf{u}t,\textbf{u},\xi),
\end{equation}
where $\textbf{x}=\textbf{x}'+\textbf{u}(t-t')$ is the particle trajectory. $f_0$ is the initial gas distribution function, $g$ is the corresponding
equilibrium state in the local space and time.
This integration solution describes the physical picture of the particle evolution. Starting with an initial state $f_0$, the particle will take free transport with a probability of $e^{-\Delta t/\tau}$. Otherwise, it will suffer a series of collisions. The effect of collisions drives the system to the local Maxwellian distribution $g$, and the particles from the equilibrium propagate along the characteristics with a surviving probability of $e^{-(t-t')/\tau}$.
The components of the numerical fluxes at the cell interface can be categorized as the upwinding free transport from $f_0$ and the central difference from the integration of the equilibrium distribution.

To construct a time-evolving gas distribution function at a cell interface,
the following notations are introduced first
\begin{align*}
	a_{x_i} \equiv  (\partial g/\partial x_i)/g=g_{x_i}/g,
	A \equiv (\partial g/\partial t)/g=g_t/g,
\end{align*}
where $g$ is the equilibrium state.  The partial derivatives $a_{x_i}$ and $A$, denoted by $s$,
have the form of
\begin{align*}
	s=s_j\psi_j =s_{1}+s_{2}u_1+s_{3}u_2+s_{4}u_3
	+s_{5}\displaystyle \frac{1}{2}(u_1^2+u_2^2+u_3^2+\xi^2).
\end{align*}
The initial gas distribution function in Eq.~\eqref{integral1} can be modeled as
\begin{equation*}
	f_0=f_0^l(\textbf{x},\textbf{u})\mathbb{H} (x_1)+f_0^r(\textbf{x},\textbf{u})(1- \mathbb{H}(x_1)),
\end{equation*}
where $\mathbb{H}(x_1)$ is the Heaviside function. Here $f_0^l$ and $f_0^r$ are the
initial gas distribution functions on the left and right sides of a cell
interface, which can be fully determined by the initially
reconstructed macroscopic variables. The first-order
Taylor expansion for the gas distribution function in space around
$\textbf{x}=\textbf{0}$ can be expressed as
\begin{align}\label{flux-3d-1}
	f_0^k(\textbf{x})=f_G^k(\textbf{0})+\frac{\partial f_G^k}{\partial x_i}(\textbf{0})x_i
	=f_G^k(\textbf{0})+\frac{\partial f_G^k}{\partial x_1}(\textbf{0})x_1
	+\frac{\partial f_G^k}{\partial x_2}(\textbf{0})x_2
	+\frac{\partial f_G^k}{\partial x_3}(\textbf{0})x_3,
\end{align}
for $k=l,r$.
According to Eq.~\eqref{ce-ns}, $f_{G}^k$ has the form
\begin{align}\label{flux-3d-2}
	f_{G}^k(\textbf{0})=g^k(\textbf{0})-\tau(u_ig_{x_i}^{k}(\textbf{0})+g_t^k(\textbf{0})),
\end{align}
where $g^k$ is the equilibrium state with the form of a Maxwell distribution.
$g^k$ can be fully determined from the
reconstructed macroscopic variables $\textbf{W}
^l, \textbf{W}
^r$ at the left and right sides of a cell interface
\begin{align}\label{get-glr}
	\int\pmb{\psi} g^{l}\text{d}\Xi=\textbf{W}
	^l,\int\pmb{\psi} g^{r}\text{d}\Xi=\textbf{W}
	^r.
\end{align}
Substituting Eq.~\eqref{flux-3d-1} and Eq.~\eqref{flux-3d-2} into Eq.~\eqref{integral1},
the kinetic part of the integral solution can be written as
\begin{equation}\label{dis1}
	\begin{aligned}
		e^{-t/\tau}f_0^k(-\textbf{u}t,\textbf{u},\xi)
		=e^{-t/\tau}g^k[1-\tau(a_{x_i}^{k}u_i+A^k)-ta^{k}_{x_i}u_i],
	\end{aligned}
\end{equation}
where the coefficients $a_{x_1}^{k},...,A^k, k=l,r$ are defined according
to the expansion of $g^{k}$.
After determining the kinetic part
$f_0$, the equilibrium state $g$ in the integral solution
Eq.~\eqref{integral1} can be expanded in space and time as follows
\begin{align}\label{equli}
	g(\textbf{x},t)= g^{c}(\textbf{0},0)+\frac{\partial  g^{c}}{\partial x_i}(\textbf{0},0)x_i+\frac{\partial  g^{c}}{\partial t}(\textbf{0},0)t,
\end{align}
where $ g^{c}$ is the Maxwellian equilibrium state at an interface.
Similarly, $\textbf{W}^c$ are the macroscopic flow variables for the determination of the
equilibrium state $ g^{c}$
\begin{align}\label{compatibility2}
	\int\pmb{\psi} g^{c}\text{d}\Xi=
	\int_{u>0}\pmb{\psi} g^{l}\text{d}\Xi+
	\int_{u<0}\pmb{\psi} g^{r}\text{d}\Xi=\textbf{W}^c.
\end{align}
Substituting Eq.~\eqref{equli} into Eq.~\eqref{integral1}, the collision part in the integral solution
can be written as
\begin{equation}\label{dis2}
	\begin{aligned}
		\frac{1}{\tau}\int_0^t
		g&(\textbf{x}',t',\textbf{u},\xi)e^{-(t-t')/\tau}\text{d}t'
		=C_1 g^{c}+C_2 a_{x_i}^{c} u_i g^{c} +C_3 A^{c} g^{c} ,
	\end{aligned}
\end{equation}
where the coefficients
$a_{x_i}^{c},A^{c}$ are
defined from the expansion of the equilibrium state $ g^{c}$. The
coefficients $C_m, m=1,2,3$ in Eq.~\eqref{dis2}
are given by
\begin{align*}
	C_1=1-&e^{-t/\tau}, C_2=(t+\tau)e^{-t/\tau}-\tau, C_3=t-\tau+\tau e^{-t/\tau}.
\end{align*}
The coefficients in Eq.~\eqref{dis1} and Eq.~\eqref{dis2}
can be determined by the spatial derivatives of macroscopic flow
variables and the compatibility condition as follows
\begin{align}\label{co}
	&\langle a_{x_1}\rangle =\frac{\partial \textbf{W} }{\partial x_1}=\textbf{W}_{x_1},
	\langle a_{x_2}\rangle =\frac{\partial \textbf{W} }{\partial x_2}=\textbf{W}_{x_2},
	\langle a_{x_3}\rangle =\frac{\partial \textbf{W} }{\partial x_3}=\textbf{W}_{x_3},\nonumber\\
	&\langle A+a_{x_1}u_1+a_{x_2}u_2+a_{x_3}u_3\rangle=0,
\end{align}
where $\left\langle ... \right\rangle$ are the moments of a gas distribution function defined by
\begin{align}\label{co-moment}
	\langle (...) \rangle  = \int \pmb{\psi} (...) g \text{d} \Xi .
\end{align}
The specific details of constructing the second-order flux on the interfaces and the formula of numerical dissipation can be found in \cite{ji2021compact}

\subsection{Direct evolution of the cell averaged slopes} \label{slope-section}

The gas-kinetic solver provides the time evolution solution at a cell interface, distinguished from the Riemann solvers, which have a constant solution.
By recalling Eq.~\eqref{point}, the conservative variables at the Gaussian point  $\textbf{x}_{p,k}$ can be updated through the moments $\pmb{\psi}$
of the gas distribution function,
\begin{equation*}\label{point-interface}
	\begin{aligned}
		\textbf{W}_{p,k}(t^{n+1})=\int \pmb{\psi} f^n(\textbf{x}_{p,k},t^{n+1},\textbf{u},\xi) \text{d}\Xi,~ k=1,...,M.
	\end{aligned}
\end{equation*}

Then, the cell-averaged slopes within each element at $t^{n+1}$ can be evaluated based on the divergence theorem,

\begin{equation*}\label{gauss-formula}
	\begin{aligned}
		\nabla \overline{W}^{n+1} \left| \Omega \right|
		&=\int_{\Omega} \nabla \overline{W}(t^{n+1}) \text{d}V
		=\int_{\partial \Omega} \overline{W}(t^{n+1}) \textbf{n} \text{d}S
		= \sum_{p=1}^{N_f}\sum_{k=1}^{M_p} \omega_{p,k} W^{n+1}_{p,k} \textbf{n}_{p,k} \Delta S_p,
	\end{aligned}
\end{equation*}
where $\textbf{n}_{p,k}=((n_{1})_{p,k},(n_{2})_{p,k},(n_{3})_{p,k})$ is the outer unit normal direction at each Gaussian point $\textbf{x}_{p,k}$.

\section{Spatial reconstruction}

This section presents the 3rd-order compact reconstruction with cell-averaged values and cell-averaged slopes. To maintain both shock-capturing ability and robustness, WENO weights and DF are used \cite{ji2021gradient}. Further improvement has been made to make reconstruction simple and more robust \cite{zhang2023high}. Only one large stencil and one sub-stencil are involved in the WENO procedure.

Firstly, we will introduce the HWENO method  \cite{ji2020hweno}  for the large stencil.
Secondly, DF will be introduced.

\subsection{Linear 3rd-order compact reconstruction for large stencil}
Firstly, a linear reconstruction approach is presented. To achieve a third-order accuracy in space, a quadratic polynomial $p^2$ is constructed as follows.

\begin{equation}
p^2(\mathbf{x})=\overline{Q}_0+\sum\limits_{|k|=1}^2a_k[(x-x_0)^{k_1} (y-y_0)^{k_2} (z-z_0)^{k_3}-\frac{1}{\left|\Omega_0\right|} \iiint_{\Omega_0} (x-x_0)^{k_1} (y-y_0)^{k_2} (z-z_0)^{k_3}] \mathrm{~d} V,
\end{equation}
where $k=\left(k_1, k_2, k_3\right)$ is the multi-index, $k|=k_1+k_2+k_3$, $(x_0, y_0, z_0)$ is the geometric center coordinate. The $p^2$ on $\Omega_0$  is constructed on the compact stencil $S$ including $\Omega_0$ and all its von Neumann neighbors $\Omega_m$ ($m=1,\cdots,N_f$, $N_f$ is the number of $\Omega_0$'s faces). The cell averages $\overline{Q}$ on $\Omega_0$ and $\Omega_m$ together with cell-averaged slopes $\overline{Q}_x,\overline{Q}_y $ and $\overline{Q}_z$ on $\Omega_m$ are used to obtain $p^2$.

The polynomial $p^2$ naturally satisfies cell averages over $\Omega_0$
\begin{equation}\label{self-constrain}
	\iiint_{\Omega_0} p^2 \text{d}V = \overline{Q}_0|\Omega_0|,
\end{equation}
Meanwhile, the $p^2$ is also required to exactly satisfy cell averages over the target cell's neighbors $\Omega_m$
\begin{equation}\label{neighbor-constrain}
 \iiint_{\Omega_m} p^2 \text{d}V = \overline{Q}_m|\Omega_m|\\
\end{equation}
Then, the following conditions are satisfied in a least-square sense
\begin{equation*}
	\begin{aligned}
		\iiint_{\Omega_m}\frac{\partial}{\partial x} p^2 \text{d}V = \left(\overline{Q}_x\right)_m|\Omega_m|\\
		\iiint_{\Omega_m}\frac{\partial}{\partial y} p^2 \text{d}V = \left(\overline{Q}_y\right)_m|\Omega_m|\\
		\iiint_{\Omega_m}\frac{\partial}{\partial z} p^2 \text{d}V = \left(\overline{Q}_z\right)_m|\Omega_m|.\\
	\end{aligned}
\end{equation*}
The constrained least-square method is used to meet the above requirements.

\subsection{Green-Gauss reconstruction for the sub stencil}
The classical Green-Gauss reconstruction \cite{shima2013green} with only cell-averaged values is adopted to provide the linear polynomial $p^1$ for the sub stencil.
\begin{equation*}
	p^1=\overline{Q}+\boldsymbol{x}\cdot \sum_{m=1}^{N_f}\frac{\overline{Q}_m+\overline{Q}_0}{2|\Omega_0|}S_m\boldsymbol{n} _m,
\end{equation*}
where $S_m$ is the area of the cell's surface and $\boldsymbol{n}_m$ is the surface's normal vector. In most cases, Green-Gaussian reconstruction has only first-order precision.
\subsection{Discontinuity Feedback}
The DF was first proposed in \cite{ji2021gradient}. Here, several improvements have been made in \cite{zhang2023high}: there is no $\epsilon$ in the improved expression of DF; the difference of Mach number is added to improve the robustness under strong rarefaction waves. Denote $\alpha_i \in[0,1]$ as DF at targeted cell $\Omega_i$
\begin{equation*}
	\alpha_i = \prod_{p=1}^m\prod_{k=0}^{M_p}\alpha_{p,k},
\end{equation*}
where $\alpha_{p,k}$ is the CF obtained by the $k$th Gaussian point at the interface $p$ around cell $\Omega_i$, which can be calculated by
\begin{equation*}
	\begin{aligned}
		&\alpha_{p,k}=\frac{1}{1+D^2}, \\
		&D=\frac{|p^l-p^r|}{p^l} +\frac{|p^l-p^r|}{p^r}+(\text{Ma}^{l}_n-\text{Ma}^{r}_n)^2+(\text{Ma}^{l}_t-\text{Ma}^{r}_t)^2,
	\end{aligned}
\end{equation*}
where $p$ is pressure, $\text{Ma}_n$ and $\text{Ma}_t$ are the Mach numbers defined by normal and tangential velocity, and superscripts $l,r$ denote the left and right values of the Gaussian points.

Then, the updated slope is modified by
\begin{equation*}
	\widetilde{\overline{\nabla \boldsymbol{W}}}_i^{n+1} = \alpha_i\overline{\nabla \boldsymbol{W}}_i^{n+1},
\end{equation*}
and the Green-Gauss reconstruction is modified as
\begin{equation*}
	p^1=\overline{Q}+\alpha \boldsymbol{x}\cdot \sum_{m=1}^{N_f}\frac{\overline{Q}_m+\overline{Q}_0}{2\left| \Omega_0 \right|}S_m\boldsymbol{n} _m.
\end{equation*}
The nonlinear reconstruction procedure can be found in our previous work \cite{ji2021compact,zhang2023high}.

\section{Matrix-free form of multi-color LU-SGS on GPU}
LU-SGS effectively eliminates high-frequency waves, and is widely used as the smoother of the multi-grid solvers.
However, as previously discussed, GPU data dependency hinders the parallel implementation of the standard LU-SGS algorithm.

\subsection{Multi-color painting processor}
A particular strategy is developed to overcome the data dependency problem on GPU \cite{bocharov2020implicit, zhang2018gpu}.
The main idea of this method is to partition the grid cells into several independent groups, where all cells within each group share the same color designation.
The fundamental principle of this coloring algorithm is that adjacent cells are assigned different colors.
In the computer program, we represent different colors using integer numbers. For instance, the color red is denoted by index 1, while the color black is represented by index 2.
The procedure to paint all the computational points is described in Algorithm\ref{coloring-algorithm}.

\begin{algorithm}\label{coloring-algorithm}
	\caption{ The procedure of coloring method.
	}
	
	\label{coloring-algorithm}
	\KwIn {The original mesh cells $\left\{C_i \mid i \in \Omega\right\}$ and a start cell $v_0$.}
	\KwOut {The array of cell colors colorArray $(:)$ and total number of colors $N_{\text {color }}$.}
	initialize colorArray $(:)=0$\;
	choose a cell $v_0$ as the start point of the traversal, and set $\operatorname{color}\left(v_0\right)=1$\;
	\Repeat{all cells are painted}
	{
		\For{each colored cell $\{v \mid \operatorname{color}(v)>0\}$}
		{
			\For{each uncolored cell $\left\{w \mid w \in C_v\right\}$}
			{
				paint $\operatorname{color}(w)=\min \left\{k>0 \mid k \neq \operatorname{color}(j), \forall j \in C_w\right\}$\;
			}
		}
	}

\end{algorithm}	

\begin{figure}
	\centering
	\includegraphics[width=0.35\linewidth]{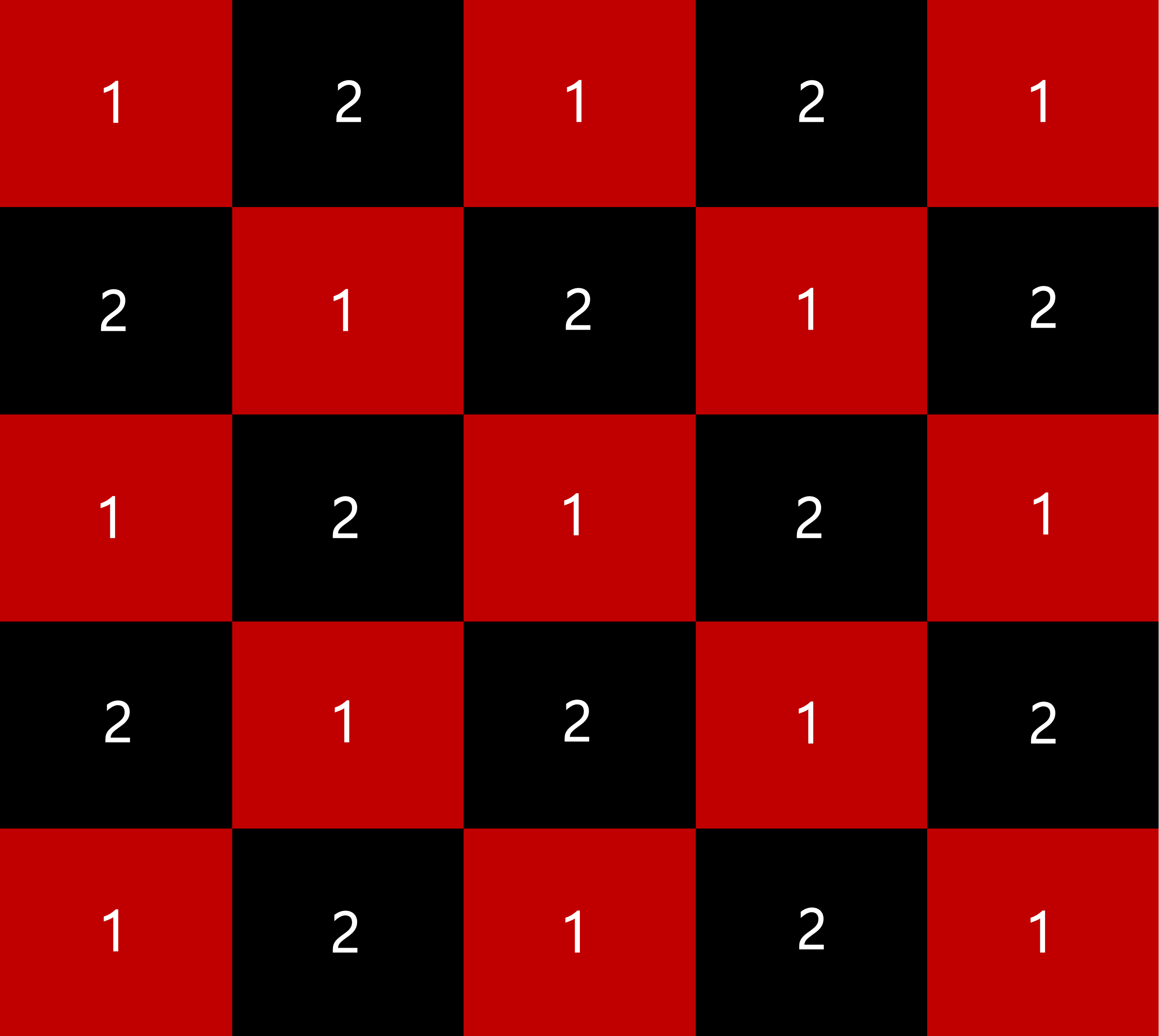}
	\includegraphics[width=0.35\linewidth]{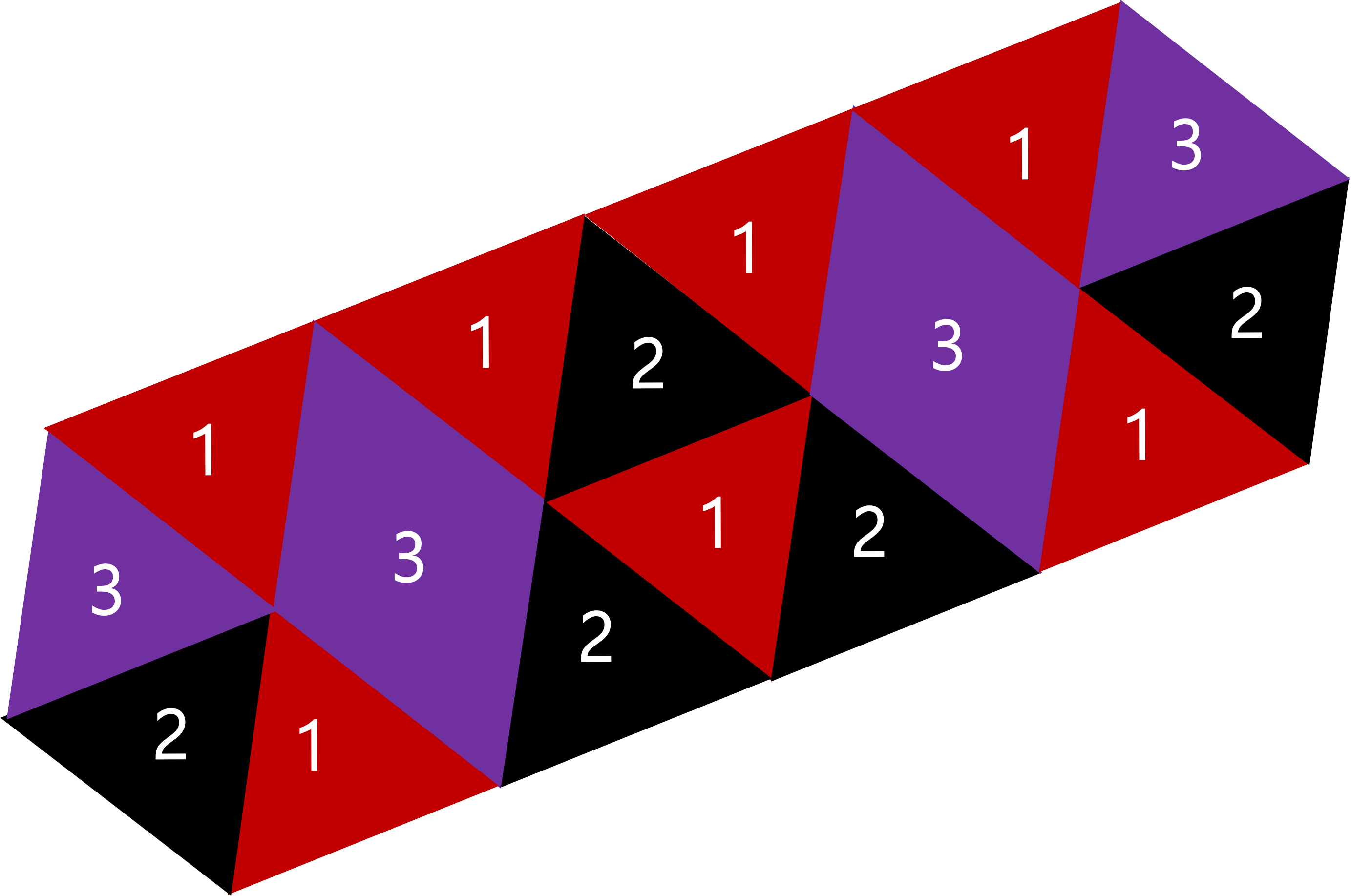}
	\caption{Mesh color distribution. Left: Regular mesh. Right: Irregular mesh.}
	\label{color-mesh}
\end{figure}

The painting procedure outlined in Algorithm\ref{coloring-algorithm} begins by selecting a starting cell $v_0$ within the computational domain. 
Subsequently, the color graph corresponding to this starting cell is established. 
Random start cells were experimented with to assess the impact of varying start cells on computational efficiency, revealing minimal influence. For consistency and convenience, the initial cell in the global array is consistently designated as the start node in this study.
As examples, the color graphs generated for regularly structured and irregularly unstructured mesh are depicted in Fig.~\ref{color-mesh}.

The left part of Fig.~\ref{color-mesh} displays a simple graph with two colors generated using Algorithm\ref{coloring-algorithm} for a structured mesh. 
Each red cell's implicit computation (color index 1) depends solely on itself and the adjacent black cells (color index 2) within its local color group. 
Similarly, the implicit computation of each black cell relies only on itself and the neighboring red cells. 
Consequently, algebraic operations on cells of the same color are independent, facilitating straightforward parallelization. 
The approach can be extended to irregularly unstructured meshes. However, a larger number of colors may be required to assign to cells due to the intricate distribution, as illustrated in the right part of Fig.~\ref{color-mesh}.
\subsection{Robust multi-color matrix-free form LU-SGS on GPU}
This section will introduce the traditional serial LU-SGS algorithm to facilitate the subsequent discussion of the GPU-based LU-SGS algorithm.
The residual form of conservation law can be written as Eq.~\eqref{resform}
\begin{equation}\label{resform}
	\frac{d \mathbf{W}_i}{d t}=-\sum_{j \in N_i} \mathbf{F}_{i j}\mathbf{S}_{i j}=\mathbf{-R}^{n},
\end{equation}
where $\mathbf{F}_{i j}$ is the total flux per unit area per unit time, $\mathbf{W}_{i}$ is the cell average conservative flow variable, and $=\mathbf{-R}^{n}$ is the residual using current step's variable.
Then, the governing equations can be rewritten in delta form as
\begin{equation}
	\frac{V_i}{\Delta t} \Delta \mathbf{W}_i^{n+1}+\sum_{j \in N(i)} S_{i j} \Delta \mathbf{F}_{i j}^{n+1}=-\sum_{j \in N(i)} S_{i j} \mathbf{F}_{i j}^n,
\end{equation}
where $\Delta \mathbf{W}_i^{n+1}=\mathbf{W}_i^{n+1}-\mathbf{W}_i^n$ and $\Delta \mathbf{F}_{i j}^{n+1}=\mathbf{F}_{i j}^{n+1}-\mathbf{F}_{i j}^n$.

To save GPU memory usage and computational cost, we approximate $\Delta \mathbf{F}_{i j}^{n+1}$ using flux splitting based on the Euler equations.
\begin{equation}
	\Delta \mathbf{F}_{i j}^{n+1}=\frac{1}{2}\left[\Delta \mathbf{T}_i^{n+1}+\Delta \mathbf{T}_j^{n+1}+r_{i j}\left(\Delta \mathbf{W}_i^{n+1}-\Delta \mathbf{W}_j^{n+1}\right)\right],
\end{equation}
where $\mathbf{T}$ is the Euler flux, $r_{i j}$ is the spectral radius and satisfies $r_{i j} \geq \Lambda_{i j}=\left|\vec{U}_{i j} \cdot \vec{n}_{i j}\right|+a_s$ and $a_s$ is the sound speed at the interface $ij$.


Since the relationship $\sum_{i \in N(i)} S_{i j} \Delta T_i^{n+1}=0$ is always satisfied for any closed finite volume $i$ in a steady state solution, and the specific delta form of conservation law is
\begin{equation*}
	\left(\frac{V_i}{\Delta t}+\frac{1}{2} \sum_{j \in N(i)} S_{i j} r_{i j}\right) \Delta \mathbf{W}_i^{n+1}+\frac{1}{2} \sum_{j \in N(i)} S_{i j}\left[\mathbf{T}\left(\mathbf{W}_j^{n+1}\right)-\mathbf{T}\left(\mathbf{W}_j^n\right)-r_{i j} \Delta \mathbf{W}_j^{n+1}\right]=-\mathbf{R}_i^{n}
\end{equation*}
Then, the above equation can solved by the traditional LU-SGS way. 
The forward step is:
\begin{equation*}\label{cpu-forward}
	\left(\frac{V_i}{\Delta t}+\frac{1}{2} \sum_{j \in N(i)} S_{i j} r_{i j}\right) \Delta \mathbf{W}_i^{\star}+\frac{1}{2} \sum_{j \in L(i)} S_{i j}\left[\mathbf{T}\left(\mathbf{W}_j^n+\Delta \mathbf{W}_j^{\star}\right)-\mathbf{T}\left(\mathbf{W}_j^n\right)-r_{i j} \Delta \mathbf{W}_j^{\star}\right]
	=-\mathbf{R}_i^n,
\end{equation*}
and the backward step is:
\begin{equation*}
	\begin{aligned}
		& \left(\frac{V_i}{\Delta t}+\frac{1}{2} \sum_{j \in N(i)} S_{i j} r_{i j}\right) \Delta \mathbf{W}_i^{n+1}+\frac{1}{2} \sum_{j \in U(i)} S_{i j}\left[\mathbf{T}\left(\mathbf{W}_j^n+\Delta \mathbf{W}_j^{n+1}\right)-\mathbf{T}\left(\mathbf{W}_j^n\right)-r_{i j} \Delta \mathbf{W}_j^{n+1}\right] \\
		& \quad=\left(\frac{V_i}{\Delta t}+\frac{1}{2} \sum_{j \in N(i)} S_{i j} r_{i j}\right) \Delta \mathbf{W}_i^{\star},
	\end{aligned},
\end{equation*}
where $L(i)$ and $U(i)$ are subsets of $N(i)$, and $L(i)$ is the neighboring cells of cell $i$ occupying the lower triangular area of this matrix, and $U(i)$ is the ones in the upper triangular area.

The traditional LU-SGS algorithm sweeps all computational cells sequentially, which is unsuitable for parallel computing. 
We have introduced a modification using the multi-coloring strategy to address this limitation. This new algorithm traverses all data cells group by group, starting from the first color and progressing to the last color in the forward updating step. 
Subsequently, it moves across the cells from the last color to the first color in the backward iteration.

The specific formula of the forward step is:
\begin{equation}\label{gpu-forward}
	\left(\frac{V_i}{\Delta t}+\frac{1}{2} \sum_{j \in N(i)} S_{i j} r_{i j}\right) \Delta \mathbf{W}_i^{\star}+\frac{1}{2} \sum_{\substack{j \in N(i) \\ Cl(j)<Cl(i)}} S_{i j}\left[\mathbf{T}\left(\mathbf{W}_j^n+\Delta \mathbf{W}_j^{\star}\right)-\mathbf{T}\left(\mathbf{W}_j^n\right)-r_{i j} \Delta \mathbf{W}_j^{\star}\right]
	=-\mathbf{R}_i^n,
\end{equation}
and the specific formula of the backward step is:
\begin{equation}\label{gpu-backward}
	\begin{aligned}
		& \left(\frac{V_i}{\Delta t}+\frac{1}{2} \sum_{j \in N(i)} S_{i j} r_{i j}\right) \Delta \mathbf{W}_i^{n+1}+\frac{1}{2} \sum_{\substack{j \in N(i) \\ Cl(j)>Cl(i)}} S_{i j}\left[\mathbf{T}\left(\mathbf{W}_j^n+\Delta \mathbf{W}_j^{n+1}\right)-\mathbf{T}\left(\mathbf{W}_j^n\right)-r_{i j} \Delta \mathbf{W}_j^{n+1}\right] \\
		& \quad=\left(\frac{V_i}{\Delta t}+\frac{1}{2} \sum_{j \in N(i)} S_{i j} r_{i j}\right) \Delta \mathbf{W}_i^{\star},
	\end{aligned}
\end{equation}
where, $Cl$ is the color array.

\subsection{DF based relaxation scheme for time marching}
Meanwhile, due to discontinuities within the elements, we should transition from implicit to explicit scheme based on the strength of the discontinuities within the elements to enhance robustness.
As discussed before, DF is an adaptive discontinuity feedback factor.
It naturally serves as a basis for transitioning from implicit to explicit schemes.
The explicit time marching can be written as:

\begin{equation}\label{DF-relaxation-exp}
	\frac{V_i}{\Delta t} \Delta \mathbf{W}_i^{n+1,explicit}=-\sum_{j \in N(i)} S_{i j} \mathbf{F}_{i j}^n.
\end{equation}

The implicit time marching can be written as:

\begin{equation}\label{DF-relaxation-implicit}
	\frac{V_i}{\Delta t_i} \Delta \mathbf{W}_i^{n+1,implicit}=-\sum_{j \in N(i)} S_{i j} \Delta \mathbf{F}_{i j}^{n+1}-\sum_{j \in N(i)} S_{i j} \mathbf{F}_{i j}^n.
\end{equation}

Constructing a convex combination of the explicit time marching scheme Eq.\eqref{DF-relaxation-exp} and the implicit time marching scheme Eq.\eqref{DF-relaxation-implicit}, a new hybrid time marching scheme can be obtained by multiplying DF on Eq.\eqref{DF-relaxation-implicit} and multiplying $(1-\text{DF})$ on Eq.\eqref{DF-relaxation-exp}:

\begin{equation}\label{DF-relaxation-hybrid}
	\left[\alpha \frac{V_i}{\Delta t_i}+\left(1-\alpha\right)\frac{V_i}{\Delta t}\right] \Delta \mathbf{W}_i^{n+1}+\alpha\sum_{j \in N(i)} S_{i j} \Delta \mathbf{F}_{i j}^{n+1}=-\sum_{j \in N(i)} S_{i j} \mathbf{F}_{i j}^n,
\end{equation}
where $\alpha$ is the DF.

\begin{remark}
	Limitation analysis:
	
	When ${\alpha \to 0} $, Eq.\eqref{DF-relaxation-hybrid} revert to Eq.\eqref{DF-relaxation-exp}, indicating that under strong discontinuities, the hybrid scheme revert to the explicit scheme.
	
	When ${\alpha \to 1} $, Eq.\eqref{DF-relaxation-hybrid} revert to Eq.\eqref{DF-relaxation-implicit}, indicating that in the smooth region, the hybrid scheme revert to the fully implicit scheme.
	
	This ensures sufficient acceleration efficiency in smooth regions and strong robustness in discontinuous regions.
\end{remark}

\begin{remark}
	The traditional relaxation method uses a fixed relaxation factor and the specific implicit scheme formula can be written as:
	\begin{equation}
		\beta\frac{V_i}{\Delta t_i} \Delta \mathbf{W}_i^{n+1}+\sum_{j \in N(i)} S_{i j} \Delta \mathbf{F}_{i j}^{n+1}=-\sum_{j \in N(i)} S_{i j} \mathbf{F}_{i j}^n+(1-\beta)\Delta \mathbf{W}_i^{*},
	\end{equation}
	This may leads to the slow down convergence rate in the smooth region and less robustness encountering discontinuous due to the lack of flow regime adaptive time marching method.
\end{remark}

Then the multi-color LU-SGS method to solve Eq.\eqref{DF-relaxation-hybrid} can be written as a forward step and a backward step:
The specific formula of the forward step is:
\begin{equation}\label{gpu-forward-relaxation}
	\begin{aligned}
		\left[\alpha\left(\frac{V_i}{\Delta t}+\frac{1}{2} \sum_{j \in N(i)} S_{i j} r_{i j}\right) + \left(1-\alpha\right)\frac{V_i}{\Delta t}\right] \Delta \mathbf{W}_i^{\star}+ \\
		\frac{1}{2}\alpha \sum_{\substack{j \in N(i) \\ Cl(j)<Cl(i)}} S_{i j}\left[\mathbf{T}\left(\mathbf{W}_j^n+\Delta \mathbf{W}_j^{\star}\right)-\mathbf{T}\left(\mathbf{W}_j^n\right)-r_{i j} \Delta \mathbf{W}_j^{\star}\right]
		=-\mathbf{R}_i^n,
	\end{aligned}
\end{equation}

and the backward step is:
\begin{equation}\label{gpu-backward-relaxation}
	\begin{aligned}
		& \left[\alpha\left(\frac{V_i}{\Delta t}+\frac{1}{2} \sum_{j \in N(i)} S_{i j} r_{i j}\right) + \left(1-\alpha\right)\frac{V_i}{\Delta t}\right] \Delta \mathbf{W}_i^{n+1} \\
		&+ \frac{1}{2} \alpha \sum_{\substack{j \in N(i) \\ Cl(j)>Cl(i)}} S_{i j}\left[\mathbf{T}\left(\mathbf{W}_j^n+\Delta \mathbf{W}_j^{n+1}\right)-\mathbf{T}\left(\mathbf{W}_j^n\right)-r_{i j} \Delta \mathbf{W}_j^{n+1}\right] \\
		& \quad=\left(\frac{V_i}{\Delta t}+\frac{1}{2} \sum_{j \in N(i)} S_{i j} r_{i j}\right) \Delta \mathbf{W}_i^{\star},
	\end{aligned}
\end{equation}

The total procedure of forward sweep and backward sweep of multi-color LU-SGS is shown as Algorithm \ref{multi-color-LU-SGS-algorithm}.

\begin{algorithm}\label{multi-color-LU-SGS-algorithm}
	\caption{ The procedure of multi-color LU-SGS method.}
	
	\label{multi-color-LU-SGS-algorithm}
	\textbf{Forward sweep}\;
	\For{$(c=1$ to $N_{\text {color }})$}
	{
		forward sweep the cells in $c^{th}$ color group by Eq.~\eqref{gpu-forward-relaxation}\;
	}
	\textbf{Backward sweep}\;
	\For{$($$c=N_{\text {color }}$ to 1$)$}
	{
		backward sweep the cells in $c^{th}$ color group by Eq.~\eqref{gpu-backward-relaxation}\;
	}
\end{algorithm}

\section{Geometric multigrid algorithm for third-order compact gas kinetic scheme}
\subsection{Mesh coarsening algorithm}
Since the number of cells in the large-scale simulation's mesh is usually very large, it is important to design an efficient algorithm for mesh coarsening. In this paper, a fast parallel mesh agglomeration algorithm is developed.

In this section, the face removal algorithm \cite{zhang20133d} is introduced first. The boundary face and parallel interface will not be deleted to implement boundary conditions easily and efficiently transfer the information at the parallel interface. Considering a face with its left cell and right cell, its hash value is
\begin{equation} \label{hash value}
	h_f=[23(N_l+N_r)+N_lN_r]\%N_f
\end{equation}
where $h_f$ is the face's hash value, $N_l$ is the face's left cell index, $N_r$ is the face's right cell index, and $N_f$ is the total interior face number of the mesh.  For example, if the face's left cell ID is 10, the face's right cell ID is 12, and the total face number of the mesh is 100, then the hash value of the face will be 5.

With the face's hash value calculation method, interior faces will be iterated and selected to delete the collection of faces. The face's hash value will be calculated when iterating at a face. Then, whether the hash value has already been inserted in the built hash table will be determined. If the face's hash value has already been inserted in the hash table, the face will not be selected for the collection.

\begin{figure}[htp]	
	\centering
	\includegraphics[width=0.35\textwidth]
	{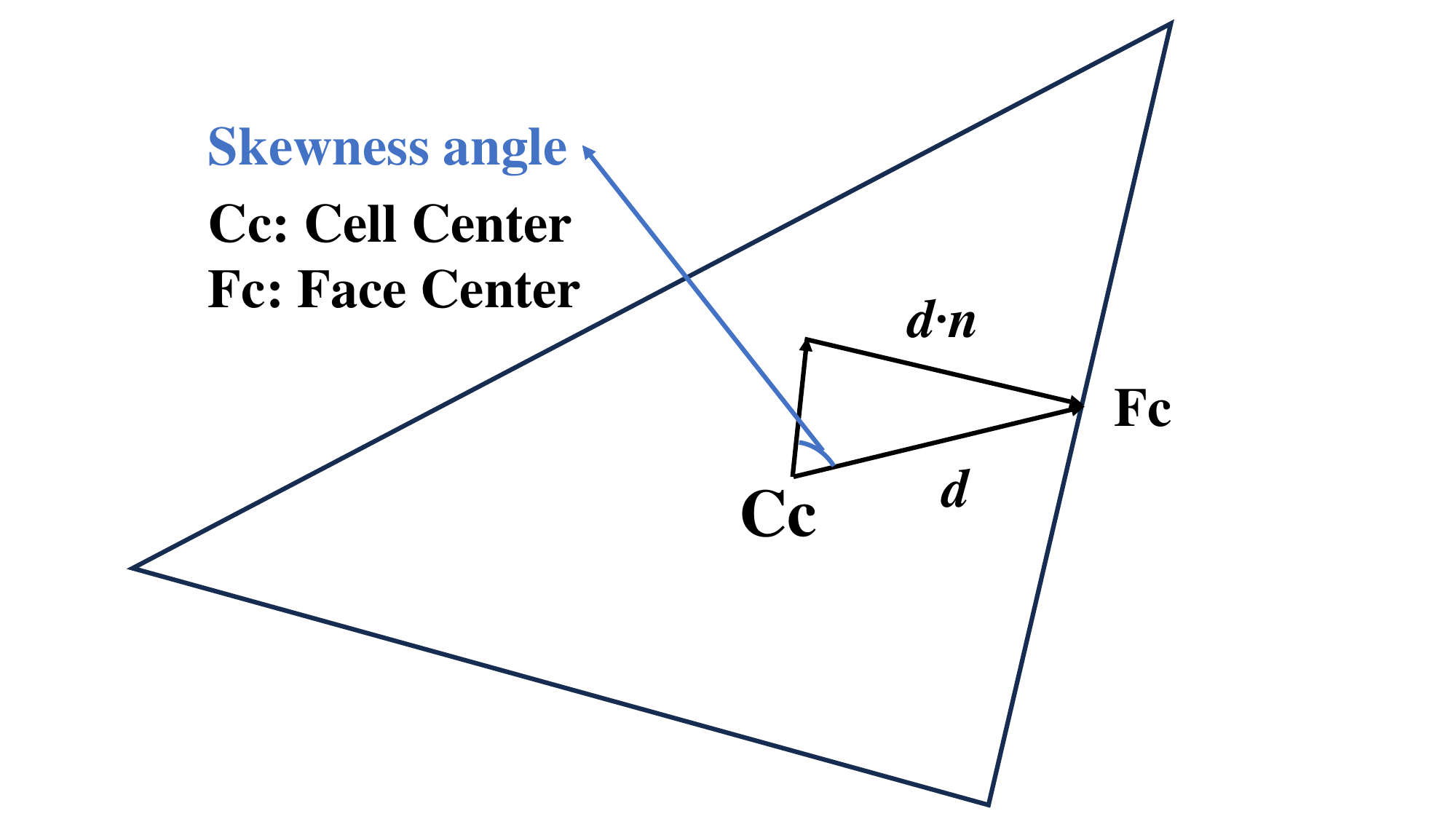}
	\caption{\label{skewsness}
		Skews angle factor .}
\end{figure}

To guarantee the coarse mesh's quality, the skewness factor will be introduced to indicate the quality of the mesh. When a face in the above collection is selected to be deleted, the virtual center of the cell, which will be constructed by the face's left and right cells, is
\begin{equation} \label{virtual center}
	C_c=\frac{V_lC_l+V_rC_r}{V_l+V_r}
\end{equation}
where $V_l$ is the volume of the left cell, $C_l$ is the center of the left cell, $V_r$ is the volume of the right cell, and $C_r$ is the center of the right cell.
The angle between the unit normal vector of the face and the vector pointed from the virtual cell center to the face center is
\begin{equation} \label{skewness factor}
	\sin \alpha=\frac{d\cdot n}{|d|}
\end{equation}
where $\alpha$ is the skewness angle representing the skewness factor, d is the distance vector pointed from the virtual cell center to the face center, and n is the face's unit normal vector shown in Fig.~\ref{skewsness}.
If the angle is smaller than a value set previously, the mesh's quality is bad, and the face will not be deleted.
To sum up, the cell merge algorithm can be summarized as Algorithm.~\ref{coarsening algorithm} shows.

\begin{algorithm}\label{coarsening algorithm}
	\caption{Cell coarsening method}
	\label{cell-merge}
	Initialization: the collection of faces to be deleted $V_d$, the collection's volume at begin is zero.\;
	\For{j=1,..,$N_f$}
	{
		
		Calculate the $face_j$ 's hash value by Eq.~\eqref{hash value}\;
		If the hash value hasn't appeared in the previous hash table, insert the face into the collection.
	}
	\For{i=1,..,$N_i$}
	{
		where $N_i$ is the element in the collection of faces to be deleted\;
		Calculate the skewness of the face and its neighbor cells by Eq.~\eqref{skewness factor}\;
		If the skewness factor is bigger than the limit, then merge the left and right cells of the face\;
	}
\end{algorithm}	

After coarsening the original mesh for the first time, the geometric variables of the coarse cell should be recalculated. The volume of the coarse cell is
\begin{equation}
	V_c=\sum_{i=0}^n V_i
\end{equation}
where $V_c$ is the volume of the coarse cell, $V_i$ is the volume of the fine cell constructing the coarse cell. The center of the cell is
\begin{equation}
	C_c=\frac{\sum_{i=0}^n V_iC_i}{V_c}
\end{equation}

\subsection{geometric multi-grid for CGKS}
This section employs a hybrid geometric multigrid and p-multigrid method in a high-order compact gas kinetic scheme. For the finest mesh, the third-order gas kinetic scheme solver is used as the pre-smoother and post-smoother \cite{mavriplis1990multigrid}. The first-order KFVS solver coupled with multi-color LU-SGS is pre-smoother and post-smoother for the other coarse meshes. For the finest mesh, the semi-discretization equation can be written as
\begin{equation}\label{smo}
	W_h^{n+1}=W_h^{n}+\frac{V_i}{\Delta t}\sum_{j} Res(W_h^n)
\end{equation}
where $W_h^{n+1}$ is the cell average conservative variable at step $n+1$, $W_h^{n}$ is the cell average conservative variable at step $n$, $Res(W_h^n)$ is the residual calculated at step $n$ using cell average conservative variable $W_h^n$, $h$ means fine level, and $j$ means the $j^{th}$ face of the current cell.

\begin{figure}[htp]
	\centering
	\includegraphics[width=1.0\textwidth]
	{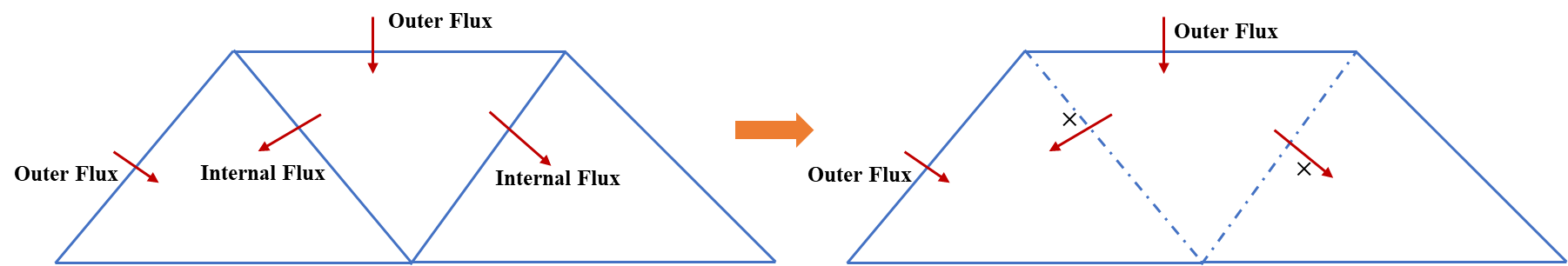}
	\caption{\label{Forcing term}
		The physical picture of the flux relationship.}
\end{figure}

Having the cell average value updated once means smoothing it once. After smoothing the cell average value, the residual used for restricting the coarse meshes residual can be updated by Eq. \eqref{smo}. The cell average conservative variable  is restricted to coarse mesh by
\begin{equation}
	W_{2h}^0=\frac{\sum_{i\in {cell2h}} V_{ifc}W_{ifc}}{V_{2h}}
\end{equation}
where $W_{2h}^0$ is the initial conservative variable in the first step of coarse mesh calculation, $V_{2h}$ is the coarse cell's volume, $W_{ifc}$ is the $i^{th}$ fine cell's conservative variable.
Residual restriction from fine mesh to coarse mesh is  guided by
\begin{equation}
	Res_{2h}^{*}=\sum_{i\in{cell2h}} Res_{ifc}^{n+1}
\end{equation}
where $Res_{2h}^{*}$ is the residual restricted from fine mesh, $Res_{ifc}^{n+1}$ is the $i^{th}$ fine cell's residual. It's important to notice that the residual restriction process's physical meaning is the coarse cell's total flux using the flux calculated from fine mesh. It is shown in Fig.~\ref{Forcing term}.

To guarantee the final result can still keep high order accuracy, a forcing term is given by
\begin{equation}
	F=Res_{2h}^{*}-Res_{2h}(W_{2h}^{0})
\end{equation}

For the coarse mesh, the macroscopic variable is updated using the forcing term by
\begin{equation}
	\frac{V_i}{\Delta t_i} \Delta \mathbf{W}_i^{n+1}+\sum_{j \in N(i)} S_{i j} \Delta \mathbf{F}_{i j}^{n+1}=-\sum_{j \in N(i)} S_{i j} \mathbf{F}_{i j}^n+F,
\end{equation}

After updating the macroscopic variable on coarse mesh, a correction process can be used to correct the fine mesh solution by
\begin{equation}\label{prolongation}
	W_h^{*}=W_h^{n+1}+\alpha(\mathbf{W},\nabla \mathbf{W})I_{2h}^h(W_{2h}^{k+1}-W_{2h}^{0})
\end{equation}

\begin{figure}[htp]
	\centering
	\includegraphics[width=0.8\textwidth]
	{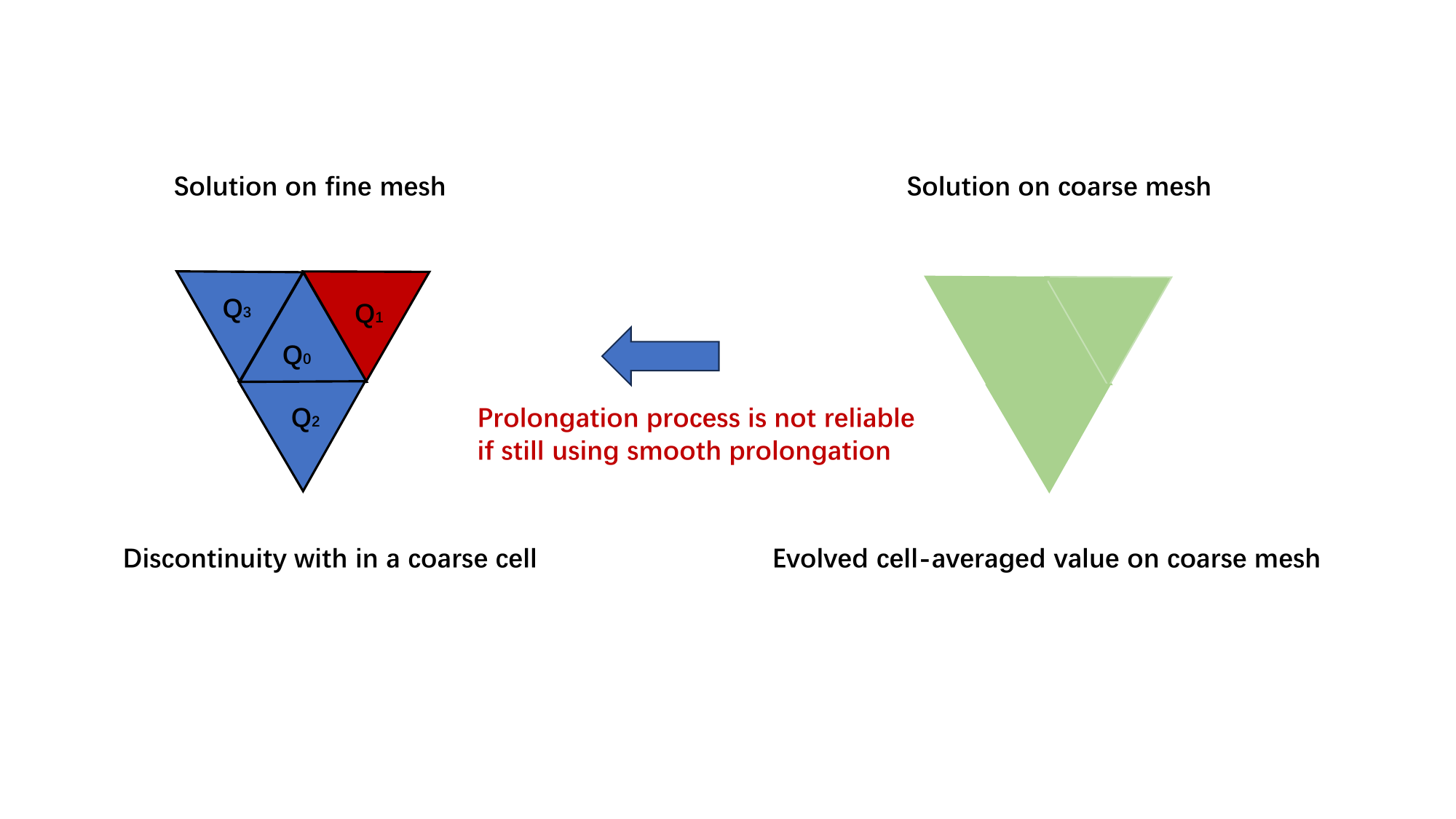}
	\caption{\label{MultiGrid-DF}
		Explain the reason why the prolongation procedure using DF.}
\end{figure}

This paper uses DF as the prolongation limiter to overcome the discontinuity within the cell.
The essence of this approach is that, due to discontinuities within the coarse grid, the solution on the coarse grid is no longer reliable. At this point, the DF must further constrain the prolongation process to limit the correction size.
Fig.~\ref{MultiGrid-DF} can help to explain the motivation.

What's more, compared with traditional upwind stencils chosen during the prolongation procedure, using DF is more compact and thus improves computational efficiency.

In this paper, the pre-smooth number is set to 1 without specifying, and the post-smooth number is set to 0.

\begin{figure}[htp]
	\centering
	\includegraphics[width=0.8\textwidth]
	{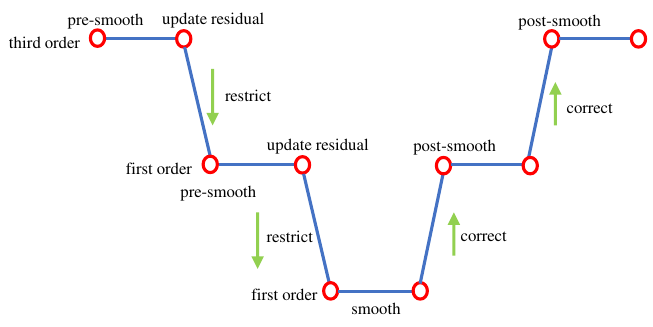}
	\caption{\label{algorithm pic}
		Hybrid method of geometric multigrid and p-multigrid.}
\end{figure}

A three-level V cycle geometric grid algorithm's picture is given by Fig.~\ref{algorithm pic}.

\subsection{A brief analysis of current geometric multi-grid CGKS}
\begin{figure}[htp]
	\centering
	\includegraphics[width=1.0\textwidth]
	{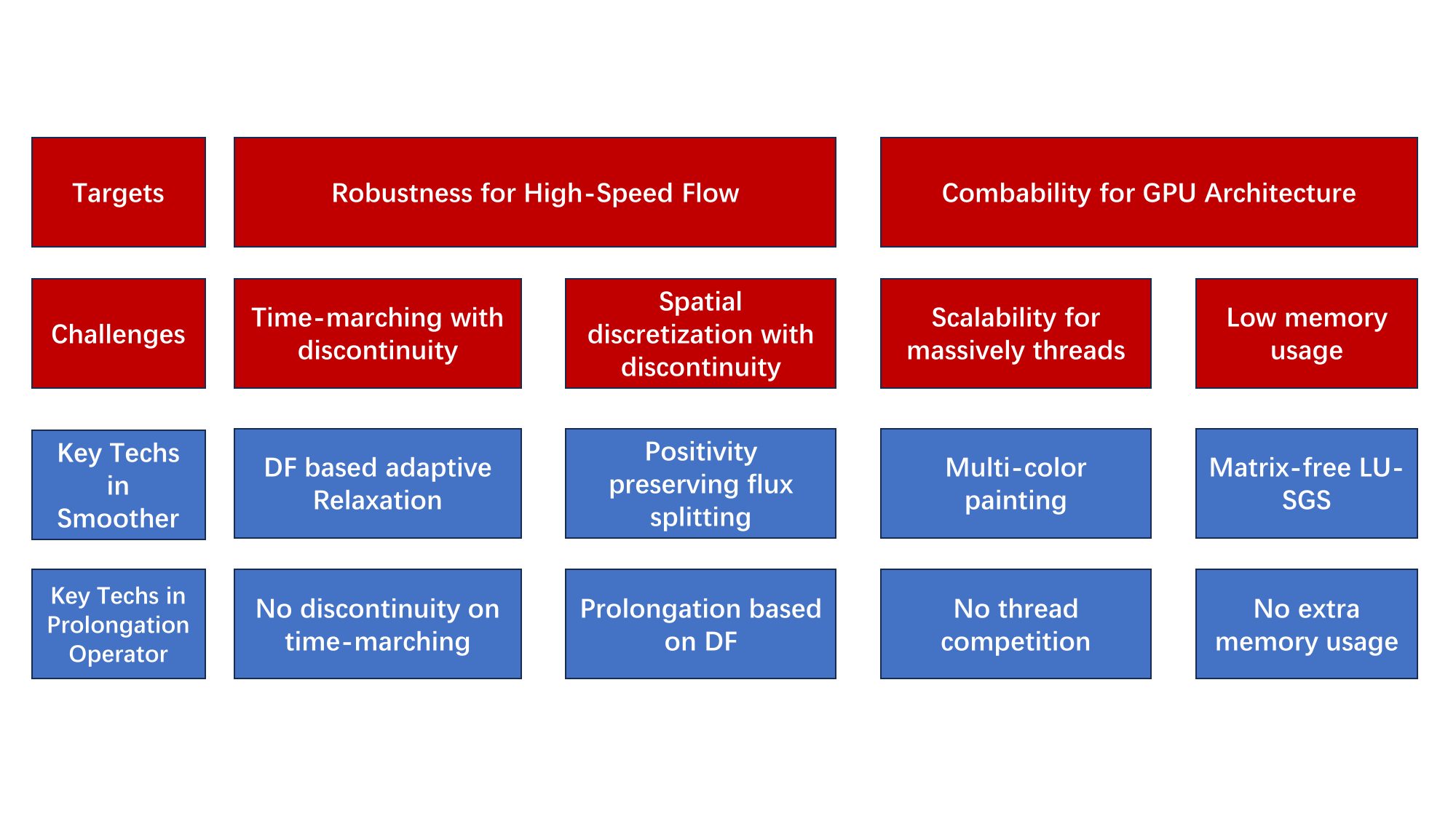}
	\caption{\label{GMG-compare}
		Comprehensive conclusion of GMG-CGKS facing discontinuity.}
\end{figure}

This paper uses DF as the prolongation limiter to overcome the discontinuity within the cell.
The essence of this approach is that, due to discontinuities within the coarse grid, the solution on the coarse grid is no longer reliable. At this point, the DF must further constrain the prolongation process to limit the correction size.
\begin{remark}
	when $\alpha \to 0$, Eq.\eqref{prolongation} becomes:
	\begin{equation}
		W_h^{*}=W_h^{n+1},
	\end{equation}
	which means under strong discontinuities, the final solution is obtained only by fine mesh time marching.
\end{remark}
Moreover, compared with traditional upwind stencils chosen during the prolongation procedure, DF-based prolongation is more compact and thus improves computational efficiency.
In multi-color LU-SGS and multi-grid processes, introducing DF ensures robustness when there are discontinuities in time and space without adding extra computational cost, as the DF has already been calculated in the original CGKS framework during the flux calculation process.
In conclusion, the properties of GMG-CGKS facing discontinuity are shown in Fig.~\ref{GMG-compare}.

\begin{figure}[htp]
	\centering
	\includegraphics[width=1.0\textwidth]
	{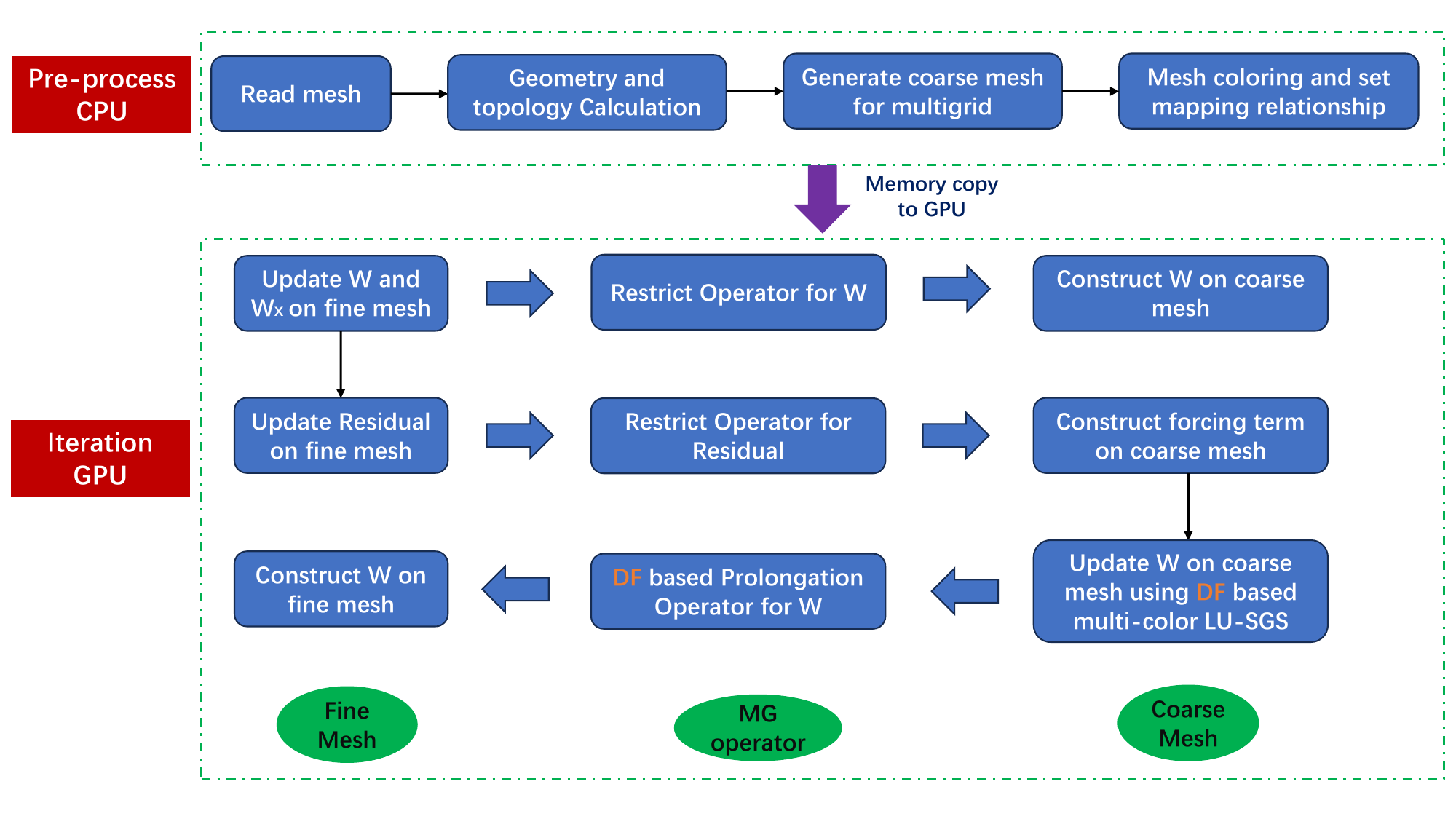}
	\caption{\label{GMG-gpu-procedure}
		An overall flowchart of GMG-CGKS solver based on GPU.}
\end{figure}

The total procedure of GMG-CGKS on GPU, including pre-process and multi-grid iteration, is shown in Fig.~\ref{GMG-gpu-procedure}.

\section{Numerical examples}
In this section, numerical tests are presented to validate the proposed scheme.
All the simulations are performed on 3-D hybrid unstructured mesh.
All simulations are performed on three-dimensional hybrid unstructured meshes, using an in-house C++ CUDA (version 12.4) solver.
The GPU hardware is the NVIDIA V100 GPU (NV-LINK version) with 7.8 TeraFLOPS for double precision and 900 GB/sec memory bandwidth.
The CPU hardware is the Intel Xeon Gold 5117 CPU with 14 cores, 2.0 GHz frequency and 115.2 GB/sec memory bandwidth.

\subsection{Subsonic flow around a cylinder}
\begin{figure}
    \centering
    \includegraphics[width=0.35\linewidth]{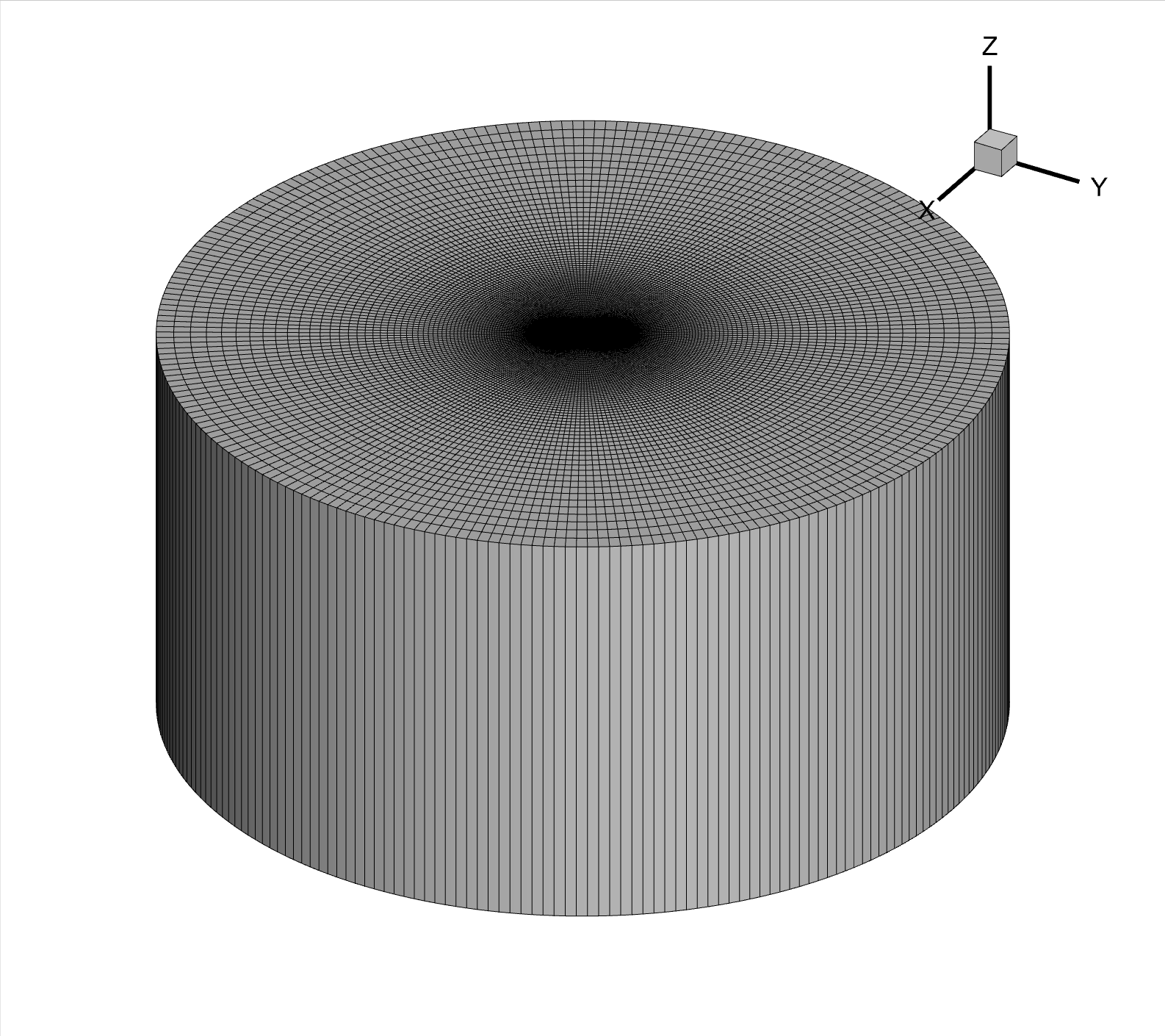}
    \includegraphics[width=0.35\linewidth]{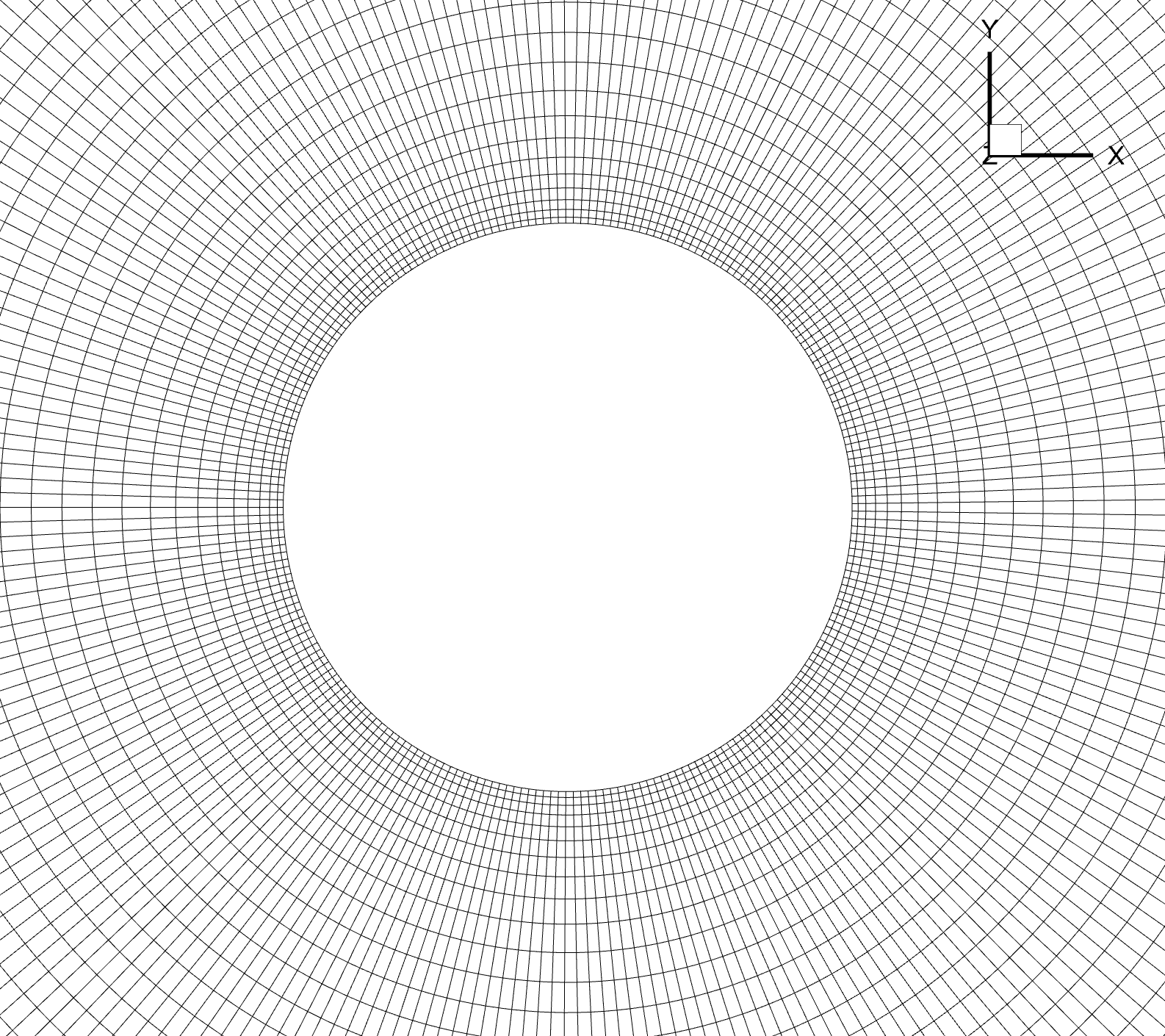}
    \caption{Mesh used in subsonic flow around a cylinder.}
    \label{cylinder-mesh}
\end{figure}

\begin{figure}
    \centering
    \includegraphics[width=0.35\linewidth]{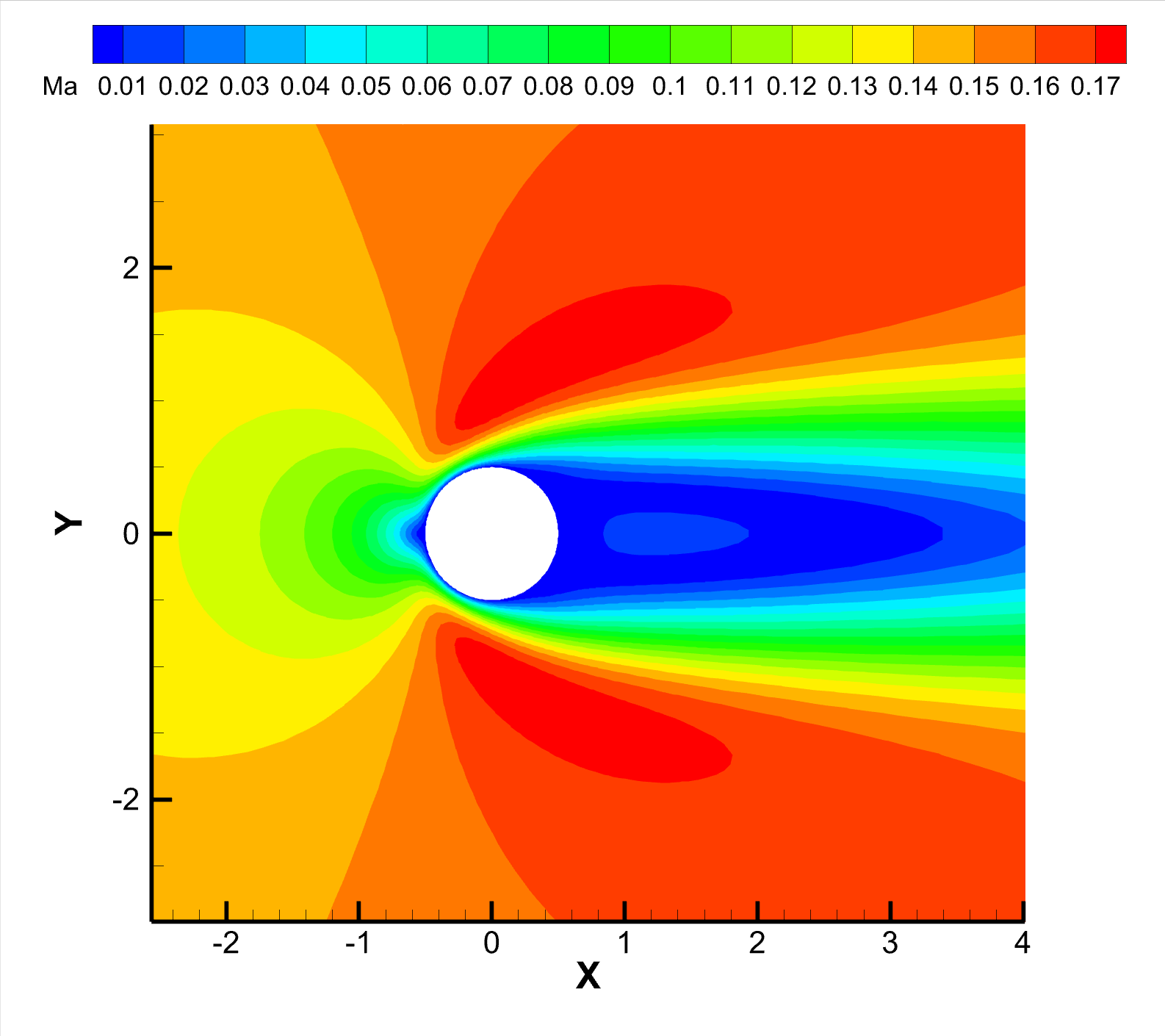}
    \includegraphics[width=0.35\linewidth]{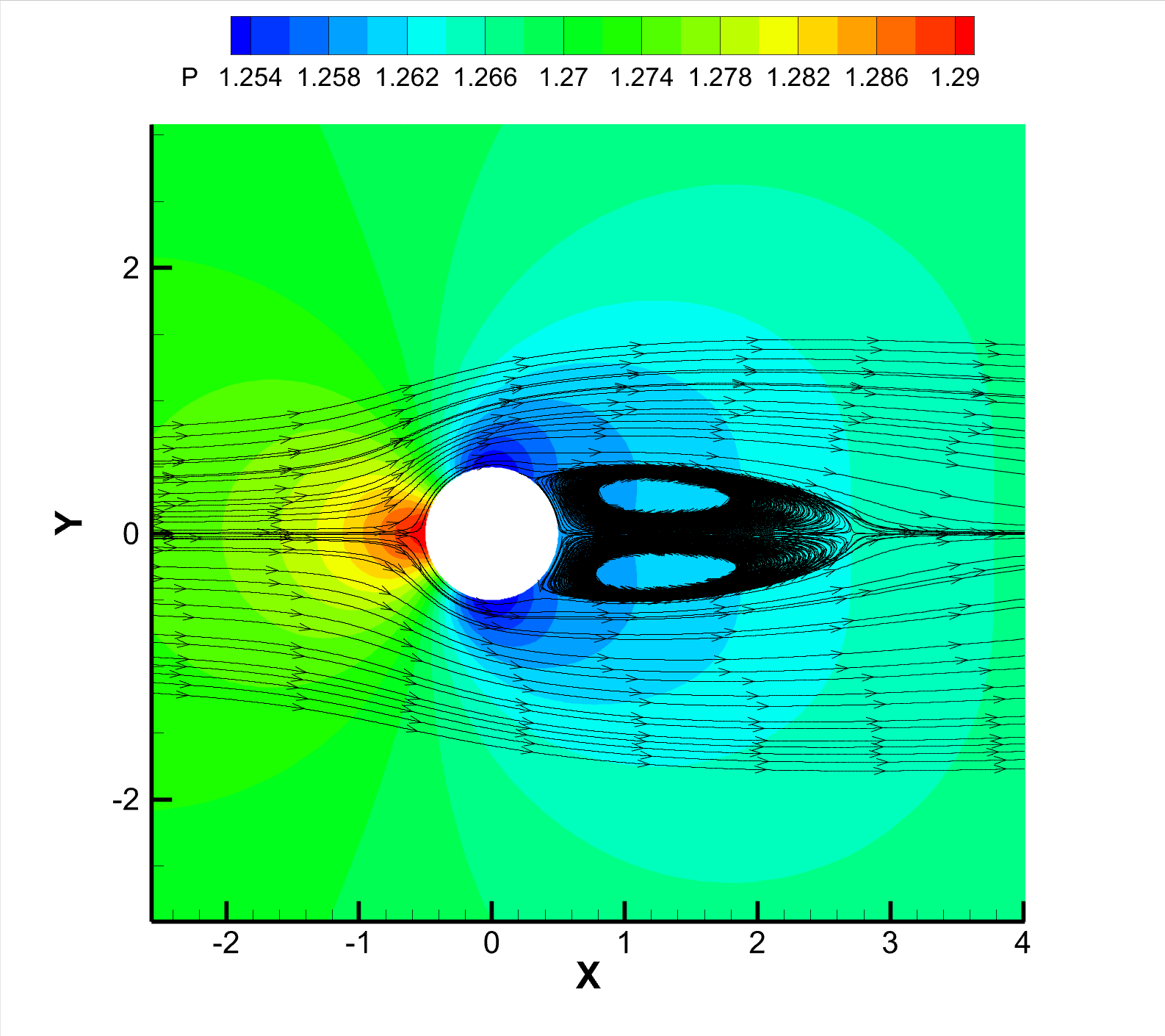}
    \caption{Subsonic flow around a cylinder. Left: Mach number contour. Right: Pressure contour.}
    \label{cylinder-contour}
\end{figure}

This section simulates subsonic laminar flow around a cylinder to validate the proposed GMG-CGKS convergence rate in nearly incompressible flow cases.
In this case, the incoming flow has a Mach number of 0.15 and a Reynolds number of 40, with the height of the first layer cell being 0.01.
This case uses a pure hexahedral mesh; the specific mesh distribution is shown in Fig .~\ref{cylinder-mesh}.
The Mach number contour, pressure contour, and streamline diagram are shown in Fig .~\ref{cylinder-contour}.
These contours demonstrate the high resolution of the GMG-CGKS method.
Table~\ref{cylinder-table} shows the quantitative results, indicating that the GMG-CGKS method aligns well with the experimental results, demonstrating its accuracy.

\begin{table}[htp]
	\small
	\begin{center}
		\def\temptablewidth{1.0\textwidth}
		{\rule{\temptablewidth}{1pt}}
		\begin{tabular*}{\temptablewidth}{@{\extracolsep{\fill}}c|c|c|c|c|c|c}
			Method & $C_D$ &  $C_L$ & $L$ &  Vortex Height & Vortex Width & $\theta$\\
			\hline
			Experiment \cite{tritton1959experiments} & $1.46-1.56$ & -- & -- & -- &  -- & --\\ 	
			Experiment \cite{coutanceau1977experimental} & -- & -- & 2.12 & 0.297 &  0.751 & $53.5^{\circ}$\\ 	
			DDG \cite{zhang2019direct} & 1.529 & -- & 2.31 & -- & -- & --\\ 	
			Current & 1.527 & $5.8e^{-14}$ & 2.22 & 0.296 &  0.714 & $53.5^{\circ}$\\ 	
		\end{tabular*}
		{\rule{\temptablewidth}{0.1pt}}
	\end{center}
	\vspace{-4mm} \caption{\label{cylinder-table} Comparison of the quantitative results of subsonic flow around a cylinder.}
\end{table}

\begin{figure}
    \centering
    \includegraphics[width=0.35\linewidth]{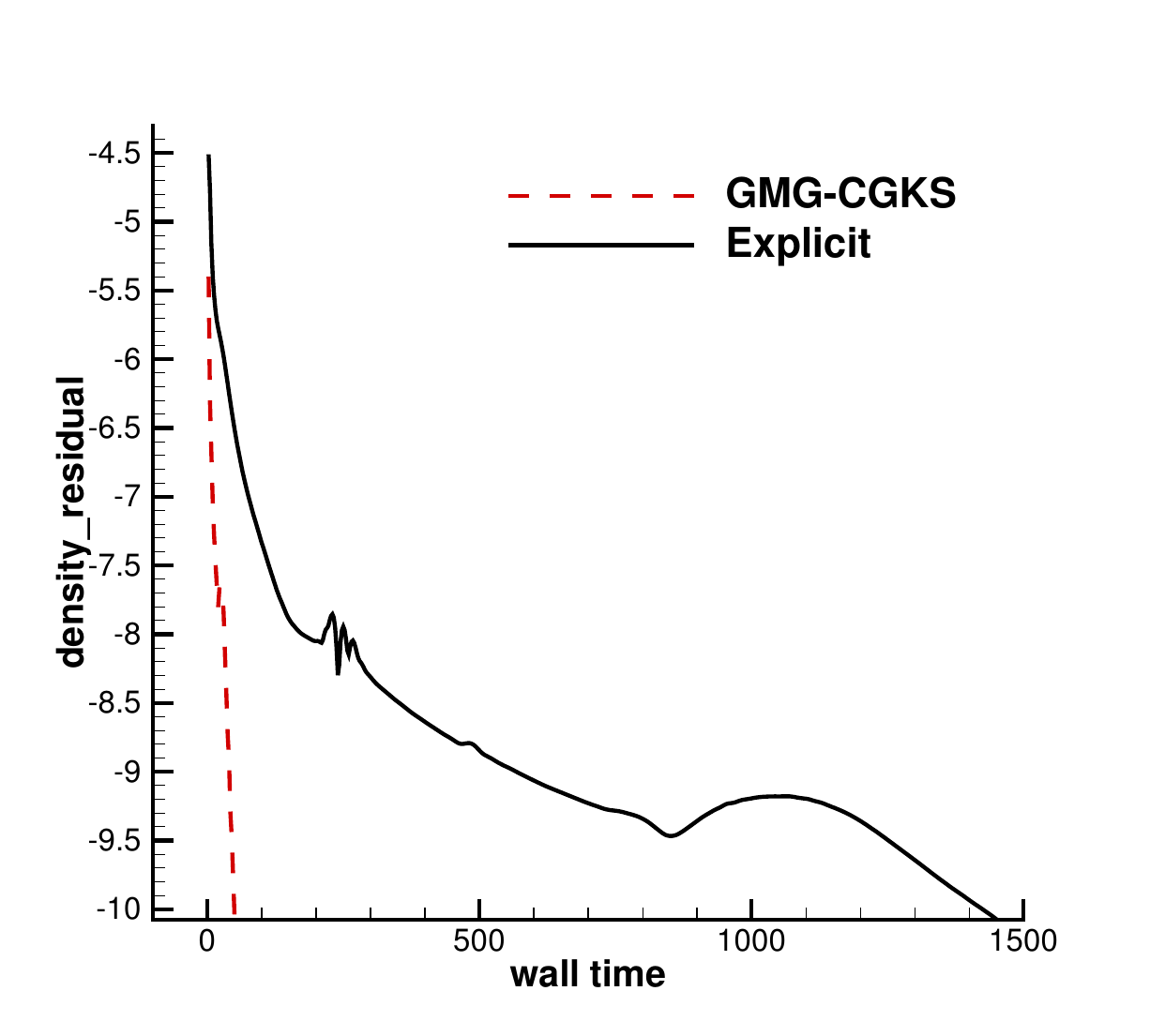}
    \includegraphics[width=0.35\linewidth]{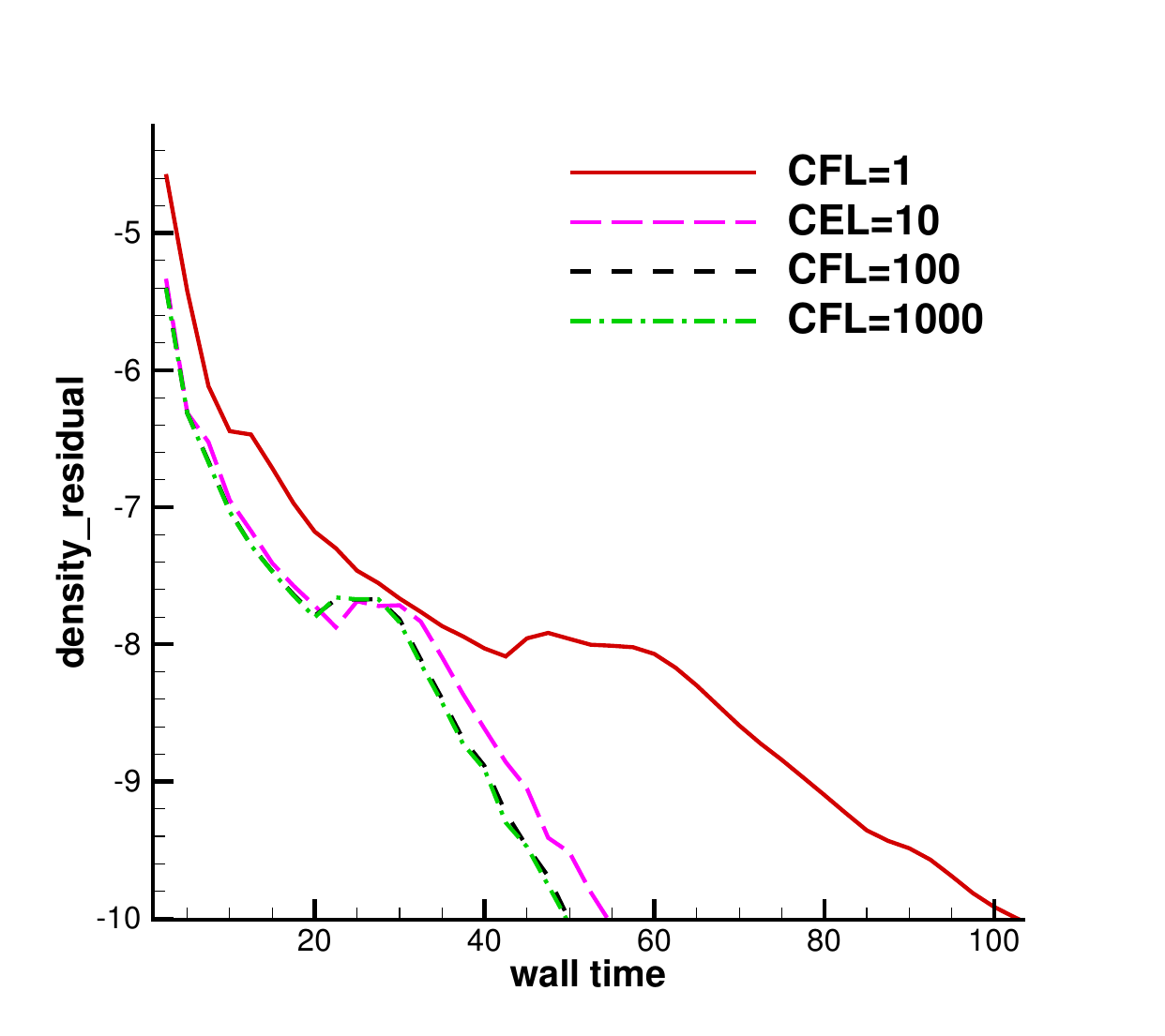}
    \caption{Subsonic flow around a cylinder. Left: The convergence history of GMG and explicit CGKS. Right: The convergence history of GMG under different CFL numbers (based on a fixed sweep number of six of LU-SGS).}
    \label{cylinder-GMG-1}
\end{figure}

\begin{figure}
    \centering
    \includegraphics[width=0.35\linewidth]{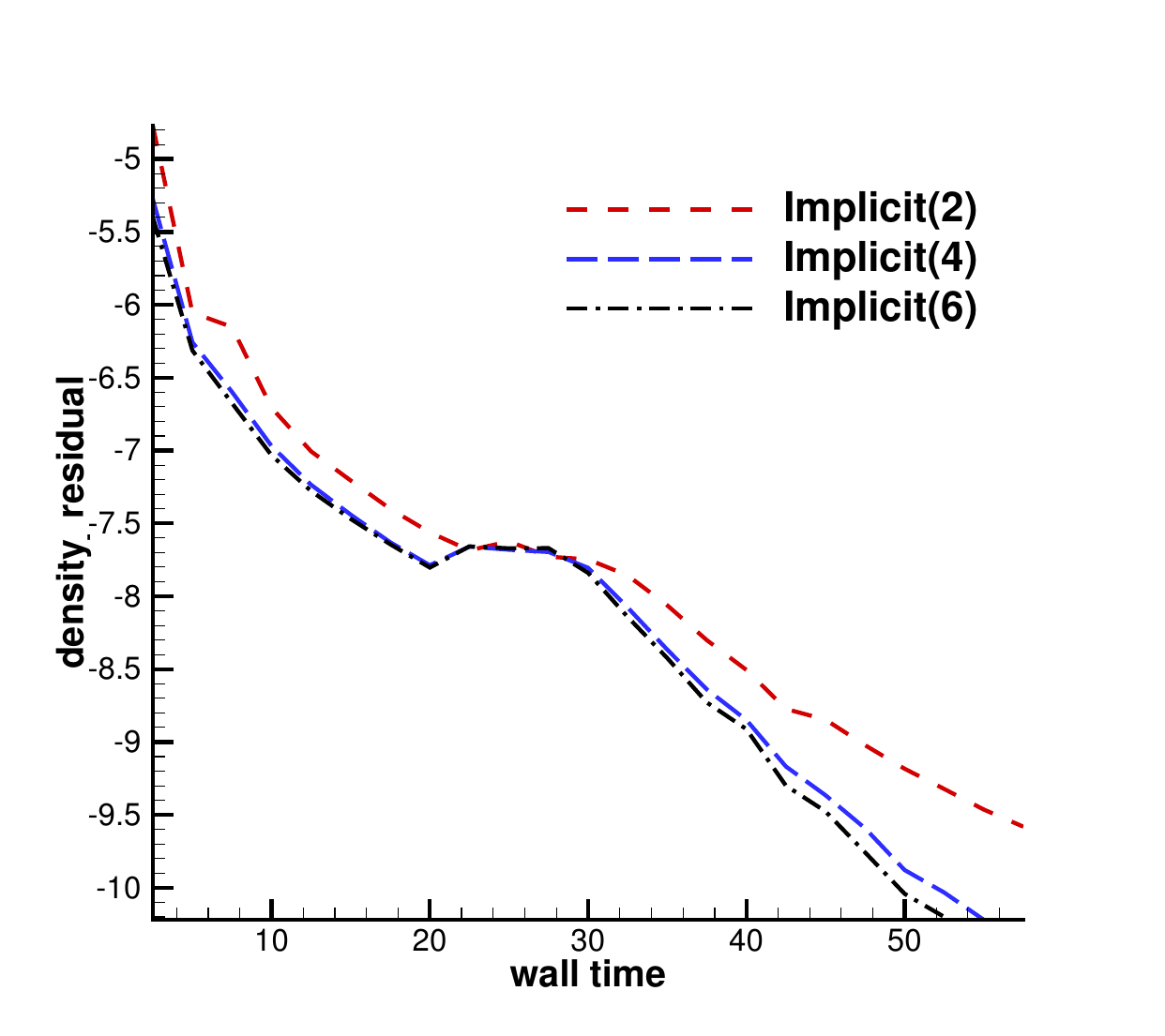}
    \includegraphics[width=0.35\linewidth]{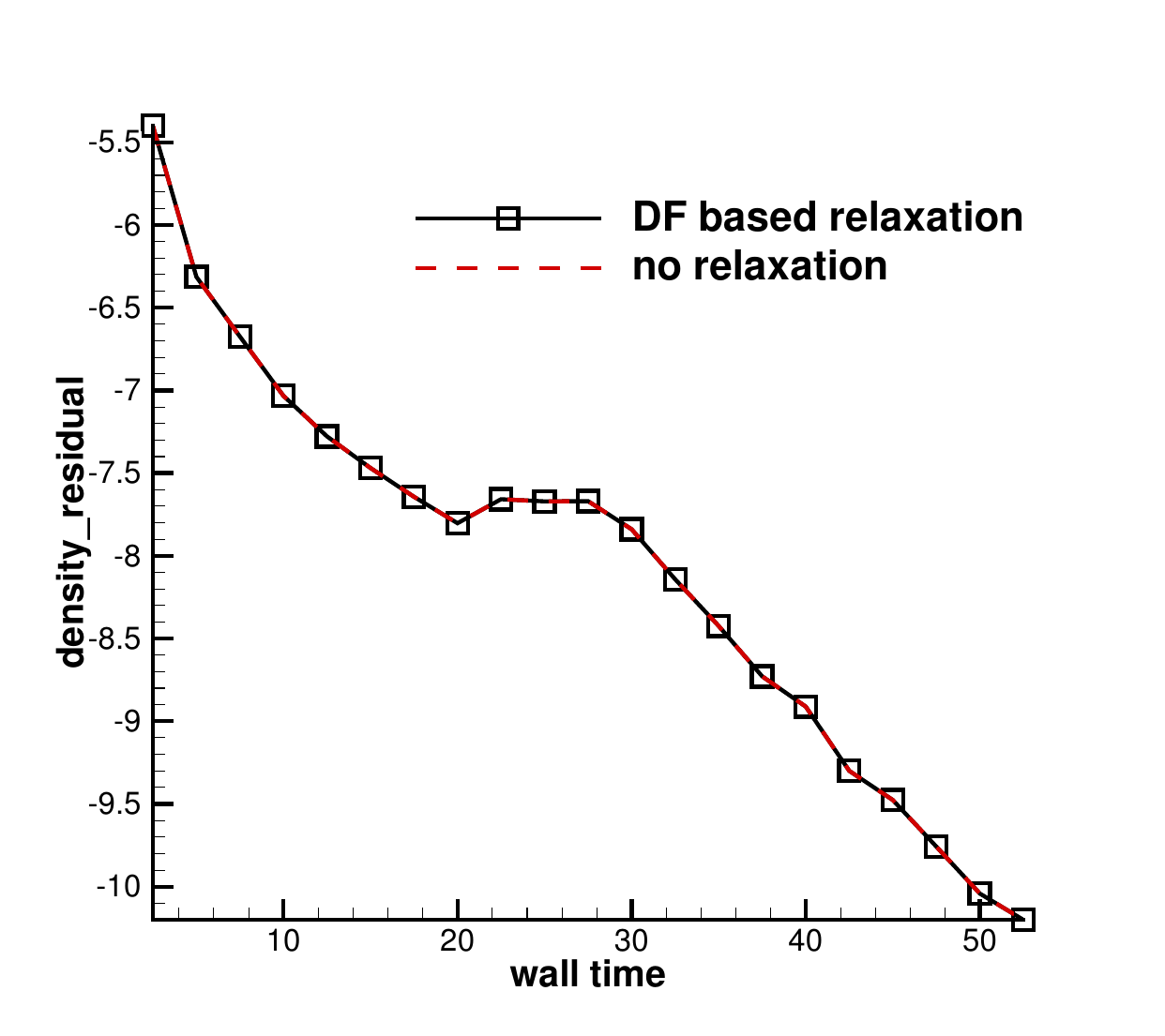}
    \caption{Subsonic flow around a cylinder. Left: The convergence history of GMG under different LU-SGS sweep numbers (based on a fixed CFL number) Right: The comparison with and without DF-based relaxation.}
    \label{cylinder-GMG-2}
\end{figure}

A series of comparisons of convergence rate for this case are shown in Fig .~\ref{cylinder-GMG-1} and Fig .~\ref{cylinder-GMG-2}, including a comparison between explicit and GMG, a comparison with a fixed CFL number while varying the LU-SGS sweep number in the coarse grid, a comparison with a fixed LU-SGS sweep number while varying the CFL number, and a comparison with and without DF-based relaxation.
Compared to the explicit scheme, the convergence rate of GMG is 32 times faster than the explicit scheme, requiring only forty seconds to converge, which is particularly notable given the already high computational speed of the GPU.
Regarding comparing CFL numbers, the GMG method does not require a high CFL value (greater than 10) to achieve a fast convergence rate.
At the same time, the comparison results indicate that the sweep count for LU-SGS should not be too low, with four to six sweeps being optimal.
In the final comparison regarding relaxation, the convergence rate of DF-based relaxation is consistent with that of the no-relaxation case. This is because the case is very smooth, allowing the DF-based relaxation to adaptively revert to the no-relaxation scheme.

\subsection{Subsonic flow around a NACA0012 airfoil}
\begin{figure}
    \centering
    \includegraphics[width=0.35\linewidth]{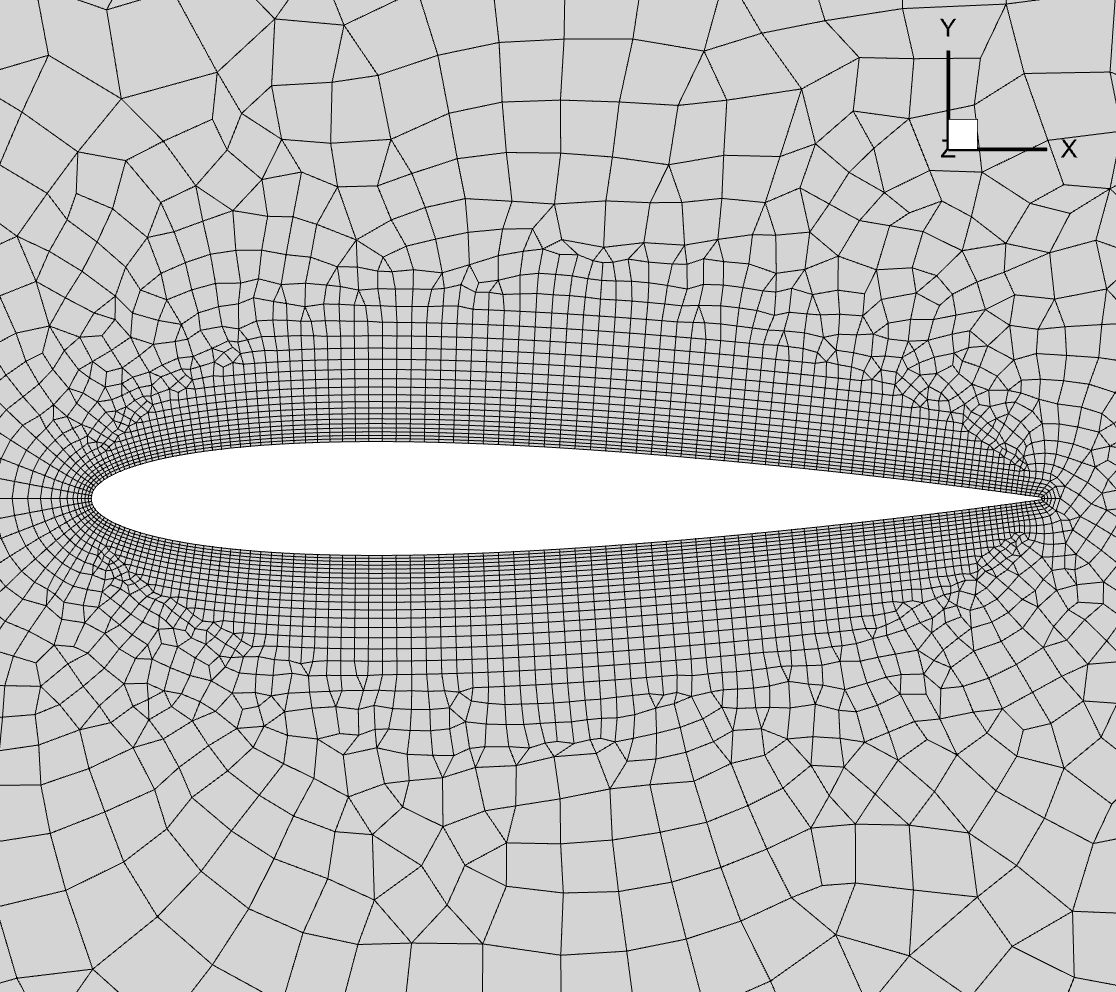}
    \includegraphics[width=0.35\linewidth]{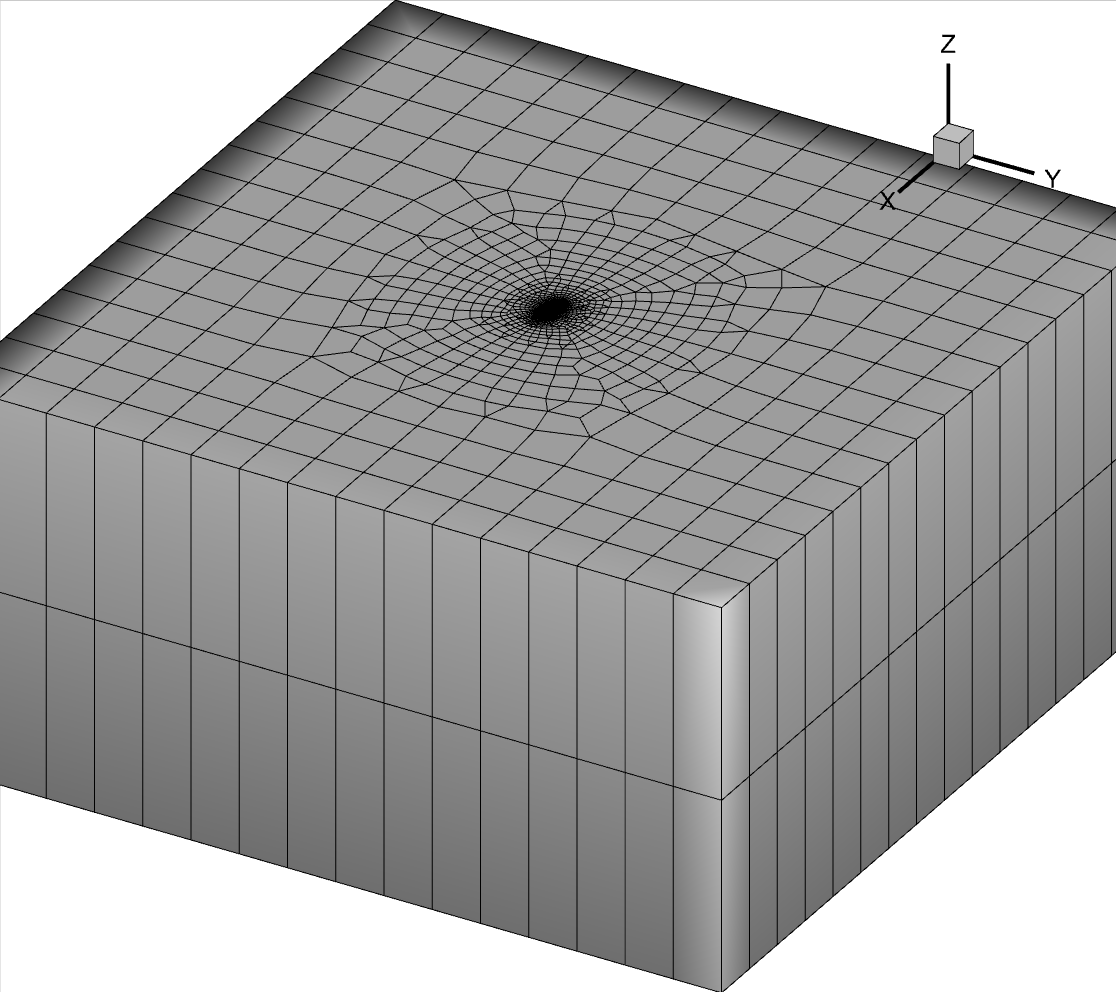}
    \caption{Mesh used in subsonic flow around a NACA0012 airfoil.}
    \label{naca-mesh}
\end{figure}

\begin{figure}
    \centering
    \includegraphics[width=0.35\linewidth]{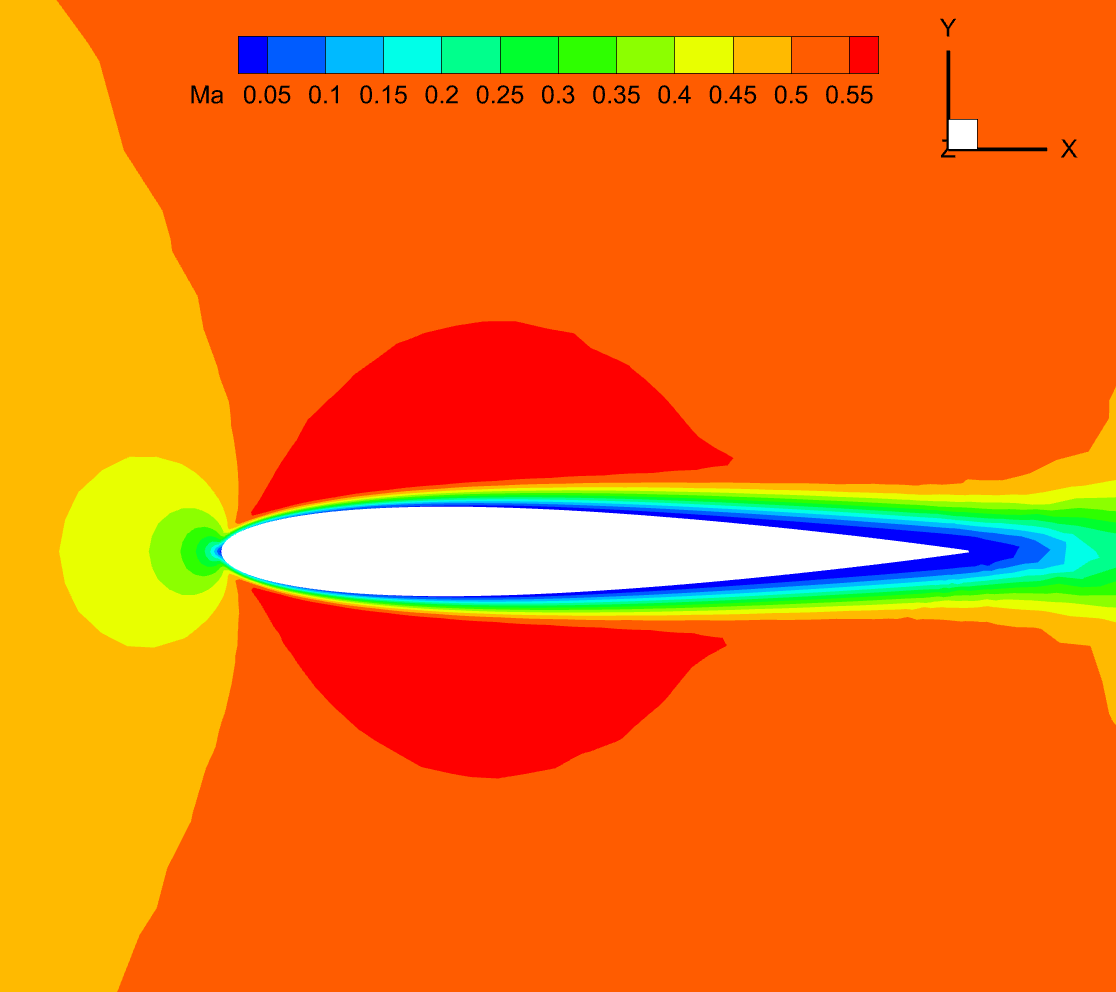}
    \includegraphics[width=0.35\linewidth]{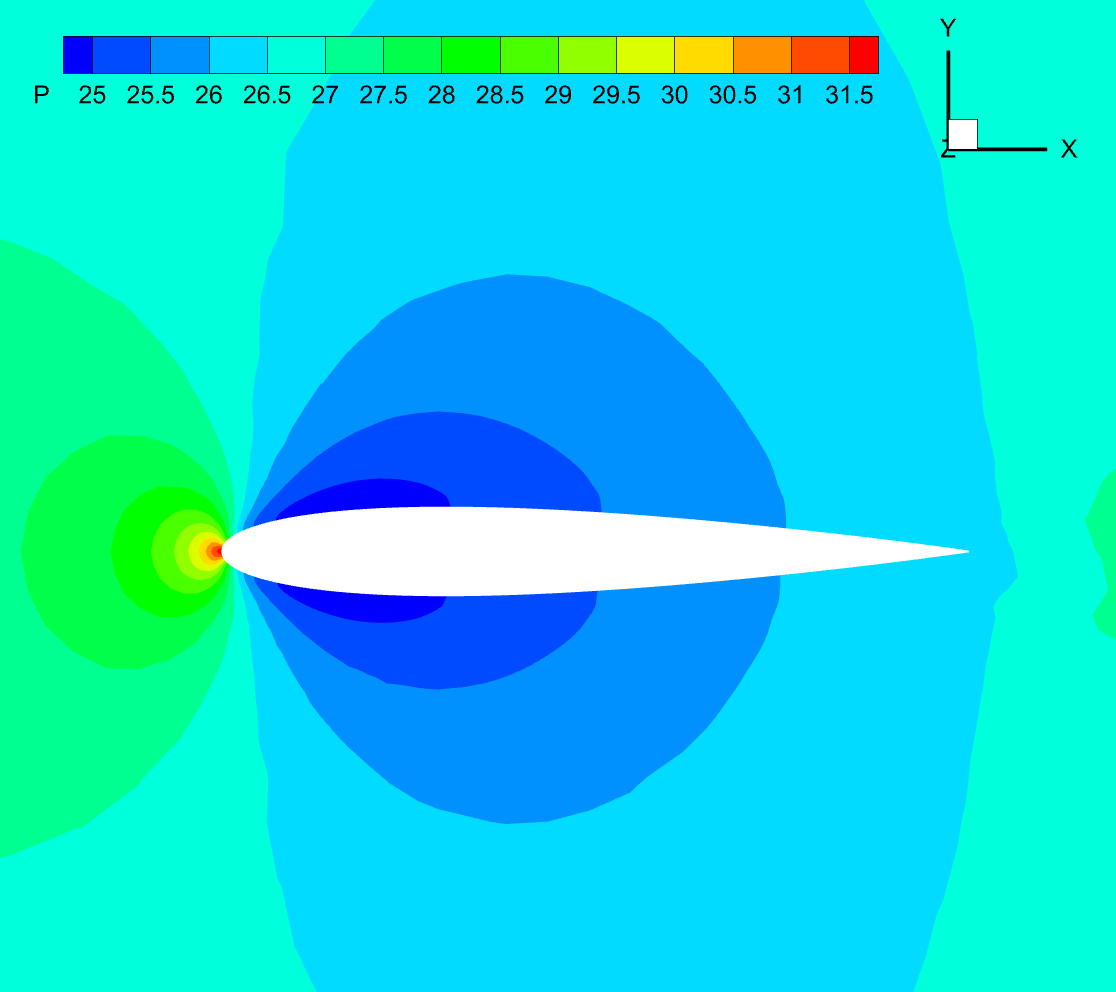}
    \caption{Subsonic flow around a NACA0012 airfoil. Left: Mach number contour. Right: Pressure contour.}
    \label{naca-contour}
\end{figure}

\begin{figure}
    \centering
    \includegraphics[width=0.35\linewidth]{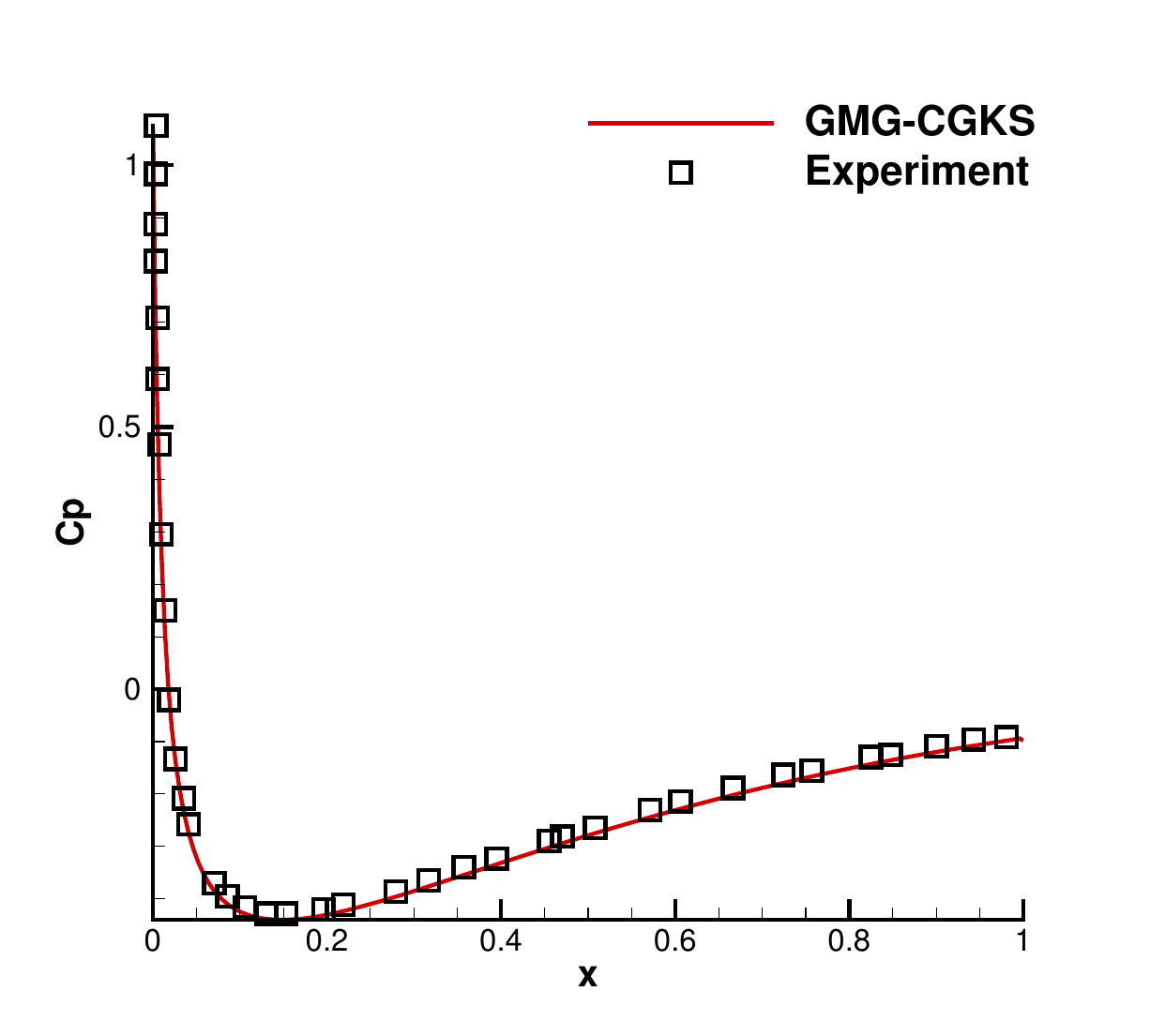}
    \caption{Subsonic flow around a NACA0012 airfoil. Surface pressure coefficient distribution.}
    \label{naca-cp}
\end{figure}

Subsonic laminar flow around a NACA0012 airfoil is a common benchmark to validate steady-state acceleration techniques' accuracy and convergence rate.
So, a subsonic flow around a NACA0012 airfoil is simulated with an incoming Mach number equal to 0.5, and a Reynolds number equal to 5000 based on the chord length L = 1 is simulated.
A hybrid unstructured mesh constructed of prisms and hexahedrons with a mesh number equal to 13076 is used for the simulation, as shown in Fig.~\ref{naca-mesh}.
Mach Number contour and pressure contour are shown in Fig.~\ref{naca-contour}.
The quantitative result, including the surface pressure coefficient, is extracted and plotted in Fig.~\ref{naca-cp}, which highly agrees with the Ref \cite{bassi1997ns}.

\begin{figure}
    \centering
    \includegraphics[width=0.35\linewidth]{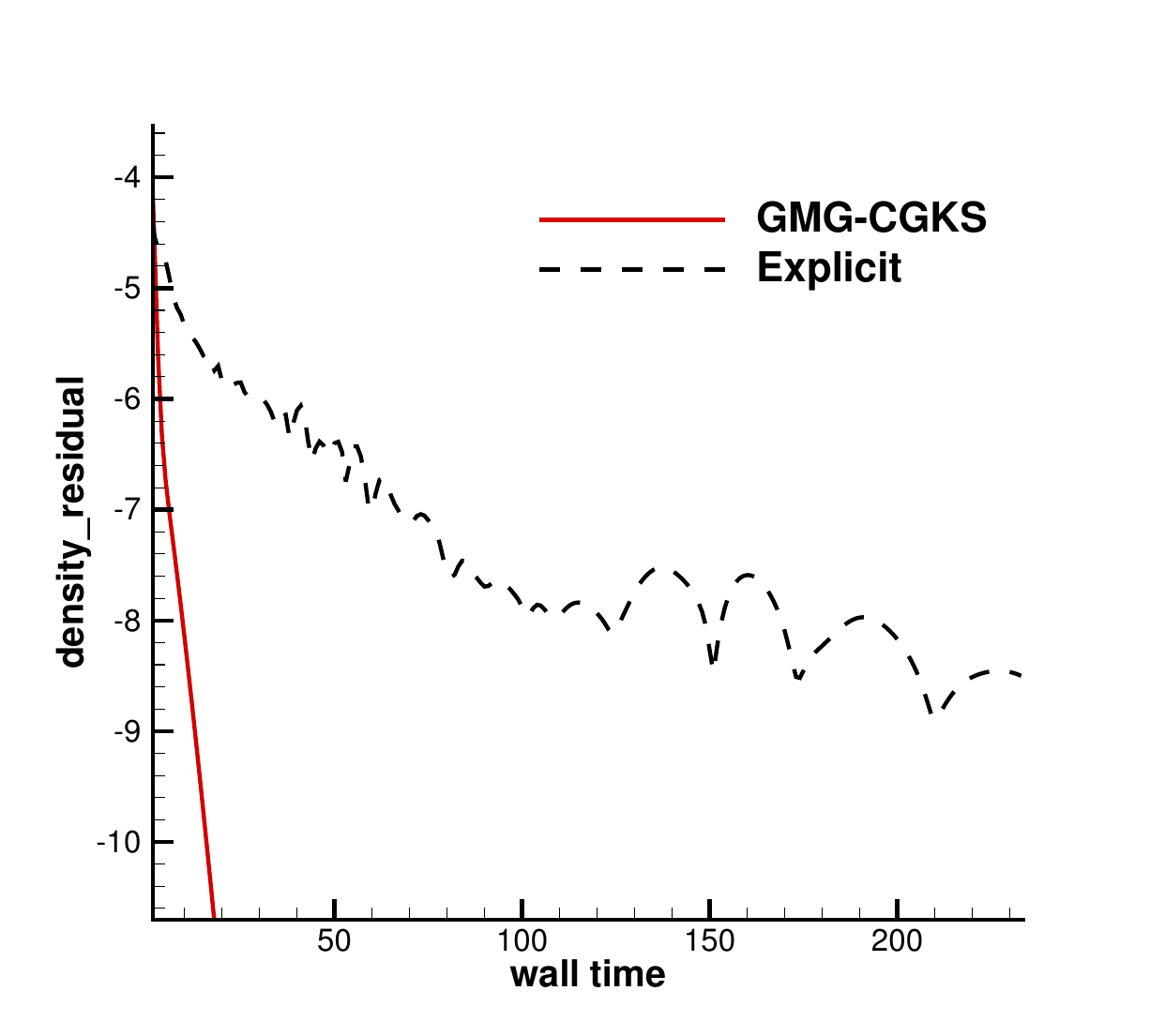}
    \includegraphics[width=0.35\linewidth]{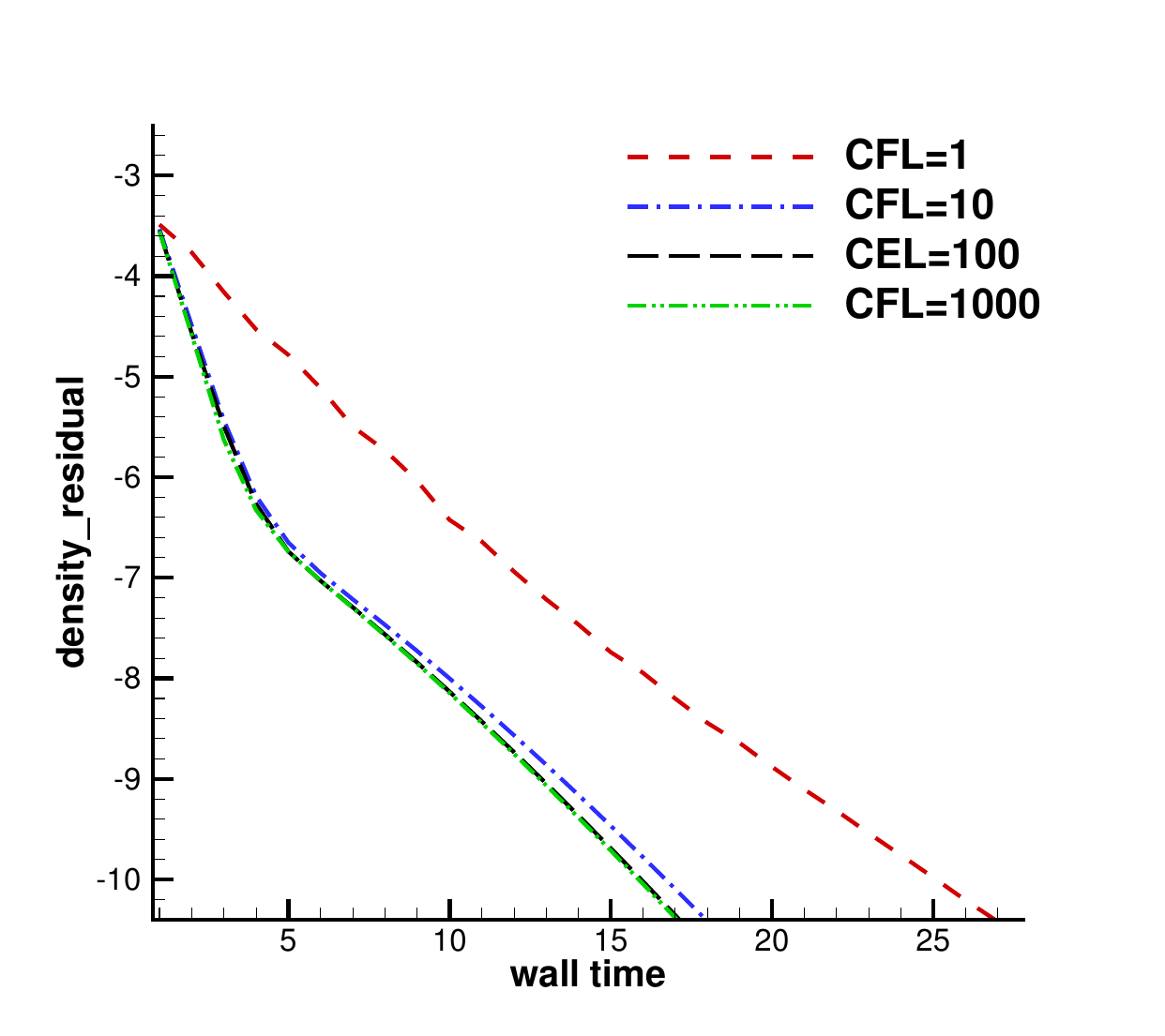}
    \caption{Subsonic flow around a NACA0012 airfoil. Left: The convergence history of GMG and explicit CGKS. Right: The convergence history of GMG under different CFL numbers (based on a fixed sweep number of six of LU-SGS).}
    \label{NACA-GMG-1}
\end{figure}

\begin{figure}
    \centering
    \includegraphics[width=0.35\linewidth]{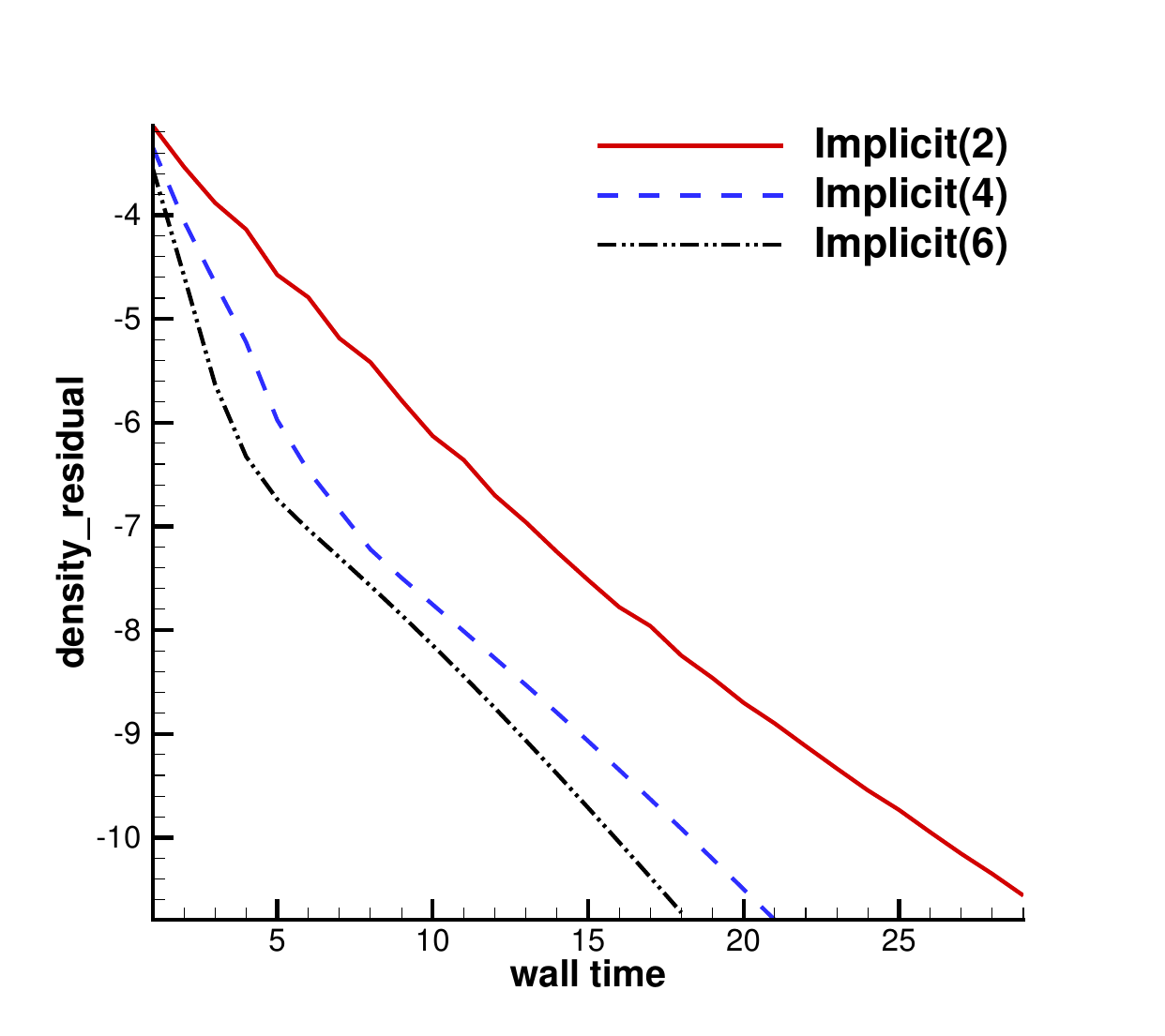}
    \includegraphics[width=0.35\linewidth]{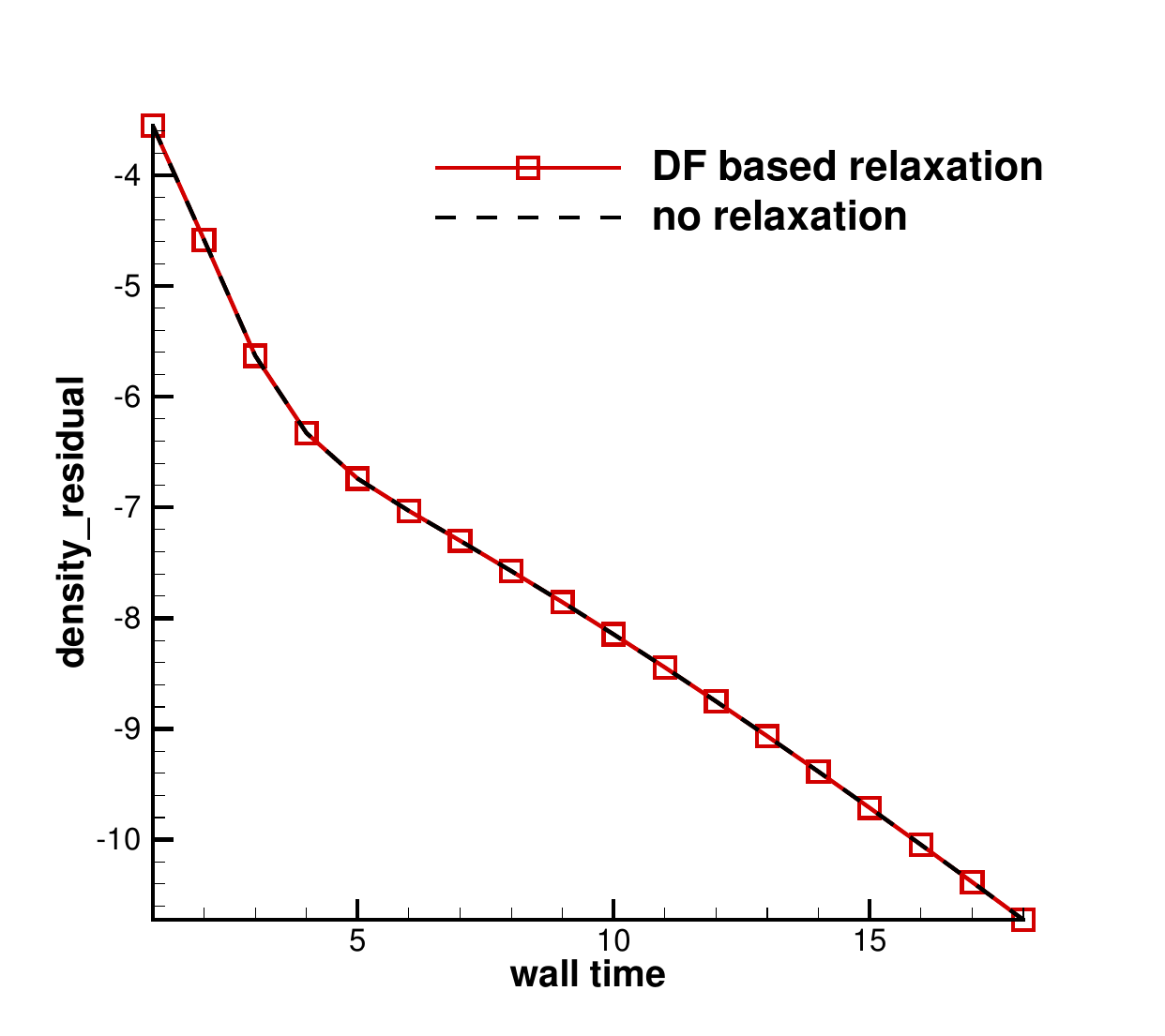}
    \caption{Subsonic flow around a NACA0012 airfoil. Left: The convergence history of GMG under different LU-SGS sweep numbers (based on a fixed CFL number) Right: The comparison with and without DF-based relaxation.}
    \label{NACA-GMG-2}
\end{figure}

A series of comparisons of convergence rate for this case are shown in Fig .~\ref{NACA-GMG-1} and Fig .~\ref{NACA-GMG-2}, including a comparison between explicit and GMG, a comparison with a fixed CFL number while varying the LU-SGS sweep number in the coarse grid, a comparison with a fixed LU-SGS sweep number while varying the CFL number, and a comparison with and without DF-based relaxation.
Compared to the explicit scheme, GMG's convergence rate is 18 times faster, requiring only 20 seconds to converge.
This is exceptionally effective for GPU-based computations.
The comparison results indicate that GMG requires only a CFL number of 10 or higher to achieve a significant speedup, which helps avoid the risk of computational instability caused by excessively high CFL numbers.
At the same time, the sweep count for LU-SGS should not be too low, as this may lead to a deterioration in convergence speed.
In the final comparison regarding relaxation, the convergence rate of DF-based relaxation is consistent with that of the no-relaxation case. This is because the case is very smooth, allowing the DF-based relaxation to adaptively revert to the no-relaxation scheme.

\subsection{Transonic flow around a dual-NACA0012 airfoil}
\begin{figure}
    \centering
    \includegraphics[width=0.35\linewidth]{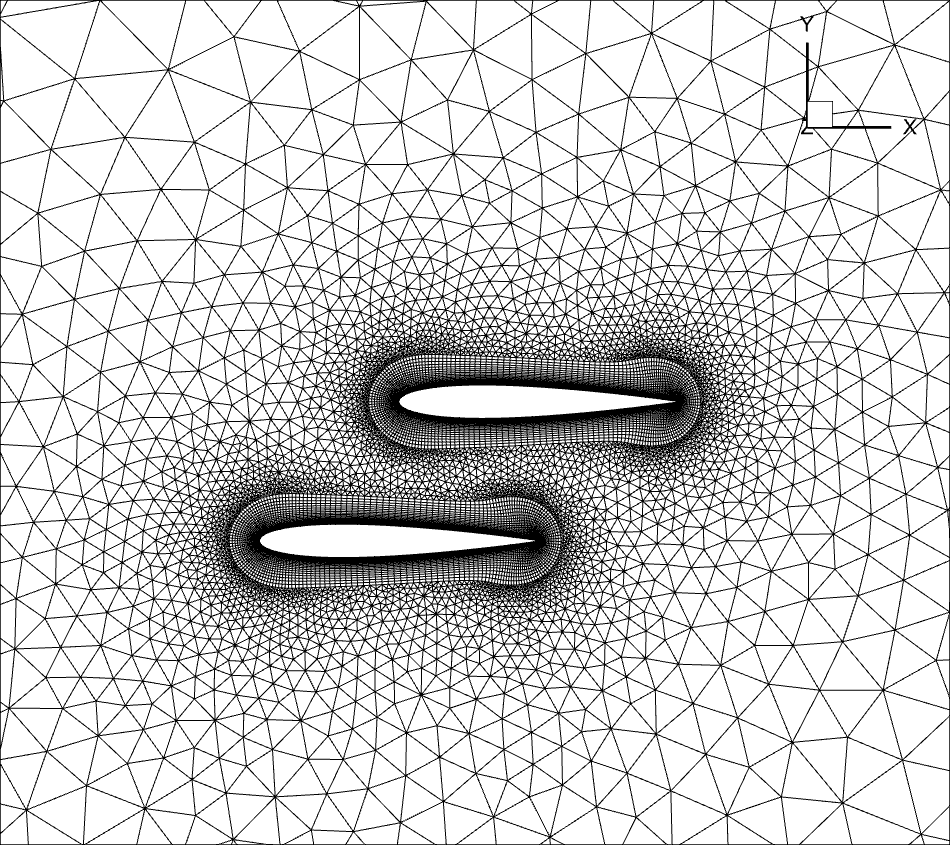}
    \includegraphics[width=0.35\linewidth]{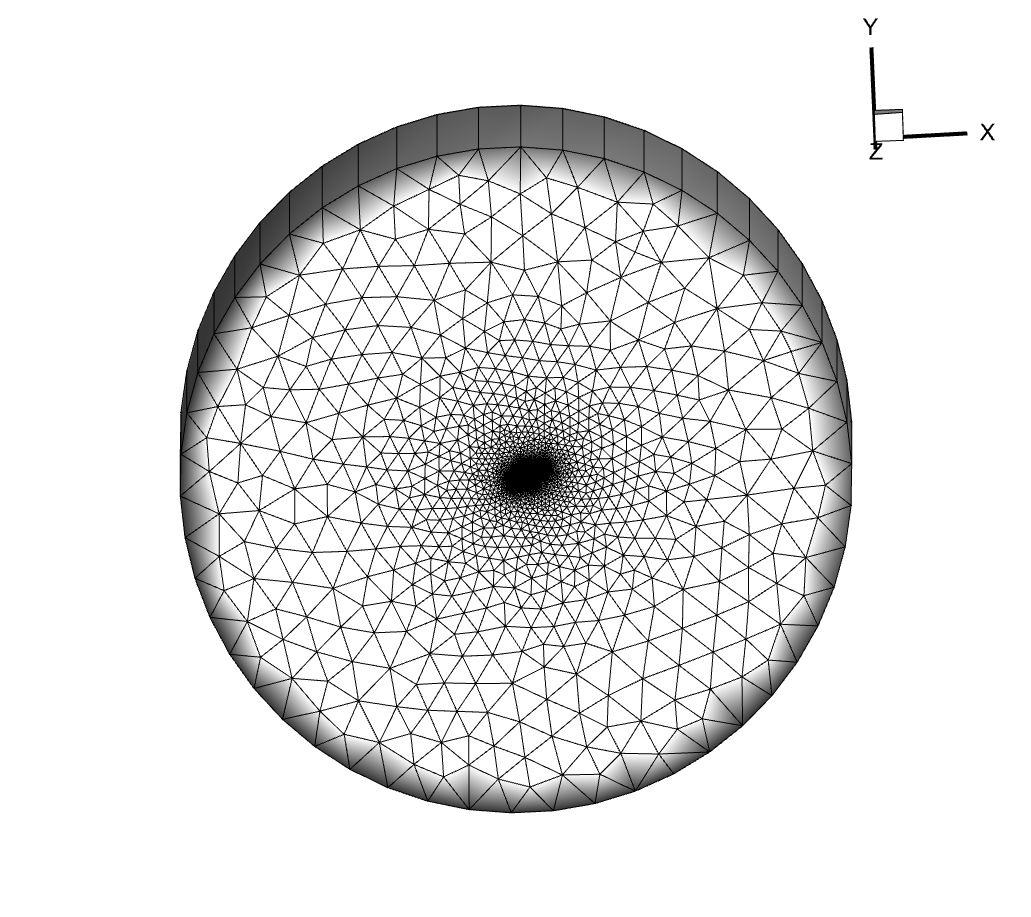}
    \caption{Mesh used in transonic flow around dual-NACA0012 airfoil.}
    \label{dualnaca-mesh}
\end{figure}

\begin{figure}
    \centering
    \includegraphics[width=0.35\linewidth]{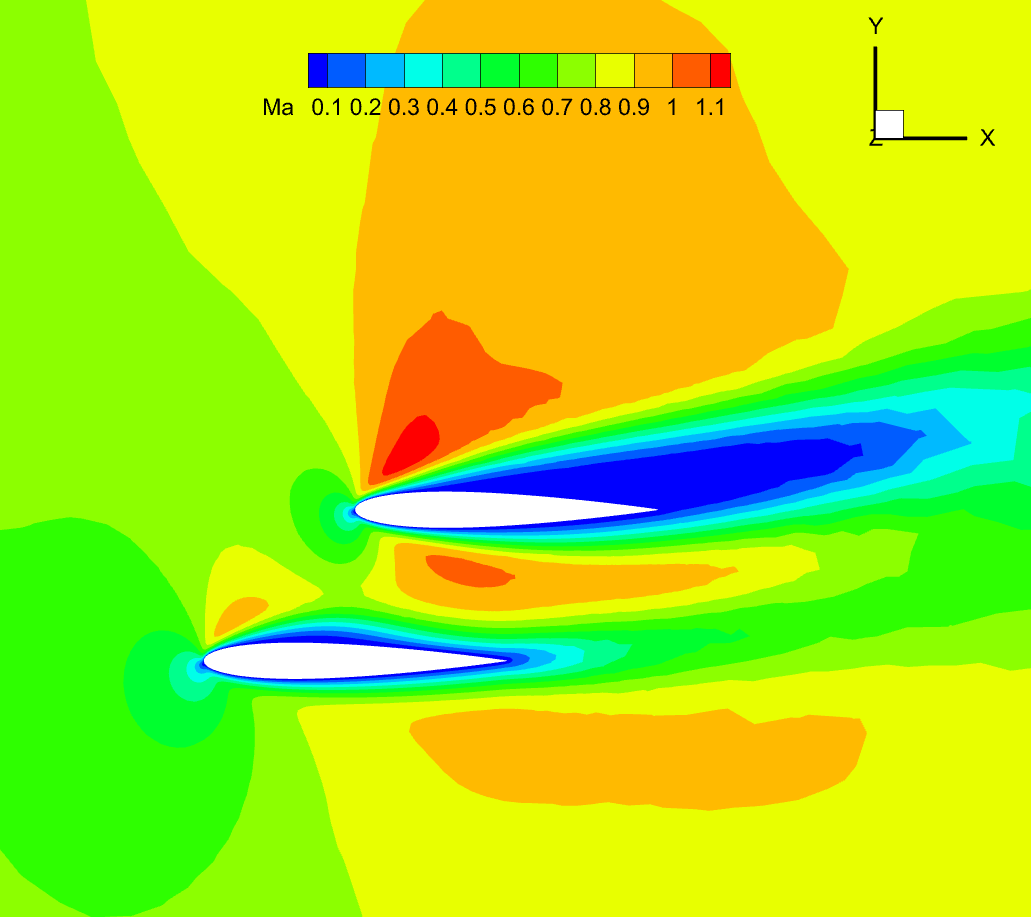}
    \includegraphics[width=0.35\linewidth]{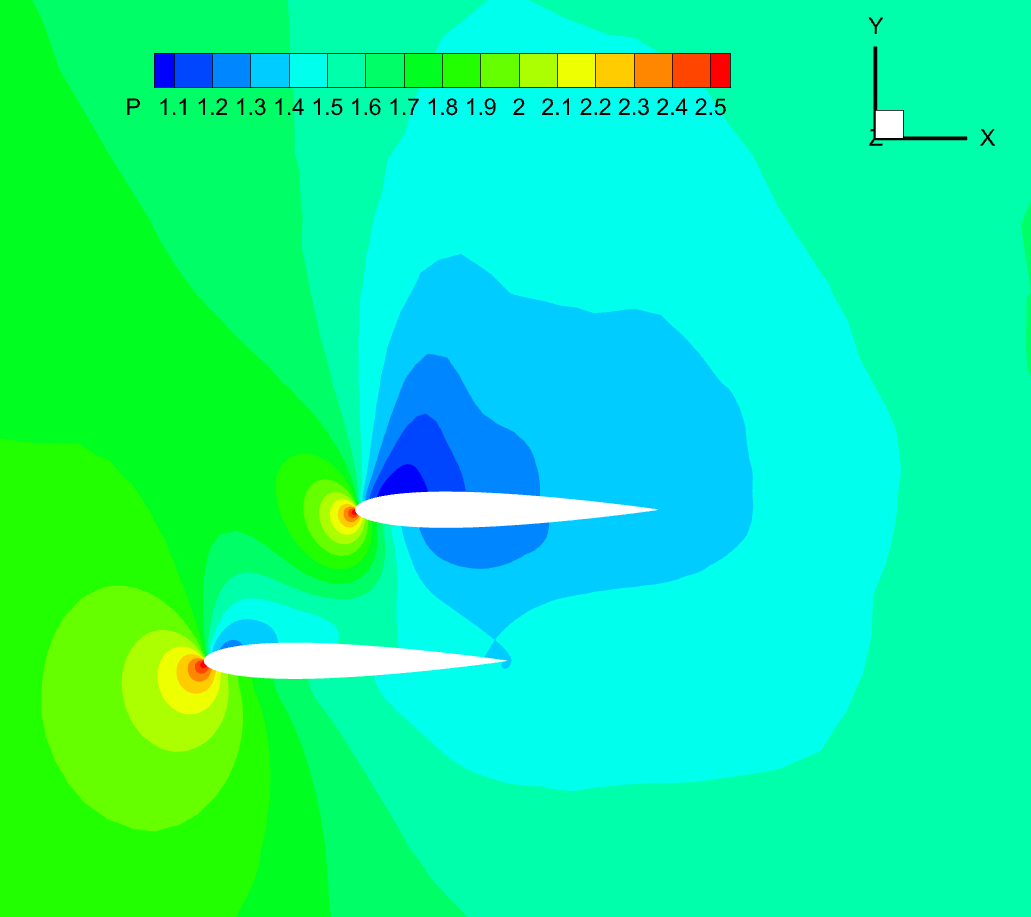}
    \caption{Transonic flow around dual-NACA0012. Left: Mach number contour. Right: Pressure contour.}
    \label{dualnaca-contour}
\end{figure}

\begin{figure}
    \centering
    \includegraphics[width=0.35\linewidth]{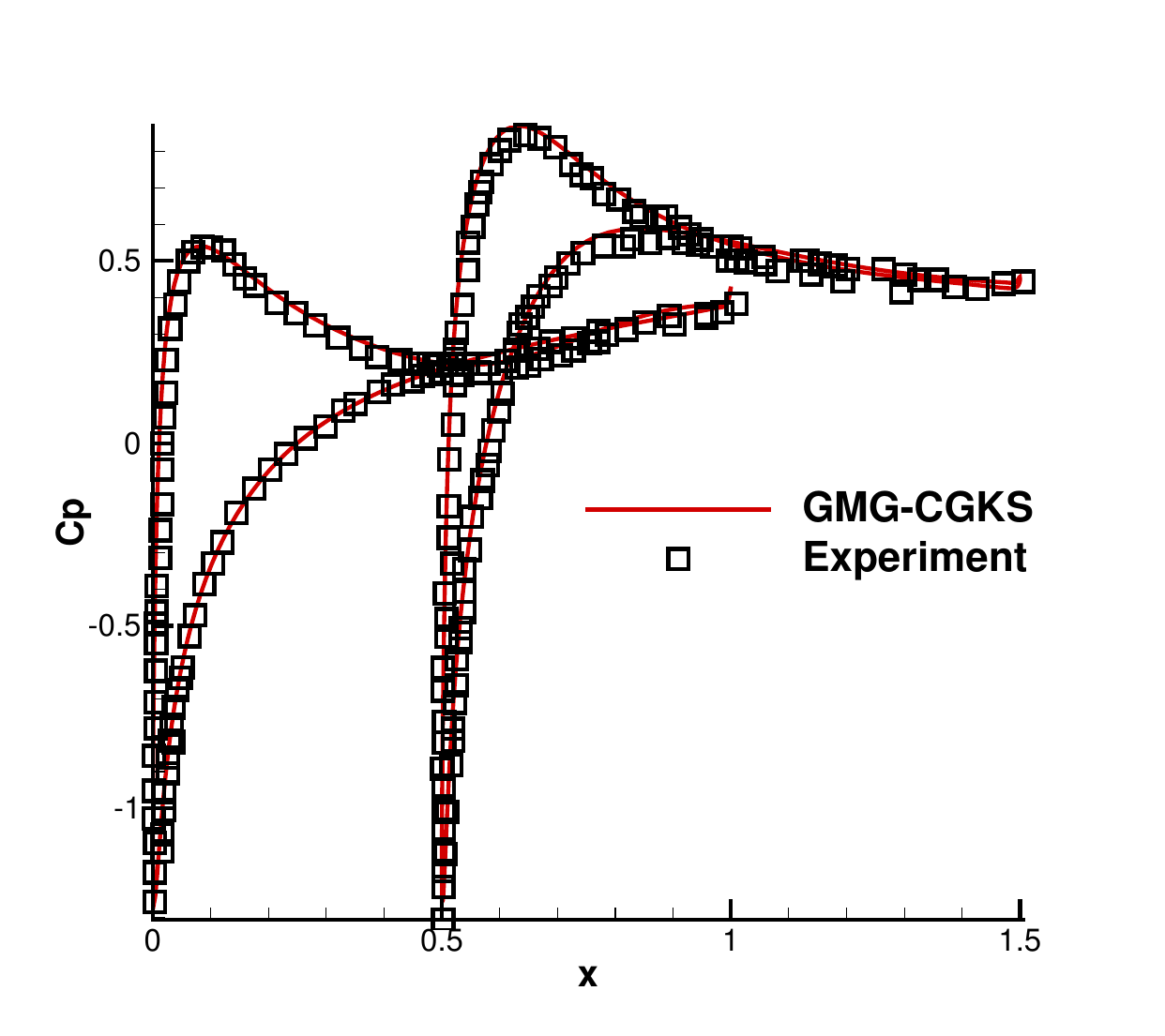}
    \caption{Transonic flow around dual-NACA0012 airfoil. Surface pressure coefficient distribution.}
    \label{dualnaca-cp}
\end{figure}

A transonic flow around the dual-NACA0012 is simulated to demonstrate that the current GMG-CGKS method accelerates the convergence rate and provides high resolution to accurately capture the flow field's structure.
The incoming Mach number is set
to be 0.8 with an angle of attack AOA = 10, and the Reynolds number is set to be 500
based on the chord length L =1.
The mesh has 28678 mixed elements, as shown in Fig.~\ref{dualnaca-mesh}.
The near wall size of the mesh is set to be h = $2 \times 10^{-3}$, which indicates that the grid Reynolds number is $2.5 \times 10^{5}$.
The far-field boundary condition is set to be subsonic inflow using Riemann invariants, and the wall is set to be a non-slip adiabatic wall.
The Mach number distribution and the pressure distribution are shown in Fig.~\ref{dualnaca-contour}.
The oblique shock wave can be observed at the front of the top airfoil.
The surface pressure coefficient is also extracted and compared with the reference data \cite{jawahar2000high}, as shown in Fig.~\ref{dualnaca-cp}. The result obtained by GMG-CGKS agrees well with the experimental data.

\begin{figure}
    \centering
    \includegraphics[width=0.35\linewidth]{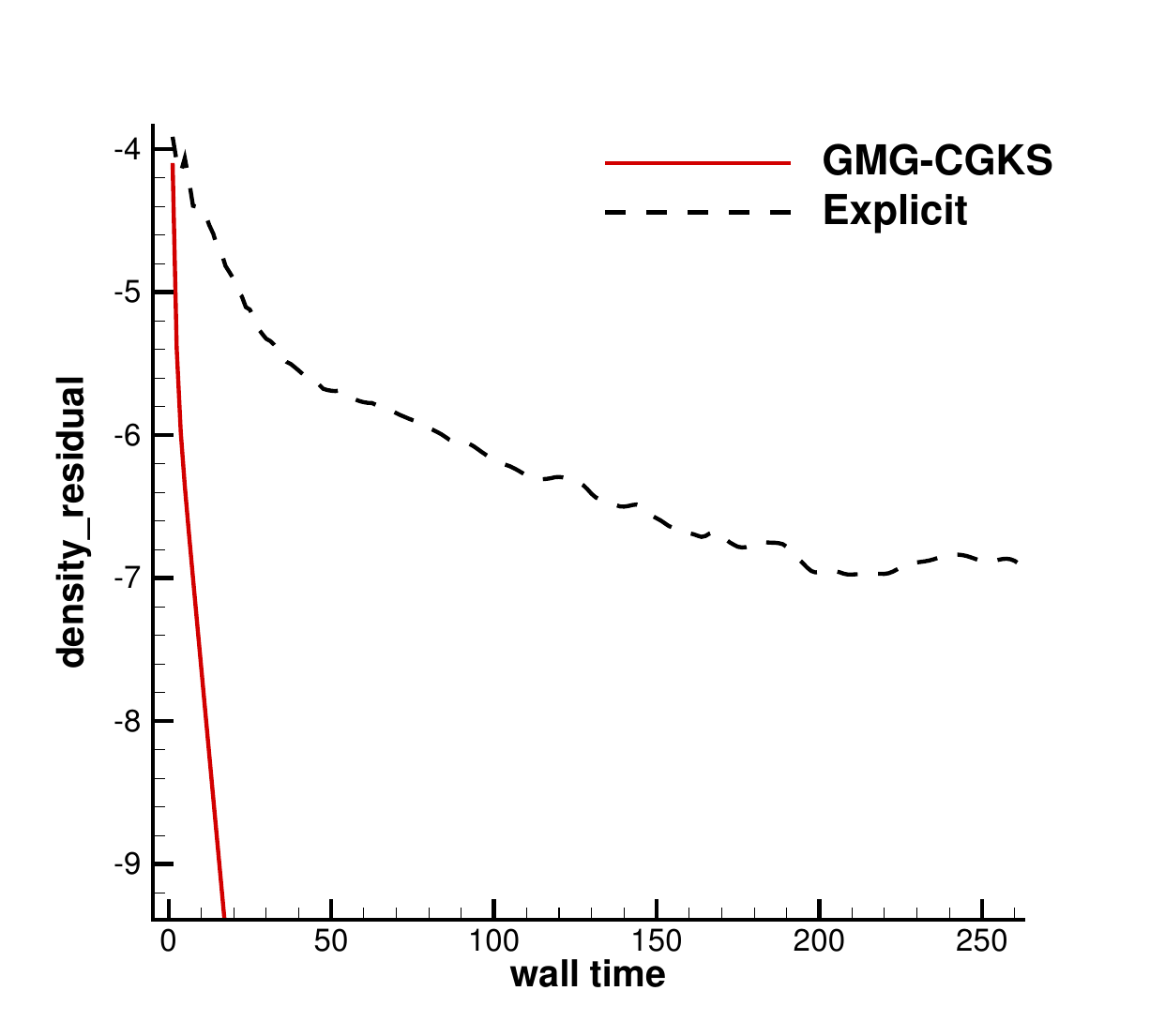}
    \includegraphics[width=0.35\linewidth]{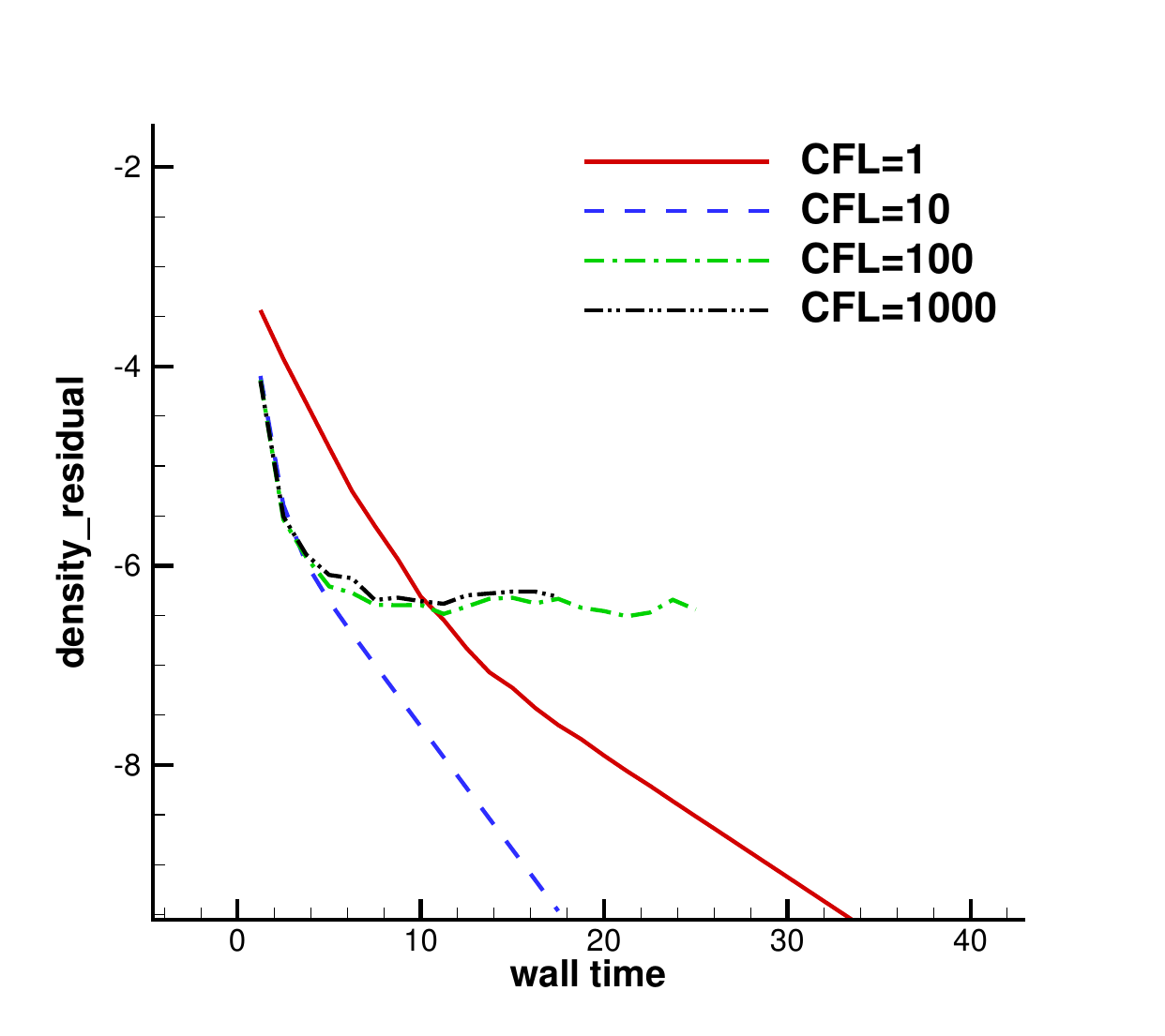}
    \caption{Transonic flow around dual-NACA0012 airfoil. Left: The convergence history of GMG and explicit CGKS. Right: The convergence history of GMG under different CFL numbers (based on a fixed sweep number of six of LU-SGS).}
    \label{dualNACA-GMG-1}
\end{figure}

\begin{figure}
    \centering
    \includegraphics[width=0.35\linewidth]{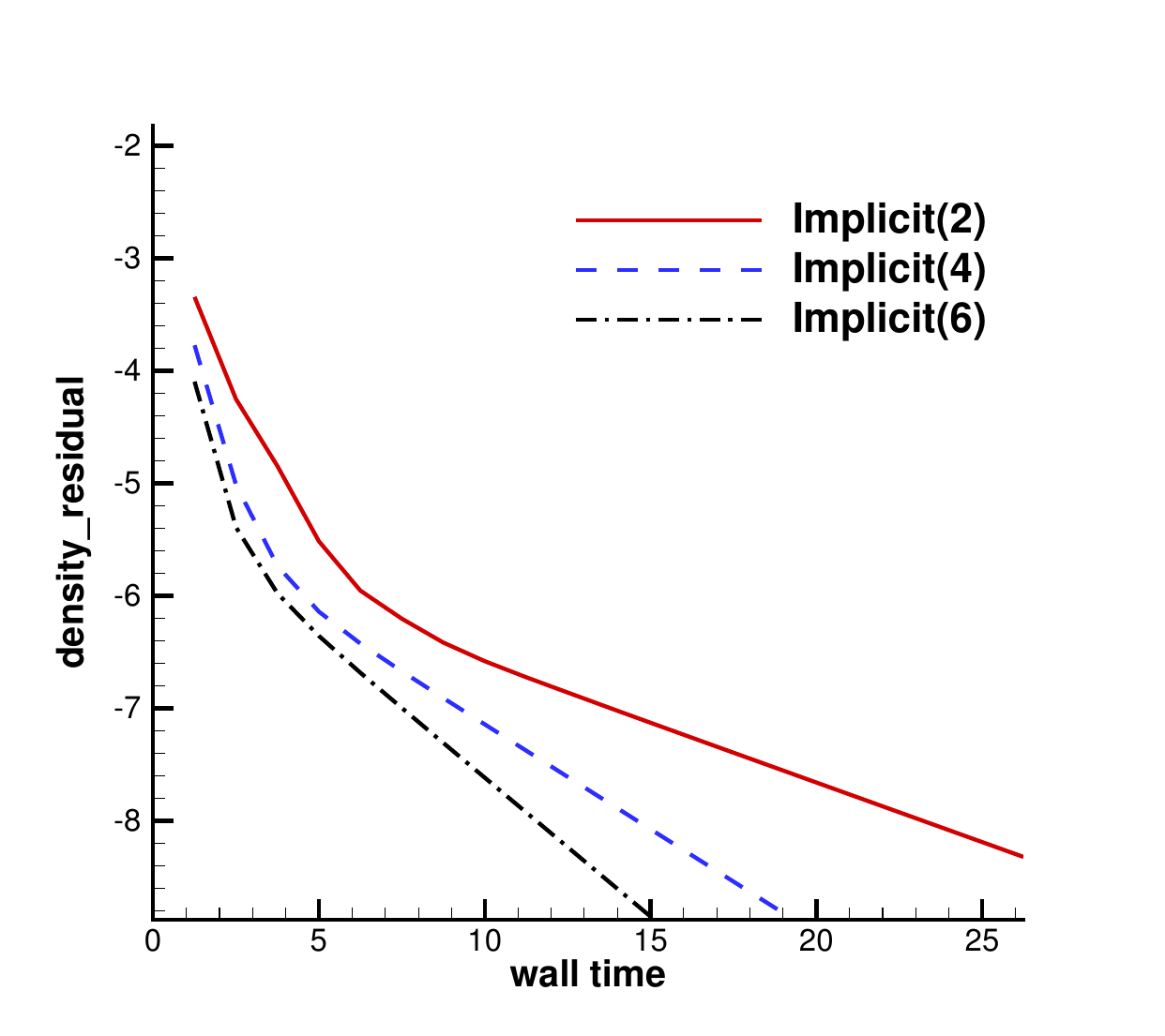}
    \includegraphics[width=0.35\linewidth]{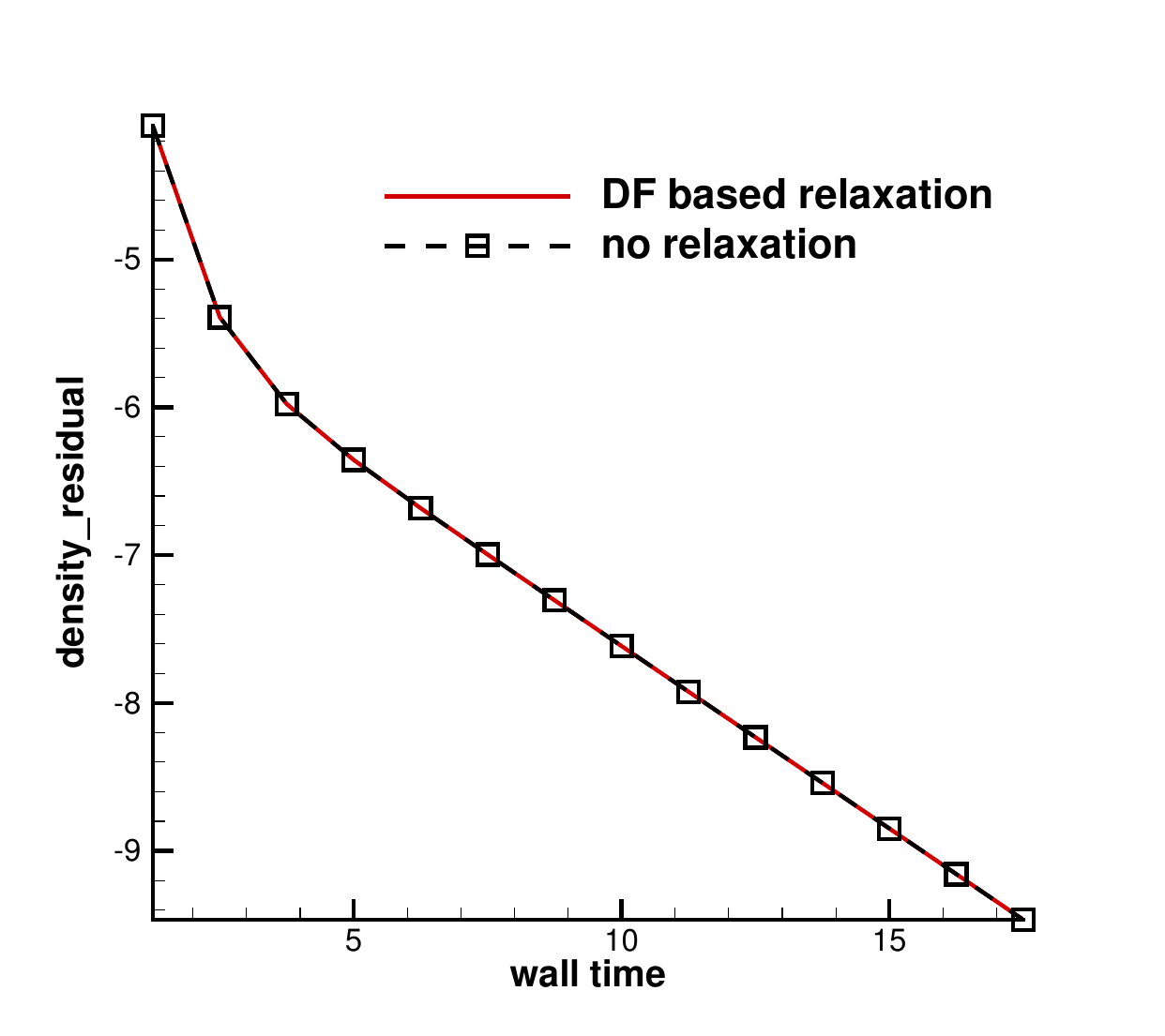}
    \caption{Transonic flow around dual-NACA0012 airfoil. Left: The convergence history of GMG under different LU-SGS sweep numbers (based on a fixed CFL number) Right: The comparison with and without DF-based relaxation.}
    \label{dualNACA-GMG-2}
\end{figure}

Similar to the previous cases, a series of comparisons of convergence rate for this case are shown in Fig .~\ref{dualNACA-GMG-1} and Fig .~\ref{dualNACA-GMG-2}.
Compared to the explicit scheme, GMG's convergence rate is 35 times faster, requiring only 15 seconds to converge. This highlights the significant efficiency of the GPU-accelerated multi-color LU-SGS-based GMG-CGKS.
The comparison of the CFL numbers shows that the GMG method does not require a high CFL number value. The CFL number greater than ten is sufficient. Values too large may prevent the residual from reaching machine zero to achieve a fast convergence rate.
Meanwhile, it is advisable to choose the sweep number for LU-SGS to be greater than two.
In the final comparison of relaxation methods, the convergence rate of the DF-based relaxation aligns with that of the no-relaxation scenario. This consistency arises from the relatively smooth property of the problem, which enables the DF-based relaxation to adaptively revert to the no-relaxation approach.

\subsection{subsonic viscous flow around a sphere}

\begin{figure}[htp]
	\centering	
	\includegraphics[height=0.35\textwidth]{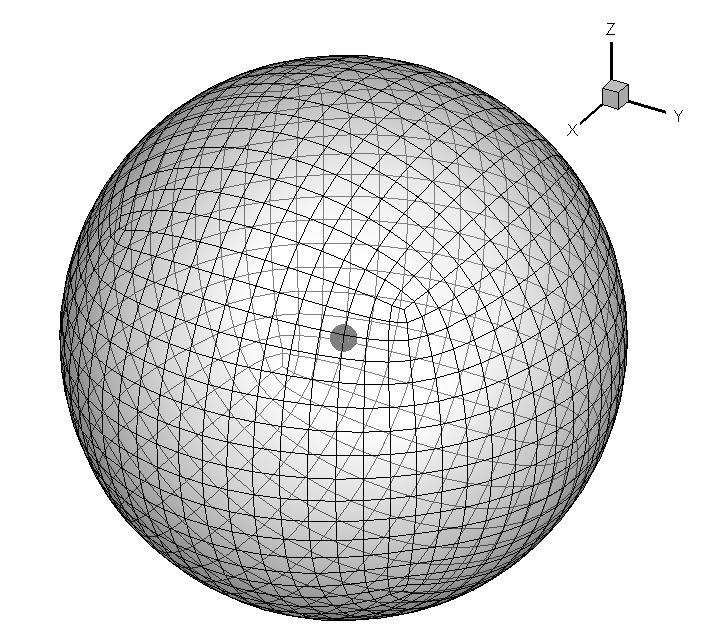}
	\includegraphics[height=0.35\textwidth]{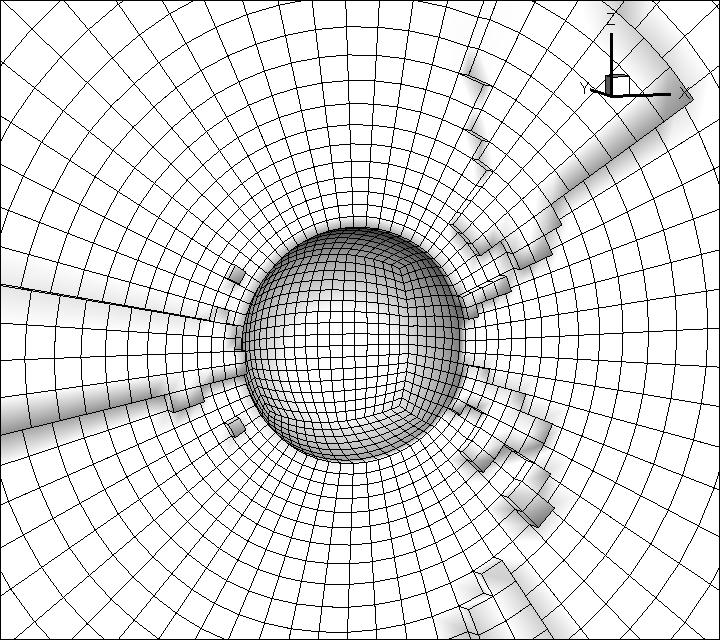}
	\caption{\label{viscous sphere mesh}
		Subsonic flow around a sphere. Mesh sample.}
\end{figure}

In this case, a subsonic flow around a sphere is simulated. The Mach number is set to 0.2535, and the Reynolds number is set to 118.0.
The surface of the sphere is set as non-slip and adiabatic.
The first mesh off the wall has the size $h = 4.5 \times 10^{-2}D$, and the total cell number is 50688.
The mesh is shown in Fig.~\ref{viscous sphere mesh}.

\begin{figure}[htp]	
	\centering	
	\includegraphics[height=0.35\textwidth]{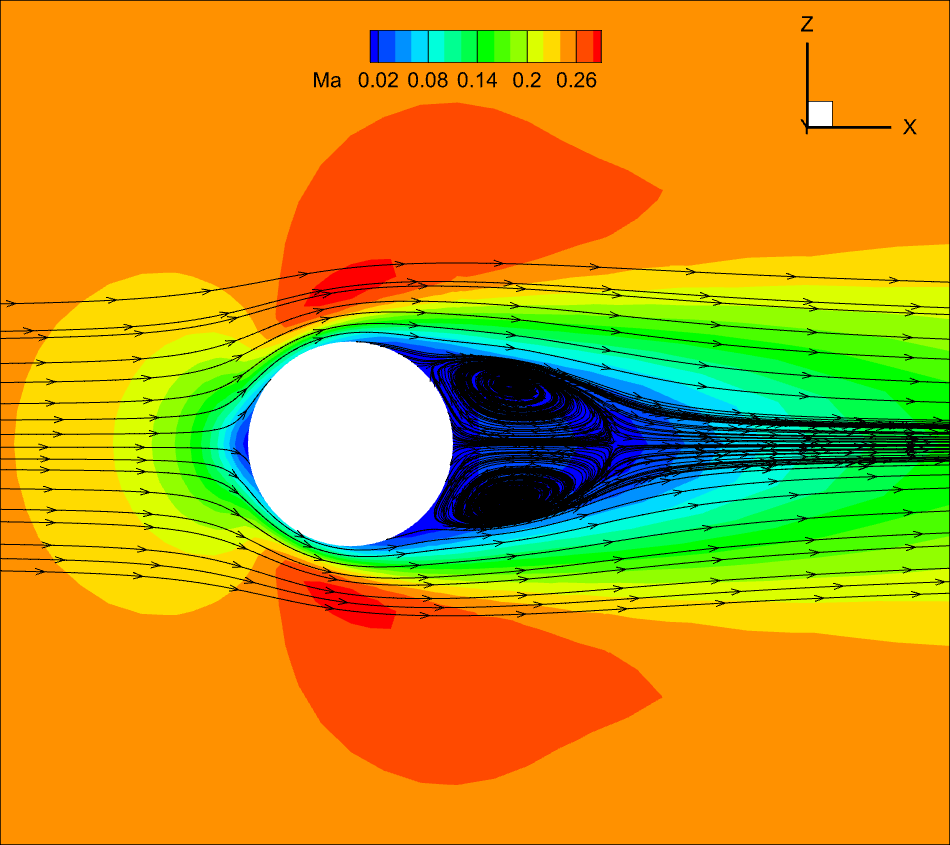}
	\includegraphics[height=0.35\textwidth]{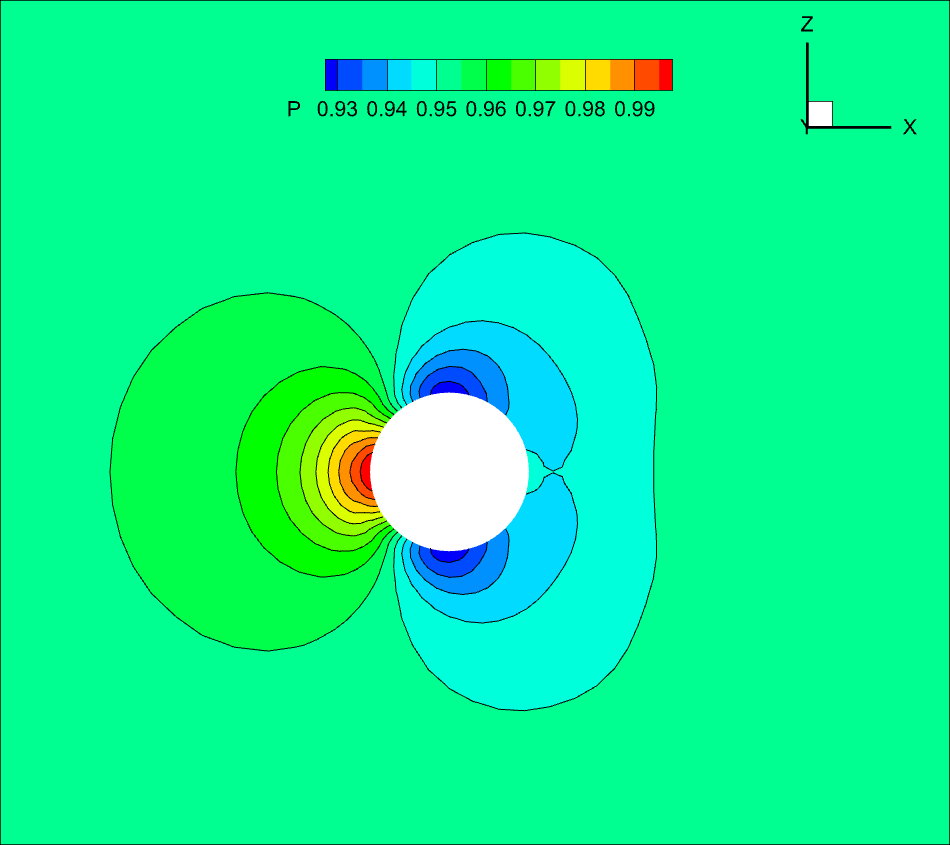}
	\caption{\label{viscous-sphere-density-contour}
		Subsonic flow around a sphere. Left: Mach contour and streamlines. Right: Pressure contour}
\end{figure}

 The Mach number contour and streamline are presented in Fig.~\ref{viscous-sphere-density-contour} to show the high resolution of the GMG-CGKS.
 Quantitative results are given in Table \ref{viscous subsonic sphere}, including the drag coefficient $C_D$, the separation angle $\theta$, and the closed wake length L, as defined in \cite{ji2021compact}.
 The results above show the current scheme agrees well with the experimental and numerical references.

\begin{table}[htp]
	\small
	\begin{center}
		\def\temptablewidth{1.0\textwidth}
		{\rule{\temptablewidth}{1pt}}
		\begin{tabular*}{\temptablewidth}{@{\extracolsep{\fill}}c|c|c|c|c|c}
			Scheme & Mesh number & $C_D$  & $\theta$  &L &Cl\\
			\hline
			Experiment \cite{taneda1956experimental}	&-- & 1.0  & 151 & 1.07 & -- \\ 	
			Third-order DDG \cite{cheng2017parallel} & 160,868 & 1.016 & 123.7 & 0.96 & --\\
			Fourth-order VFV \cite{wang2017thesis}  & 458,915 & 1.014 & --& -- & 2.0e-5\\
			GMG-CGKS & 50688 & 1.020  & 124.5 & 0.89 & 2.26e-5\\
		\end{tabular*}
		{\rule{\temptablewidth}{0.1pt}}
	\end{center}
	\vspace{-4mm} \caption{\label{viscous subsonic sphere} Quantitative comparisons among different compact schemes for the subsonic flow around a sphere.}
\end{table}

\begin{figure}
    \centering
    \includegraphics[width=0.35\linewidth]{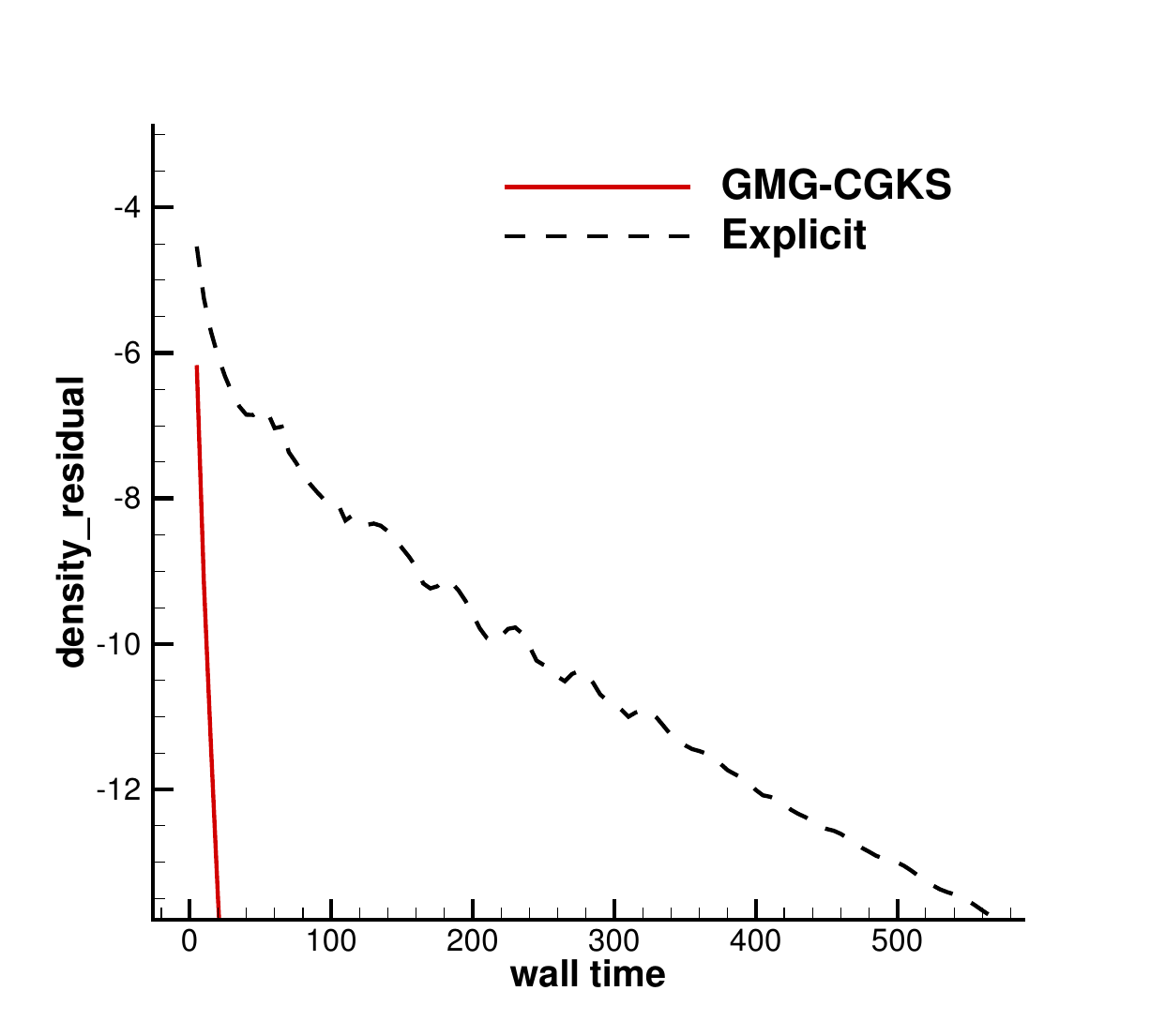}
    \includegraphics[width=0.35\linewidth]{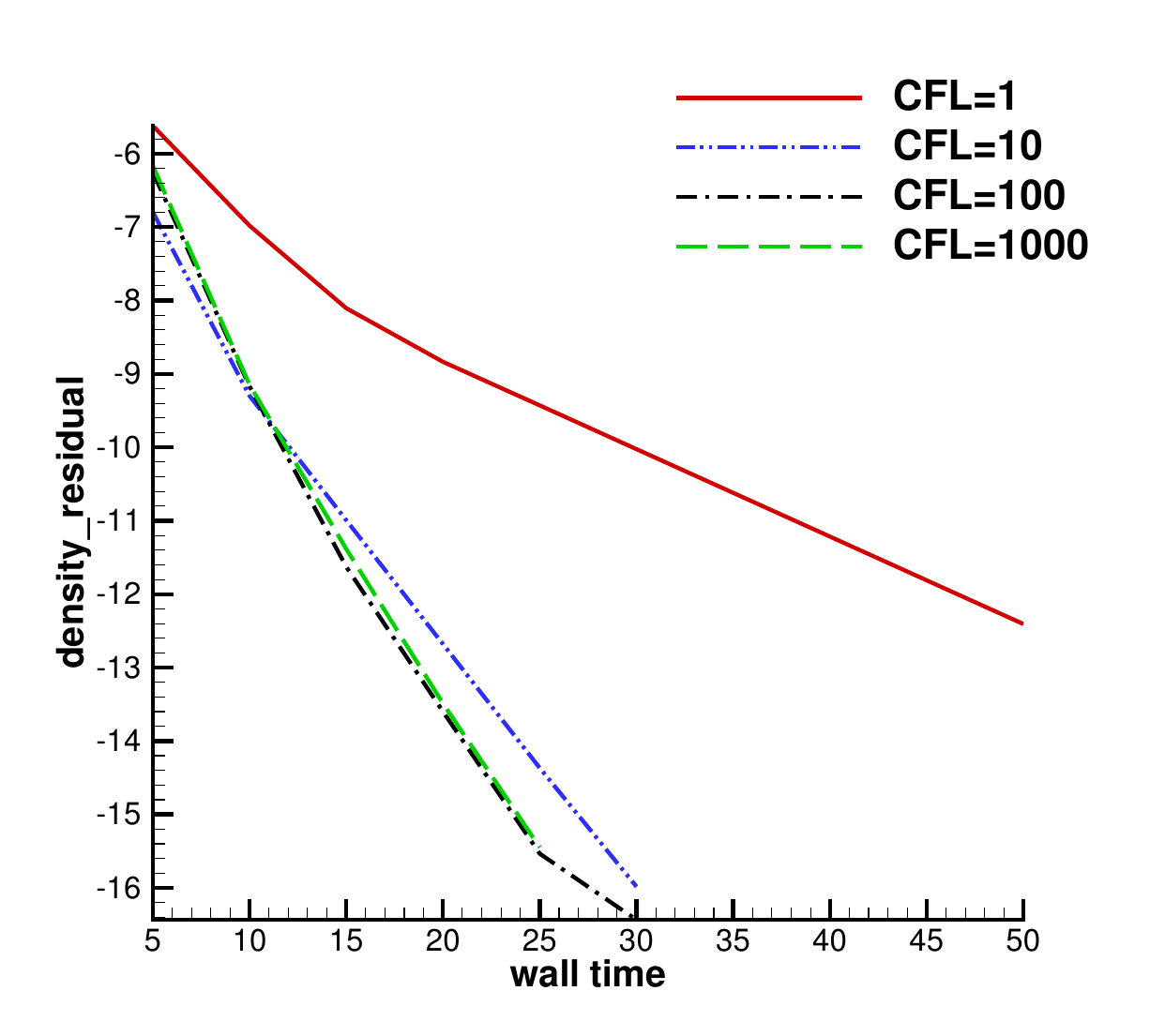}
    \caption{Subsonic flow around a sphere. Left: The convergence history of GMG and explicit CGKS. Right: The convergence history of GMG under different CFL numbers (based on a fixed sweep number of six of LU-SGS).}
    \label{sphere02535-GMG-1}
\end{figure}

\begin{figure}
    \centering
    \includegraphics[width=0.35\linewidth]{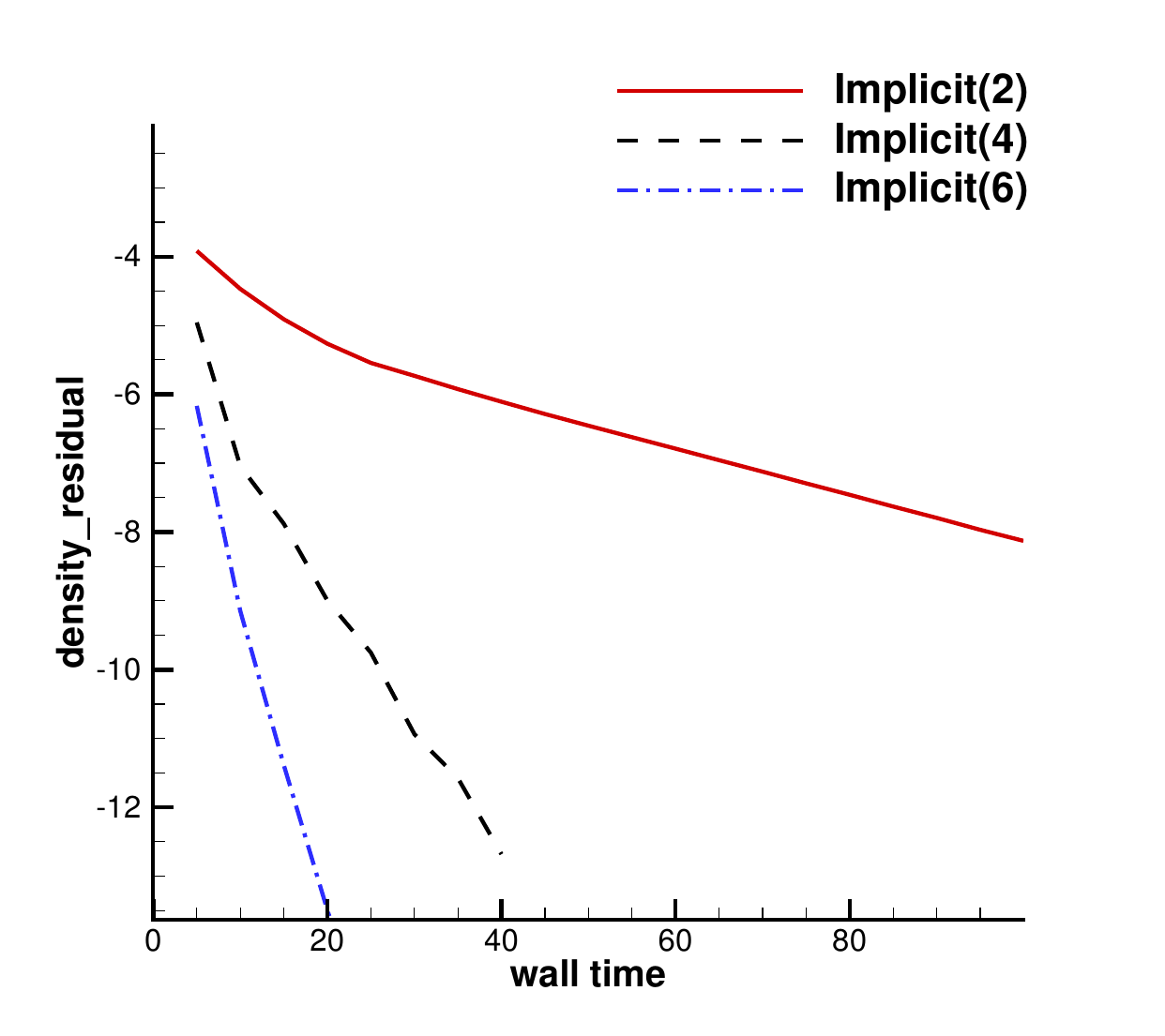}
    \includegraphics[width=0.35\linewidth]{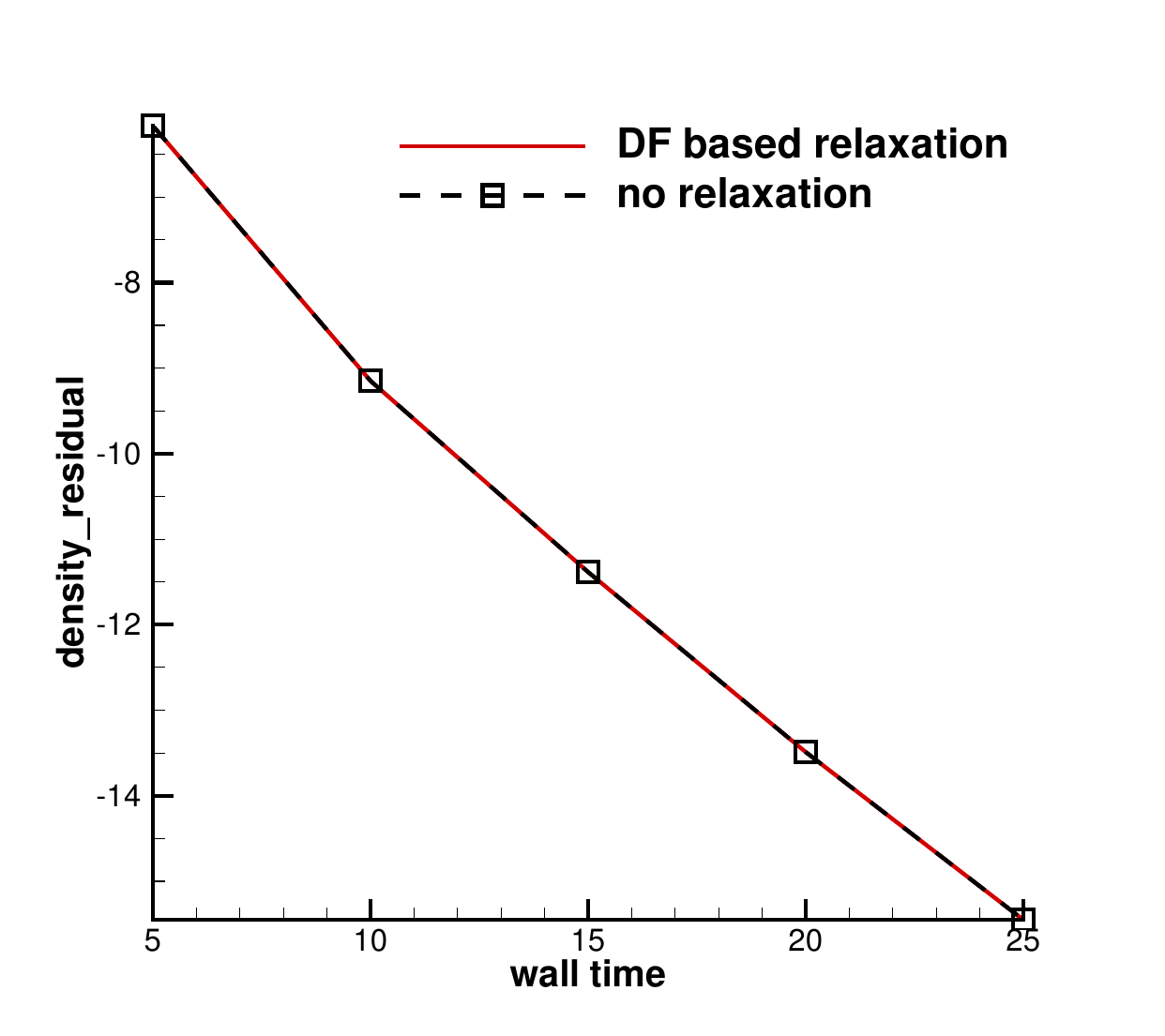}
    \caption{Subsonic flow around a sphere. Left: The convergence history of GMG under different LU-SGS sweep numbers (based on a fixed CFL number) Right: The comparison with and without DF-based relaxation.}
    \label{sphere03535-GMG-2}
\end{figure}

A series of comparisons of convergence rate for this case are shown in Fig .~\ref{sphere02535-GMG-1} and Fig .~\ref{sphere03535-GMG-2}, including a comparison between explicit and GMG, a comparison with a fixed CFL number while varying the LU-SGS sweep number in the coarse grid, a comparison with a fixed LU-SGS sweep number while varying the CFL number, and a comparison with and without DF-based relaxation.
Compared to the explicit scheme, the GMG method achieves a convergence rate that is 28 times faster, requiring only 20 seconds to converge. This is particularly noteworthy given the already high computational speed of the GPU, underscoring the significant efficiency of the GPU-accelerated multi-color LU-SGS-based GMG-CGKS.
Regarding comparing CFL numbers, the GMG method does not necessitate a high CFL value; a value greater than ten is sufficient to achieve a fast convergence rate.
Additionally, the comparison results indicate that the sweep count for LU-SGS should not be too low, with four to six sweeps being optimal.
In the final comparison of relaxation methods, the convergence rate of the DF-based relaxation is consistent with that of the no-relaxation case. This consistency is due to the relatively smooth nature of the problem, which allows the DF-based relaxation to adaptively revert to the no-relaxation scheme.

\subsection{Transonic viscous flow around a sphere}
\begin{figure}[htp]	
	\centering	
	\includegraphics[height=0.35\textwidth]{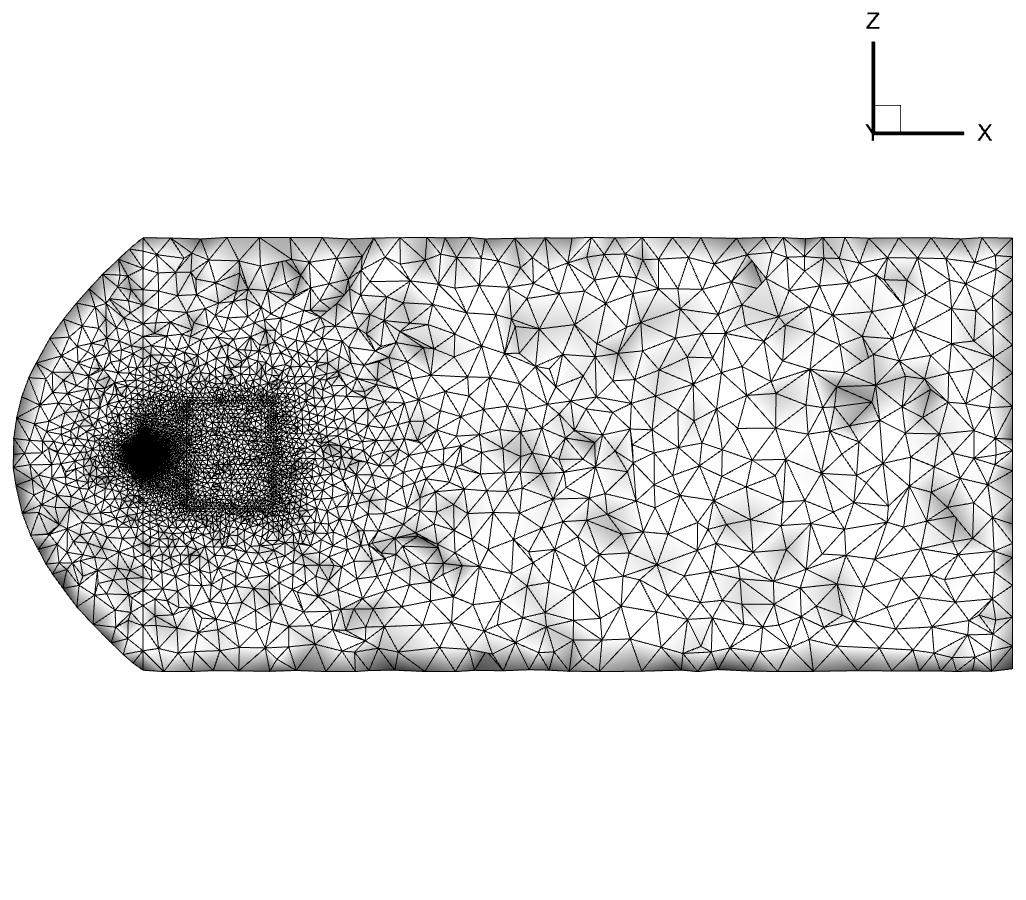}
	\includegraphics[height=0.35\textwidth]{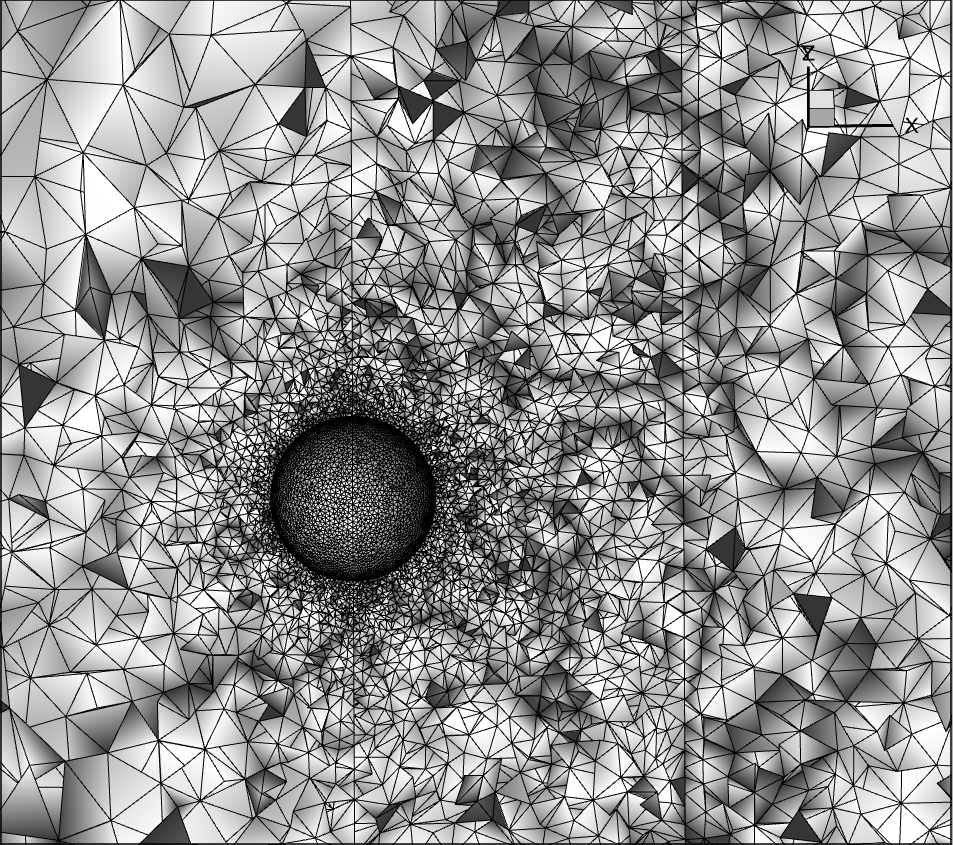}
	\caption{\label{tetrahedron-sphere-mesh}
		The mesh of Transonic flow around a sphere. Left: Global mesh. Right: Local mesh. }
\end{figure}

\begin{figure}[htp]	
	\centering	
	\includegraphics[height=0.35\textwidth]{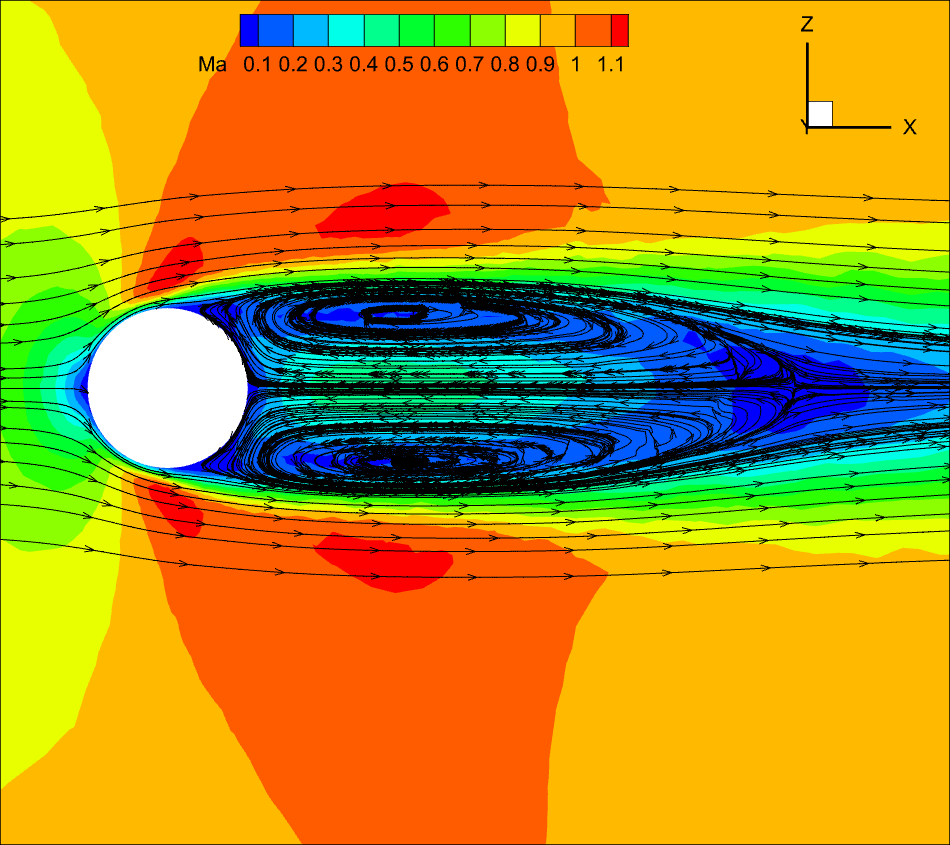}
	\includegraphics[height=0.35\textwidth]{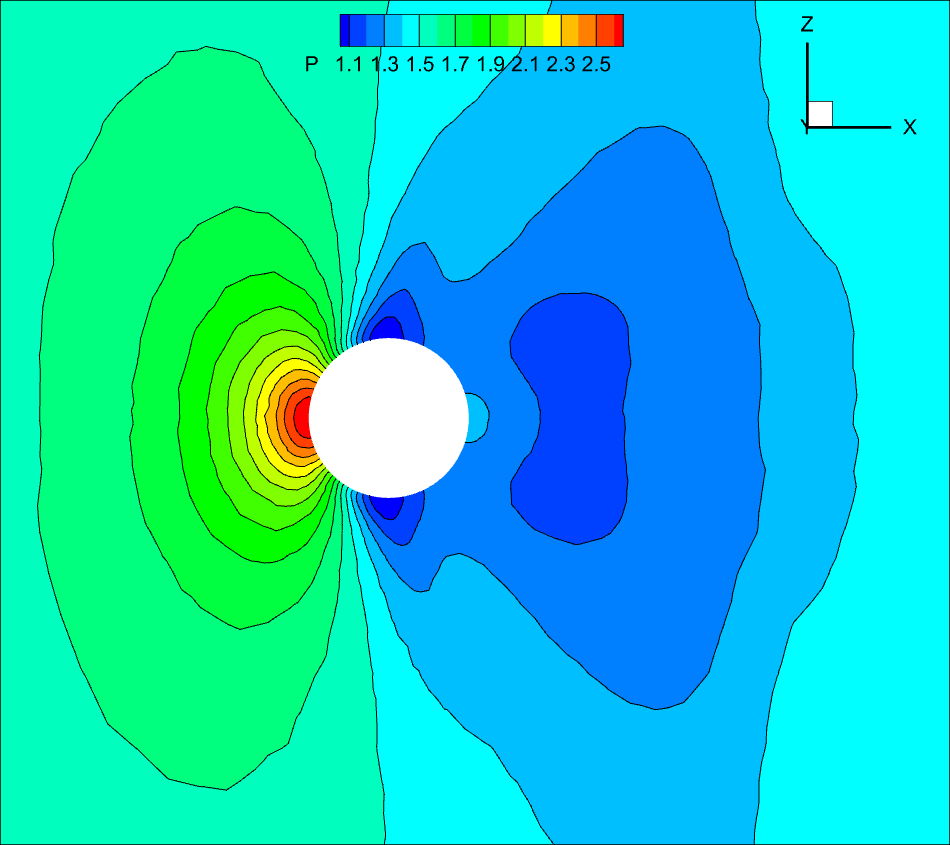}
	\caption{\label{viscous-sphere-ma0.95-contour}
		Transonic flow around a sphere. Left: Mach number contour with streamline through the sphere. Right: Pressure contour. }
\end{figure}

A transonic viscous flow around a sphere is simulated to show the performance of the GMG-CGKS for transonic viscous flow.
The Mach number is set to be 0.95, and the Reynolds number is set to be 300.0.
In this case, we use the pure tetrahedron mesh with a mesh number equal to 665914, and the wake part of the sphere is refined to capture the vortex.
The mesh used in this case is shown in Fig.~\ref{tetrahedron-sphere-mesh}.
The numerical results of the Mach number contour and streamline around a sphere are shown in Fig.~\ref{viscous-sphere-ma0.95-contour}, which indicates the high resolution of the GMG-CGKS.
Quantitative results include the drag coefficient $C_D$,  the wake length $L$, and the separation angle $\theta$ are listed in Table \ref{transsonic-sphere}.
The results above show that the current scheme agrees well with the numerical references, even when using the higher order.

\begin{table}[htp]
	\small
	\begin{center}
		\def\temptablewidth{1.0\textwidth}
		{\rule{\temptablewidth}{1pt}}
		\begin{tabular*}{\temptablewidth}{@{\extracolsep{\fill}}c|c|c|c|c}
			Scheme & Mesh Number & $C_D$  & $\theta$  &L\\
			\hline
			WENO6 \cite{Nagata2016sphere} 	&909,072 & 0.968  & 111.5 & 3.48\\ 	
			Original CGKS \cite{ji2021compact}  & 515,453 & 0.950  & 112.7 & 3.30\\
   GMG-CGKS  & 665,914 & 0.960  & 111.6 & 3.49\\
		\end{tabular*}
		{\rule{\temptablewidth}{0.1pt}}
	\end{center}
	\vspace{-4mm} \caption{\label{transsonic-sphere} Quantitative comparisons between the current scheme and the reference solution for the transonic flow around a sphere.}
\end{table}

\begin{figure}
    \centering
    \includegraphics[width=0.35\linewidth]{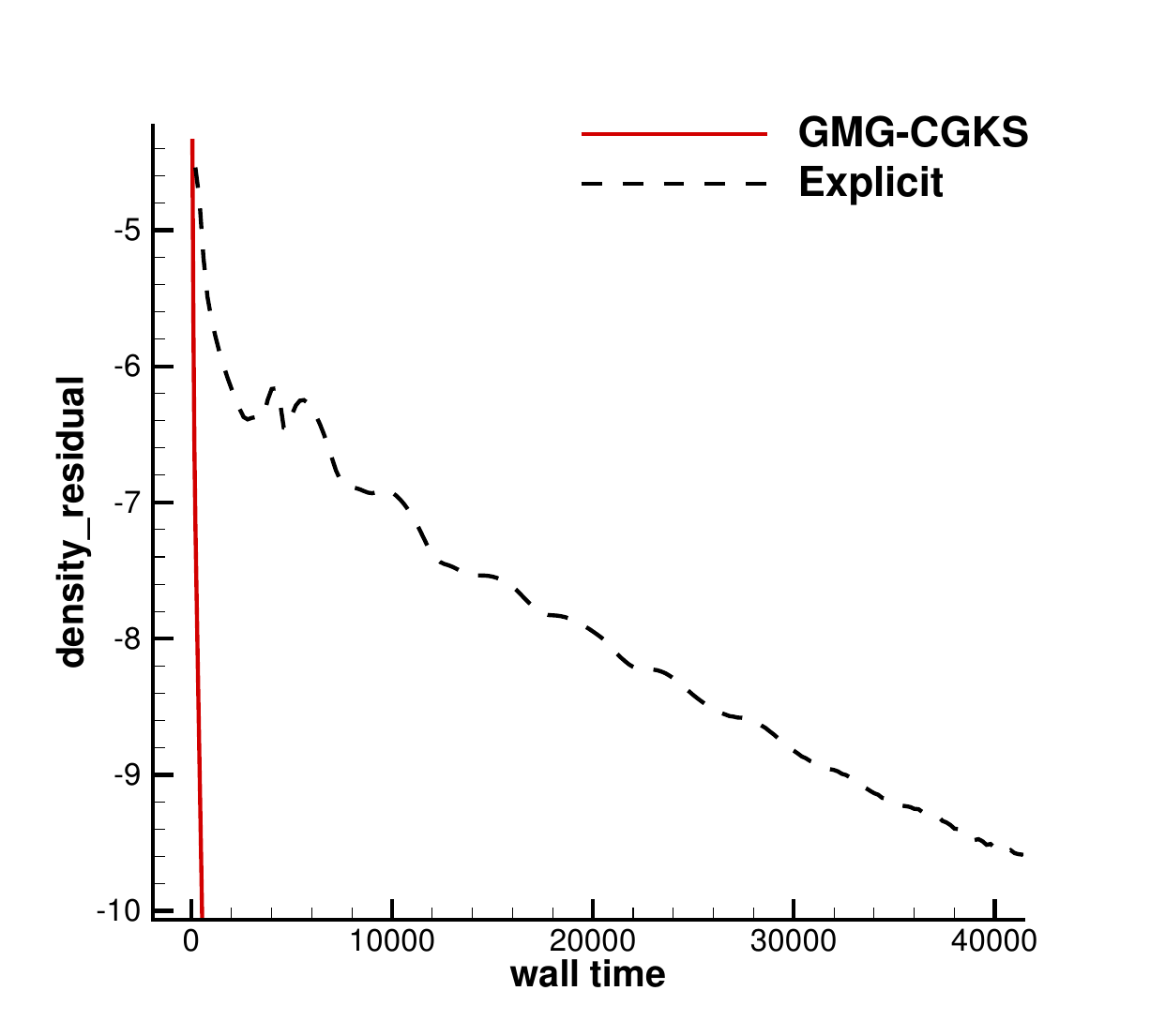}
    \includegraphics[width=0.35\linewidth]{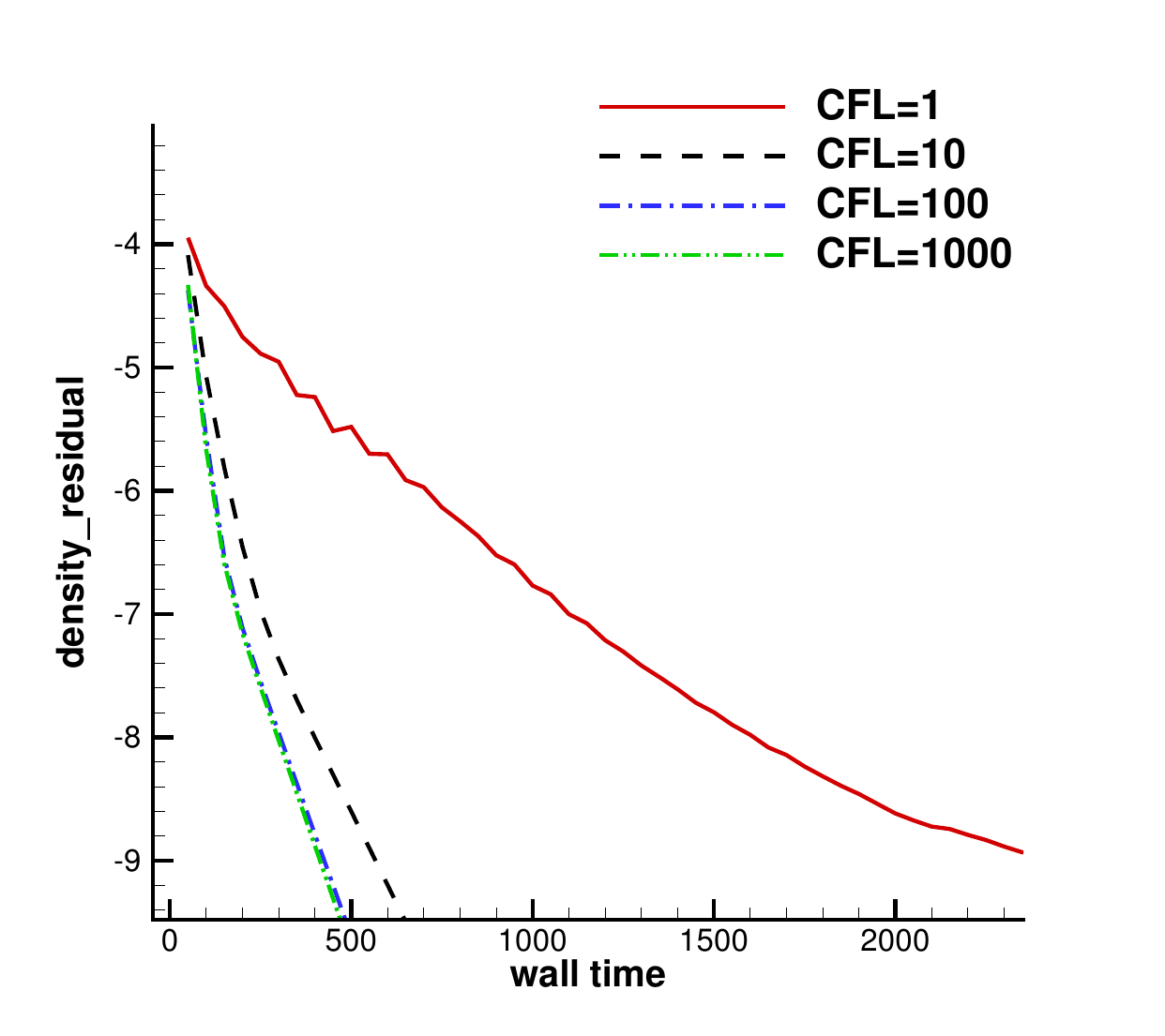}
    \caption{Transonic flow around a sphere. Left: The convergence history of GMG and explicit CGKS. Right: The convergence history of GMG under different CFL numbers (based on a fixed sweep number of six of LU-SGS).}
    \label{sphere095-GMG-1}
\end{figure}

\begin{figure}
    \centering
    \includegraphics[width=0.35\linewidth]{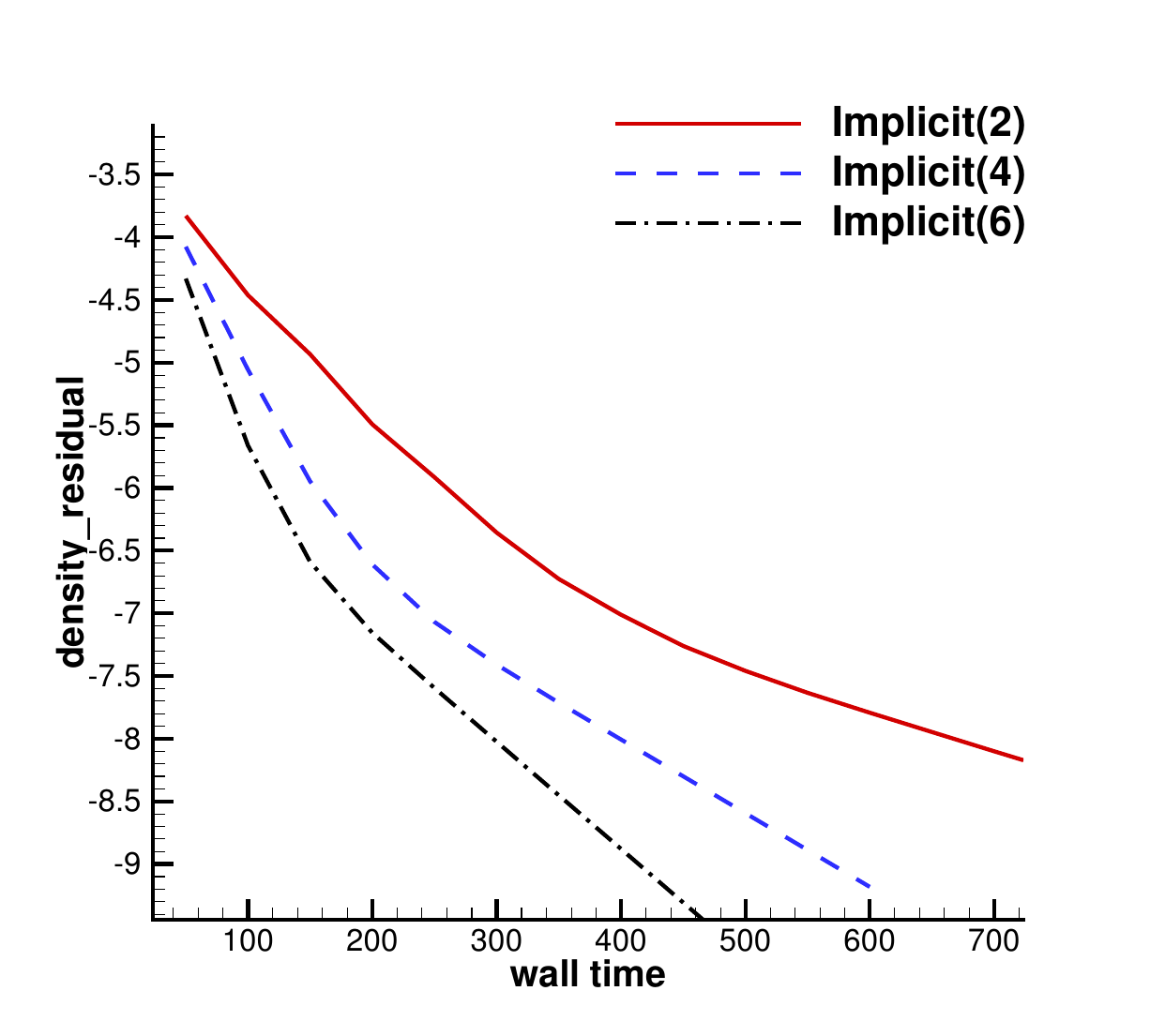}
    \includegraphics[width=0.35\linewidth]{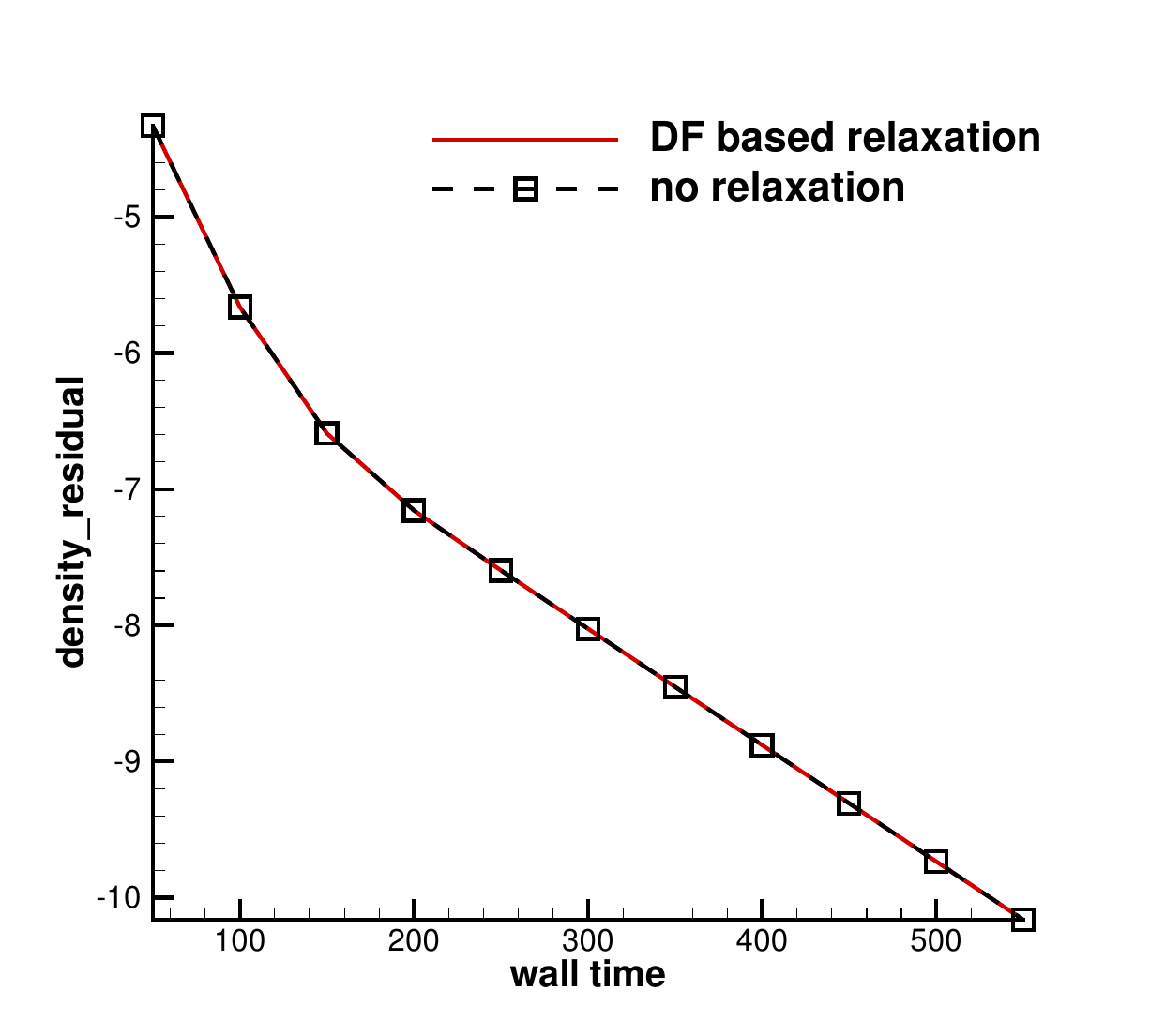}
    \caption{Transonic flow around a sphere. Left: The convergence history of GMG under different LU-SGS sweep numbers (based on a fixed CFL number) Right: The comparison with and without DF-based relaxation.}
    \label{sphere095-GMG-2}
\end{figure}

A series of comparisons of convergence rate for this case are shown in Fig .~\ref{sphere095-GMG-1} and Fig .~\ref{sphere095-GMG-2}, including a comparison between explicit and GMG, a comparison with a fixed CFL number while varying the LU-SGS sweep number in the coarse grid, a comparison with a fixed LU-SGS sweep number while varying the CFL number, and a comparison with and without DF-based relaxation.
Compared to the explicit scheme, GMG's convergence rate is 90 times faster, requiring only 8 minutes to converge. This is particularly notable given the GPU's already high computational speed, highlighting the significant efficiency of the GPU-accelerated multi-color LU-SGS-based GMG-CGKS.
Regarding comparing CFL numbers, the GMG method does not require a high CFL value (greater than ten is enough) to achieve a fast convergence rate.
At the same time, the comparison results indicate that the sweep count for LU-SGS should not be too low, with four to six sweeps being optimal.
In the final comparison regarding relaxation, the convergence rate of DF-based relaxation is consistent with that of the no-relaxation case. This is because the case is relatively smooth, allowing the DF-based relaxation to adaptively revert to the no-relaxation scheme.

\subsection{supersonic viscous flow around a sphere}

\begin{figure}[htp]	
	\centering	
	\includegraphics[height=0.35\textwidth]{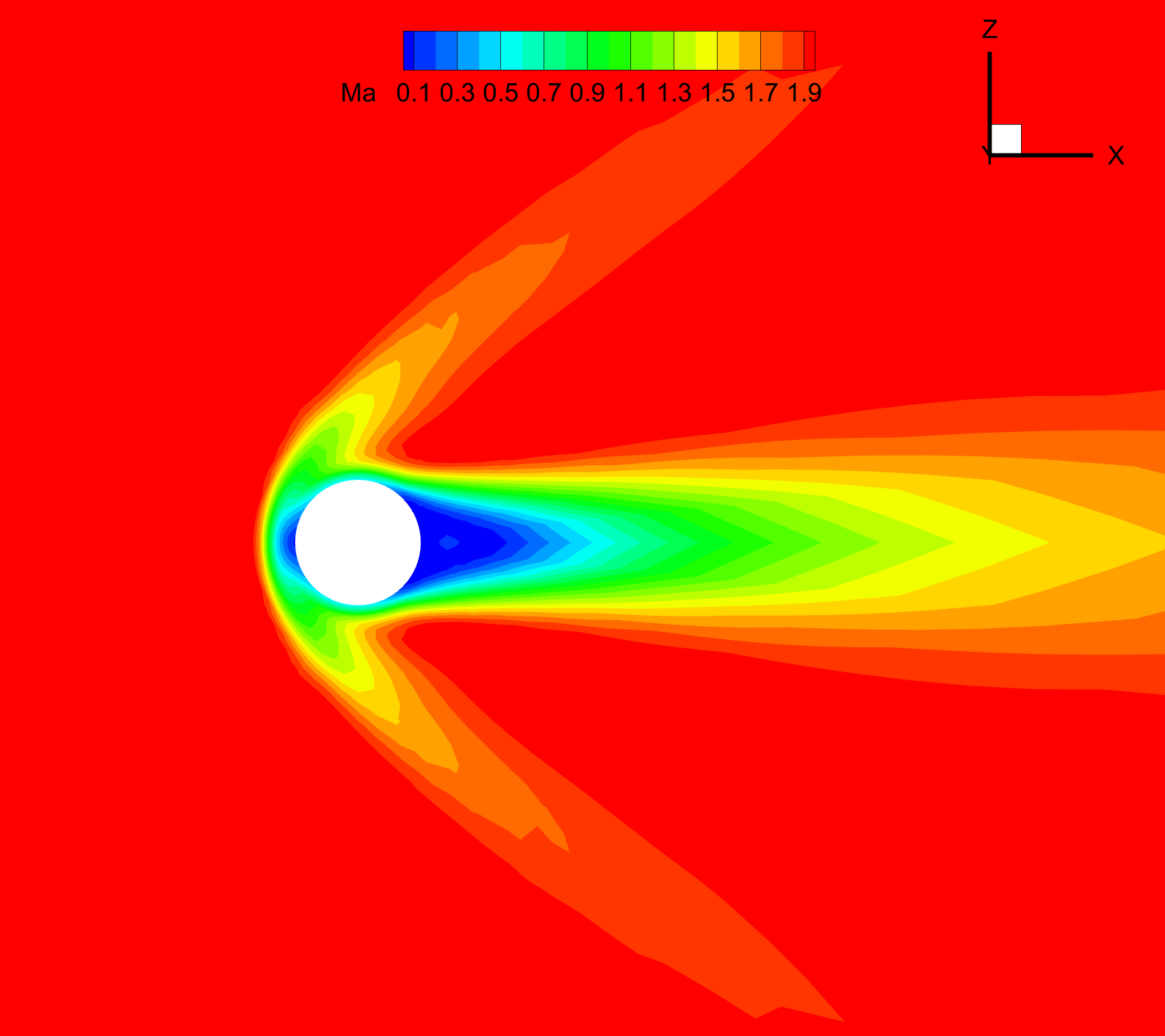}
	\includegraphics[height=0.35\textwidth]{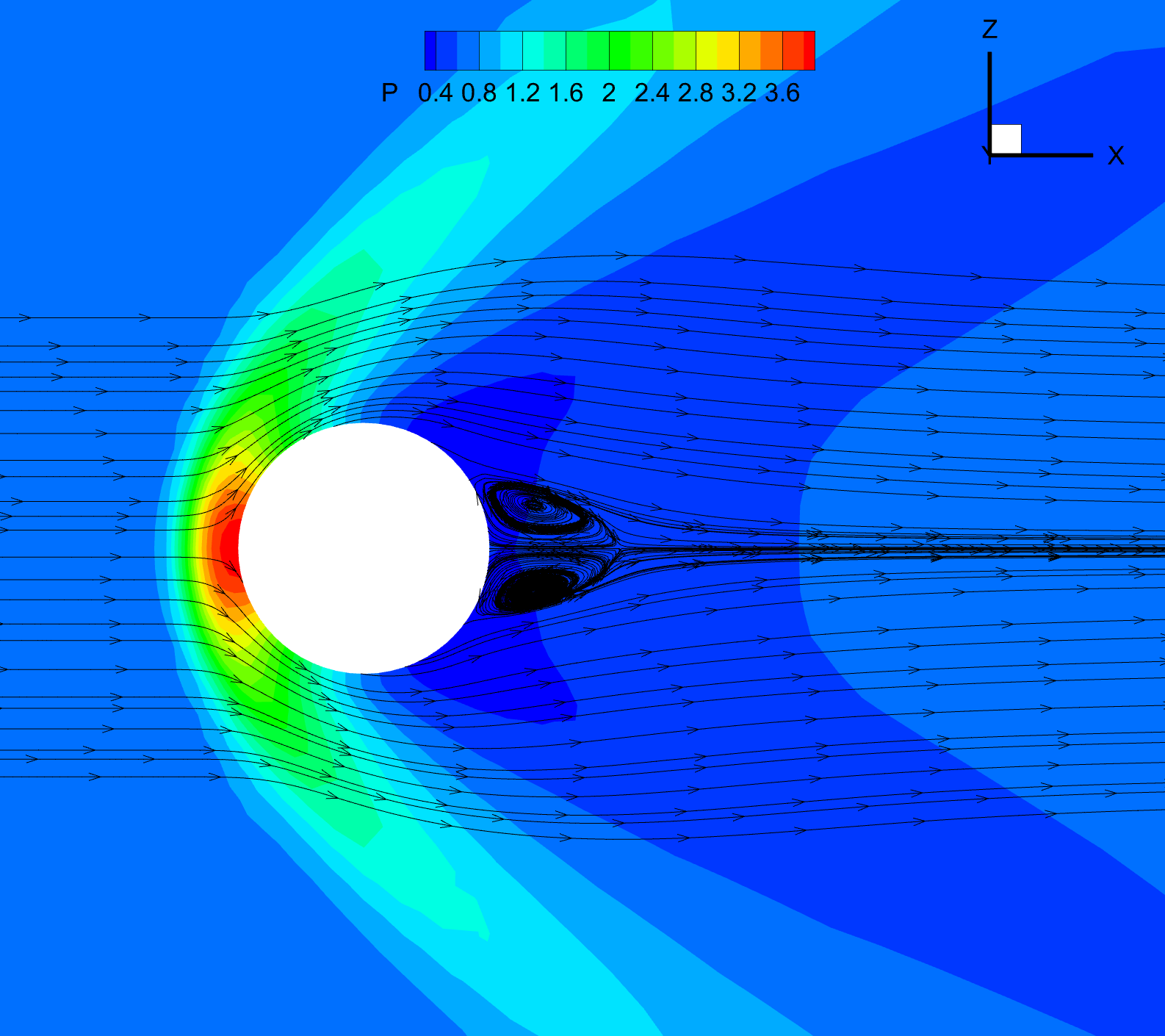}
	\caption{\label{viscous-sphere-ma2-contour}
		Supersonic flow around a sphere. Left: Mach number contour with streamline through the sphere. Right: Pressure contour. }
\end{figure}

A supersonic flow around a sphere is simulated to verify that the GMG-CGKS can perform well in the supersonic flow region.
The Mach number is set to be 2.0, and the Reynolds number is set to be 300.
The mesh used in this case is the same as the subsonic case.
The first layer mesh at the wall has a thickness $2.3 \times 10^{-2}D$.
Mach number and pressure distribution are shown in Fig .~\ref{viscous-sphere-ma2-contour}.
Quantitative results are listed in Table \ref{supersonic-sphere}, which agrees well with those given by Ref.\cite{Nagata2016sphere}.

\begin{table}[htp]
	\small
	\begin{center}
		\def\temptablewidth{1.0\textwidth}
		{\rule{\temptablewidth}{1pt}}
		\begin{tabular*}{\temptablewidth}{@{\extracolsep{\fill}}c|c|c|c|c|c}
			Scheme & Mesh Number & Cd  & $\theta$  &L & Shock stand-off\\
			\hline
			WENO6 \cite{Nagata2016sphere} 	&909,072 & 1.386  & 150.9 & 0.38 & 0.21 \\ 	
			Original CGKS \cite{JI2024112590}  & 50688 & 1.368  & 148.5 & 0.45 & 0.28 \\
                GMG-CGKS  & 50688 & 1.368  & 148.5 & 0.45 & 0.28 \\
		\end{tabular*}
		{\rule{\temptablewidth}{0.1pt}}
	\end{center}
	\vspace{-4mm} \caption{\label{supersonic-sphere} Quantitative comparisons between the current scheme and the reference solution for the supersonic flow around a sphere.}
\end{table}

\begin{figure}
    \centering
    \includegraphics[width=0.35\linewidth]{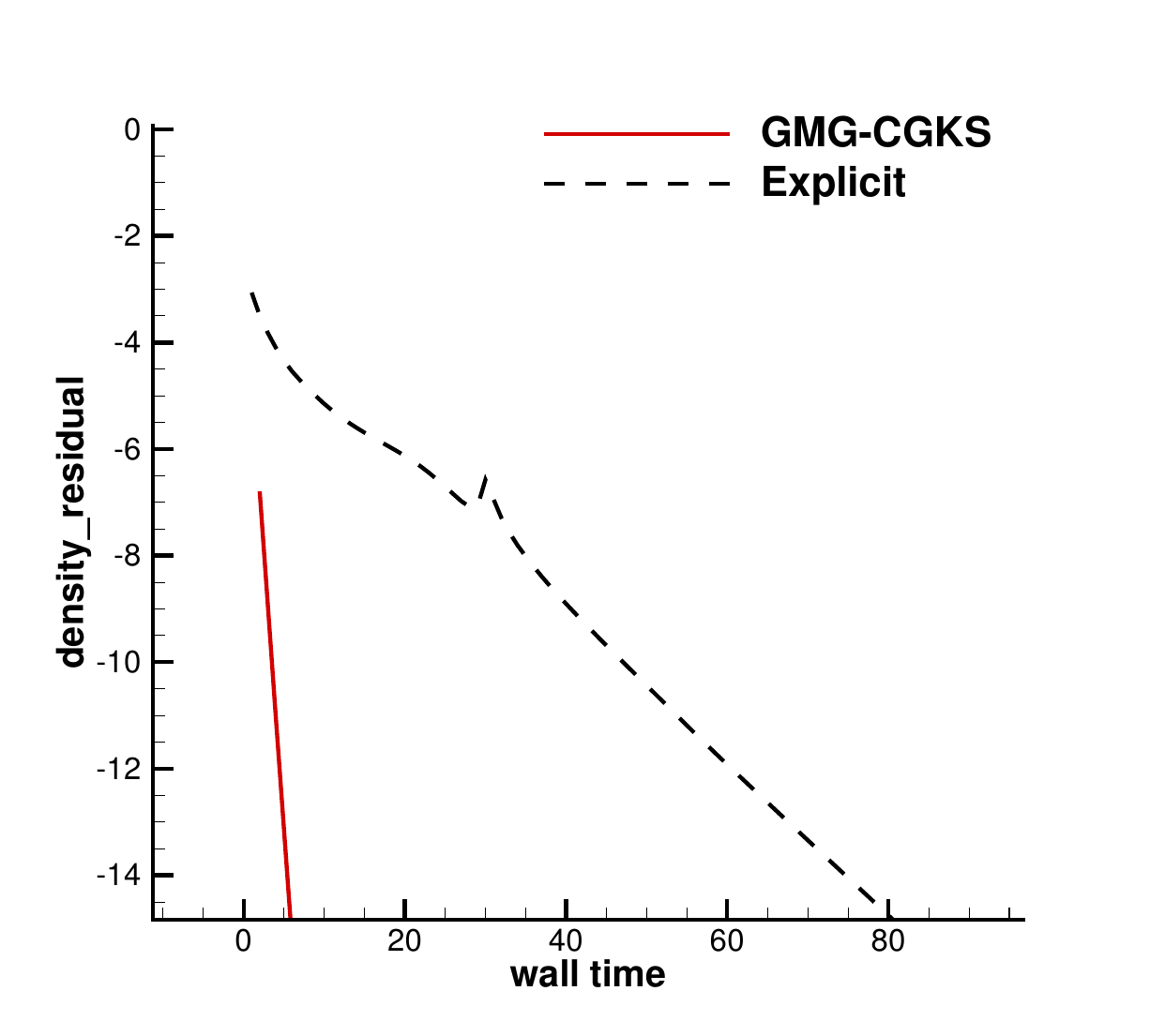}
    \includegraphics[width=0.35\linewidth]{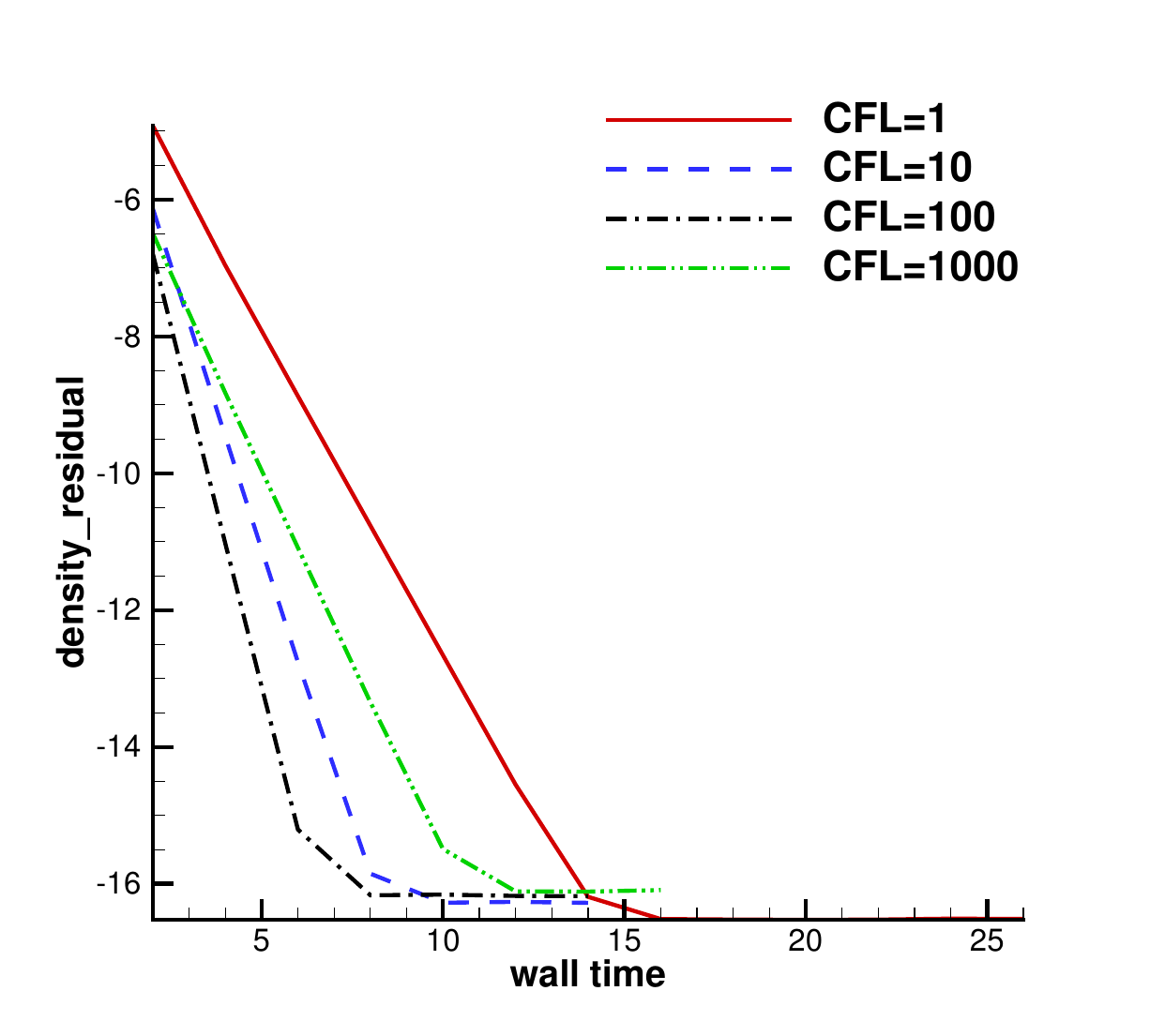}
    \caption{Supersonic flow around a sphere. Left: The convergence history of GMG and explicit CGKS. Right: The convergence history of GMG under different CFL numbers (based on a fixed sweep number of six of LU-SGS).}
    \label{sphere2-GMG-1}
\end{figure}

\begin{figure}
    \centering
    \includegraphics[width=0.35\linewidth]{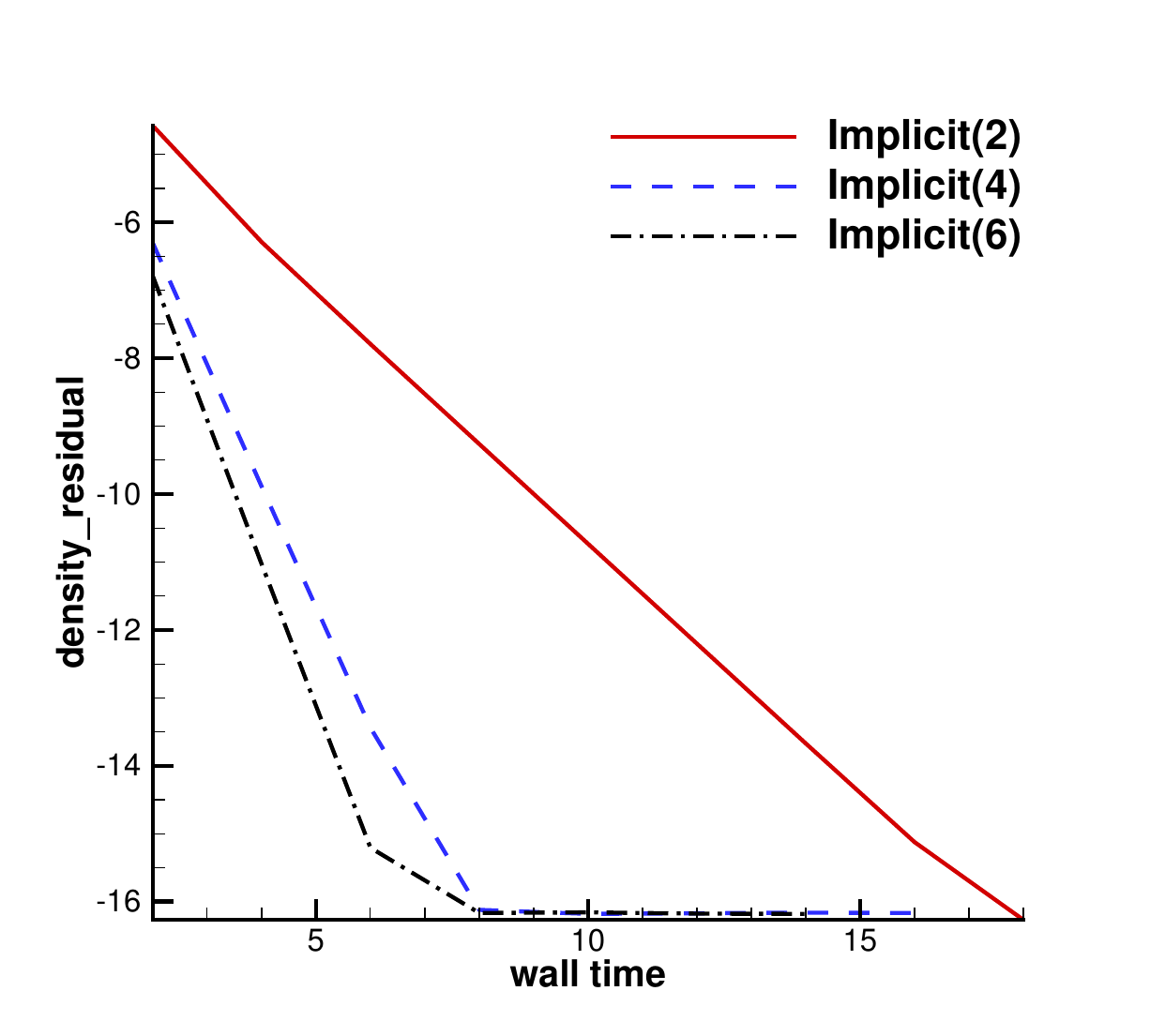}
    \caption{Supersonic flow around a sphere. The convergence history of GMG under different LU-SGS sweep numbers (based on a fixed CFL number) }
    \label{sphere2-GMG-2}
\end{figure}

A series of comparisons of convergence rate for this case are shown in Fig .~\ref{sphere2-GMG-1} and Fig .~\ref{sphere2-GMG-2}, including a comparison between explicit and GMG, a comparison with a fixed CFL number while varying the LU-SGS sweep number in the coarse grid, a comparison with a fixed LU-SGS sweep number while varying the CFL number, and a comparison with and without DF-based relaxation.
Compared to the explicit scheme, GMG's convergence rate is 13 times faster, requiring only 5 seconds to converge. This is particularly notable given the GPU's already high computational speed, highlighting the significant efficiency of the GPU-accelerated multi-color LU-SGS-based GMG-CGKS.
Regarding comparing CFL numbers, the GMG method does not require a high CFL value (greater than ten is enough) to achieve a fast convergence rate.
At the same time, the comparison results indicate that the sweep count for LU-SGS should not be too low, with four to six sweeps being optimal.
Furthermore, if DF-based relaxation is not used, this case will explode, demonstrating the importance of DF-based relaxation.

\subsection{M6-wing}
\begin{figure*}[htp]	
	\centering	
	\includegraphics[height=0.35\textwidth]{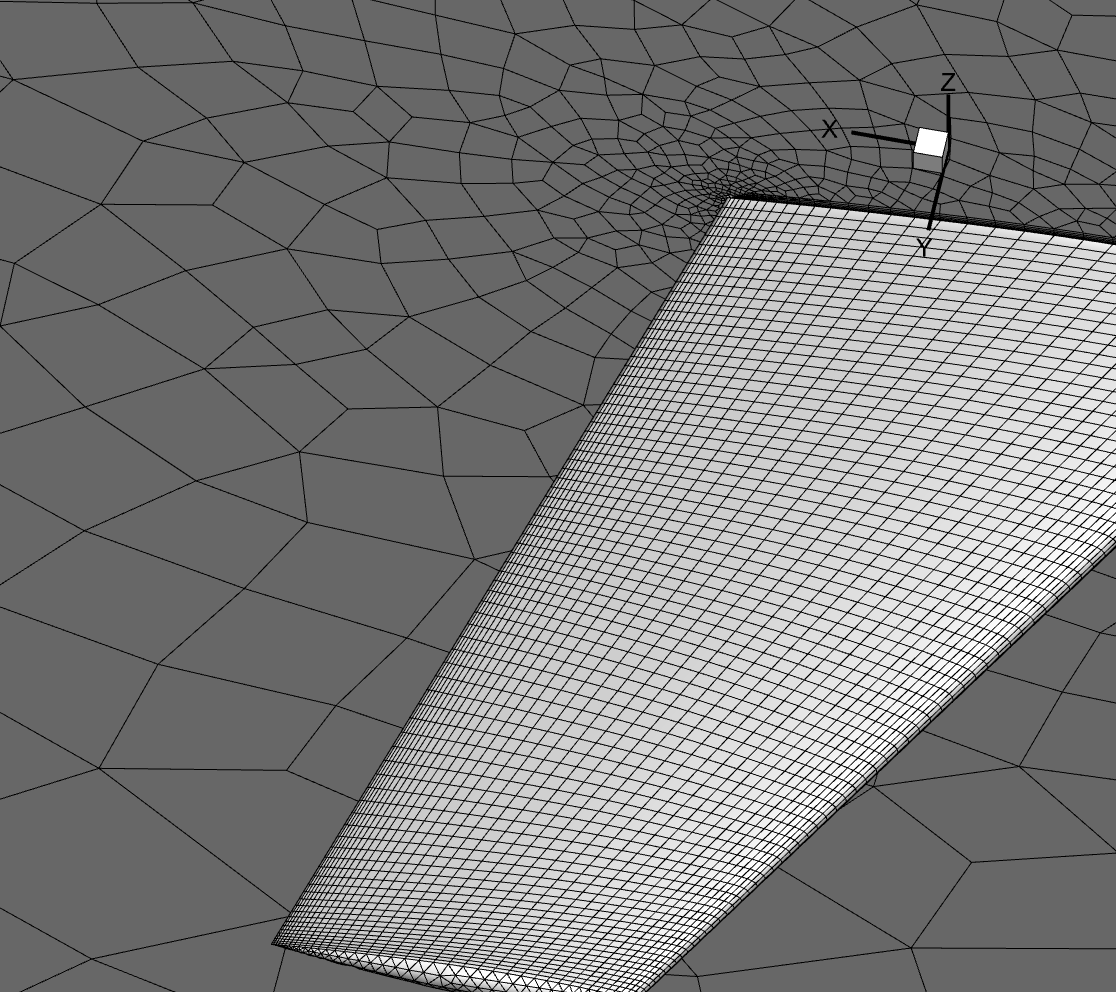}
	\includegraphics[height=0.35\textwidth]{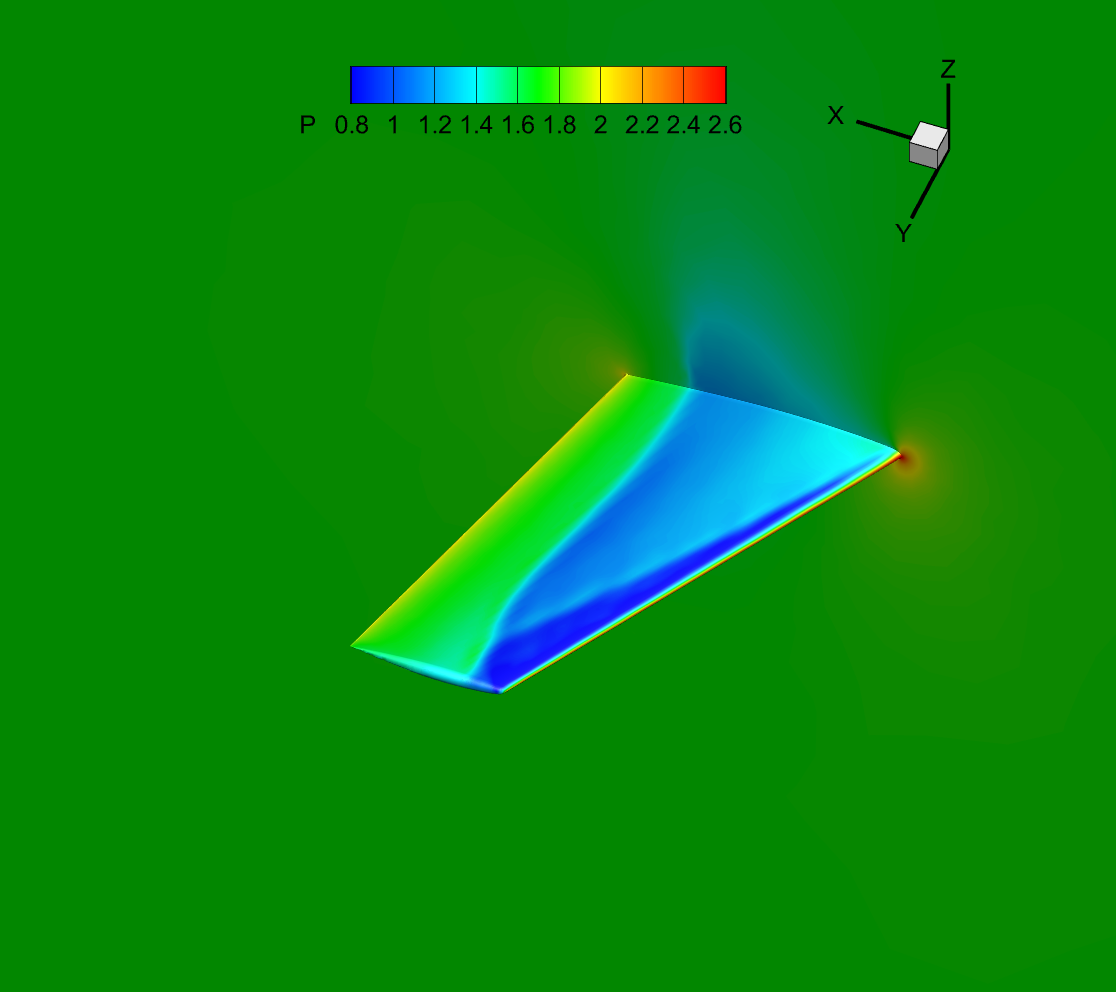}
	\caption{\label{M6 Contour}
		Transonic flow around a M6 wing. Left: Mesh. Right: Pressure contour.}
\end{figure*}

\begin{figure}[htp]	
	\centering	
	\includegraphics[height=0.35\textwidth]{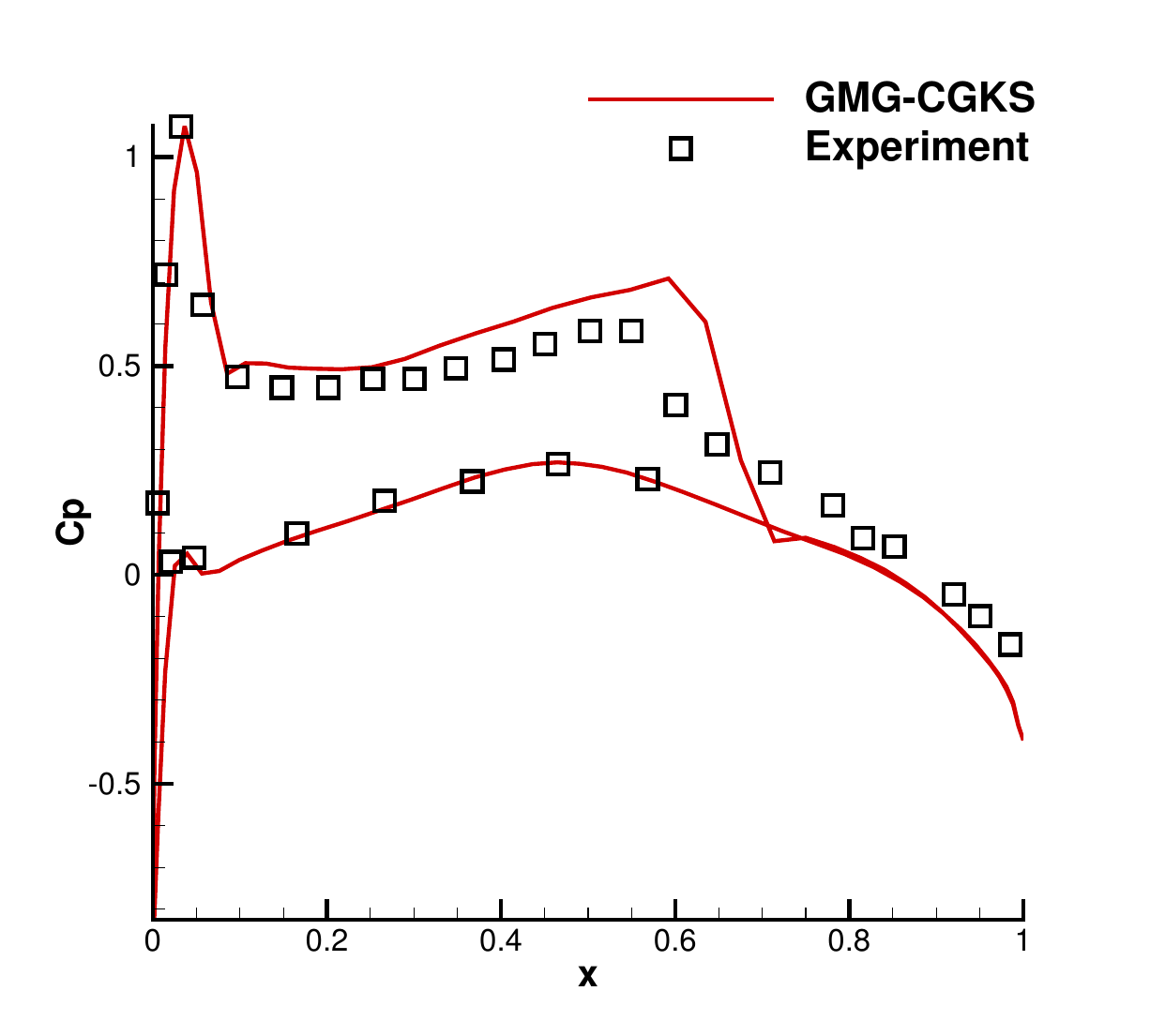}
	\includegraphics[height=0.35\textwidth]{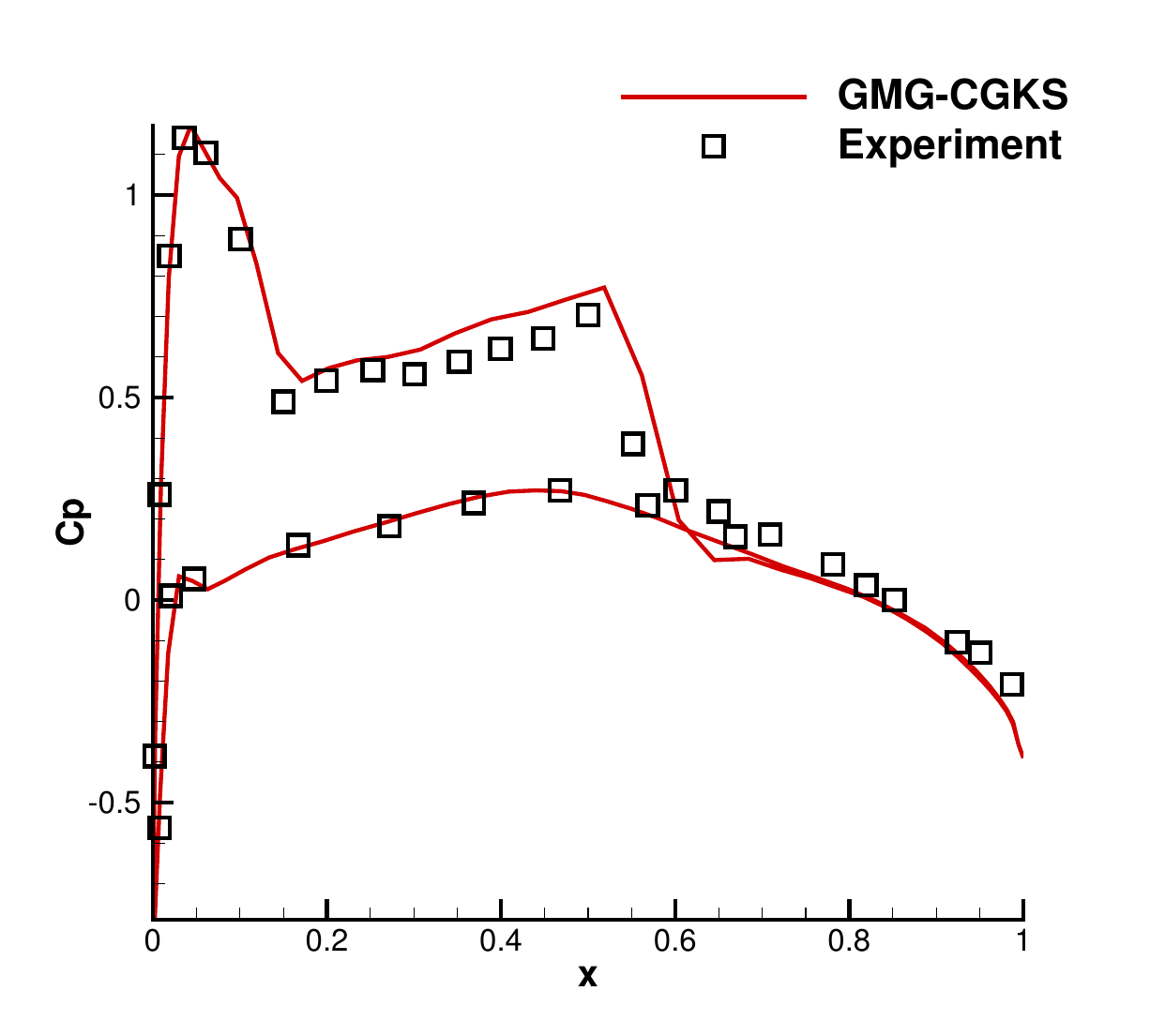}
	\includegraphics[height=0.35\textwidth]{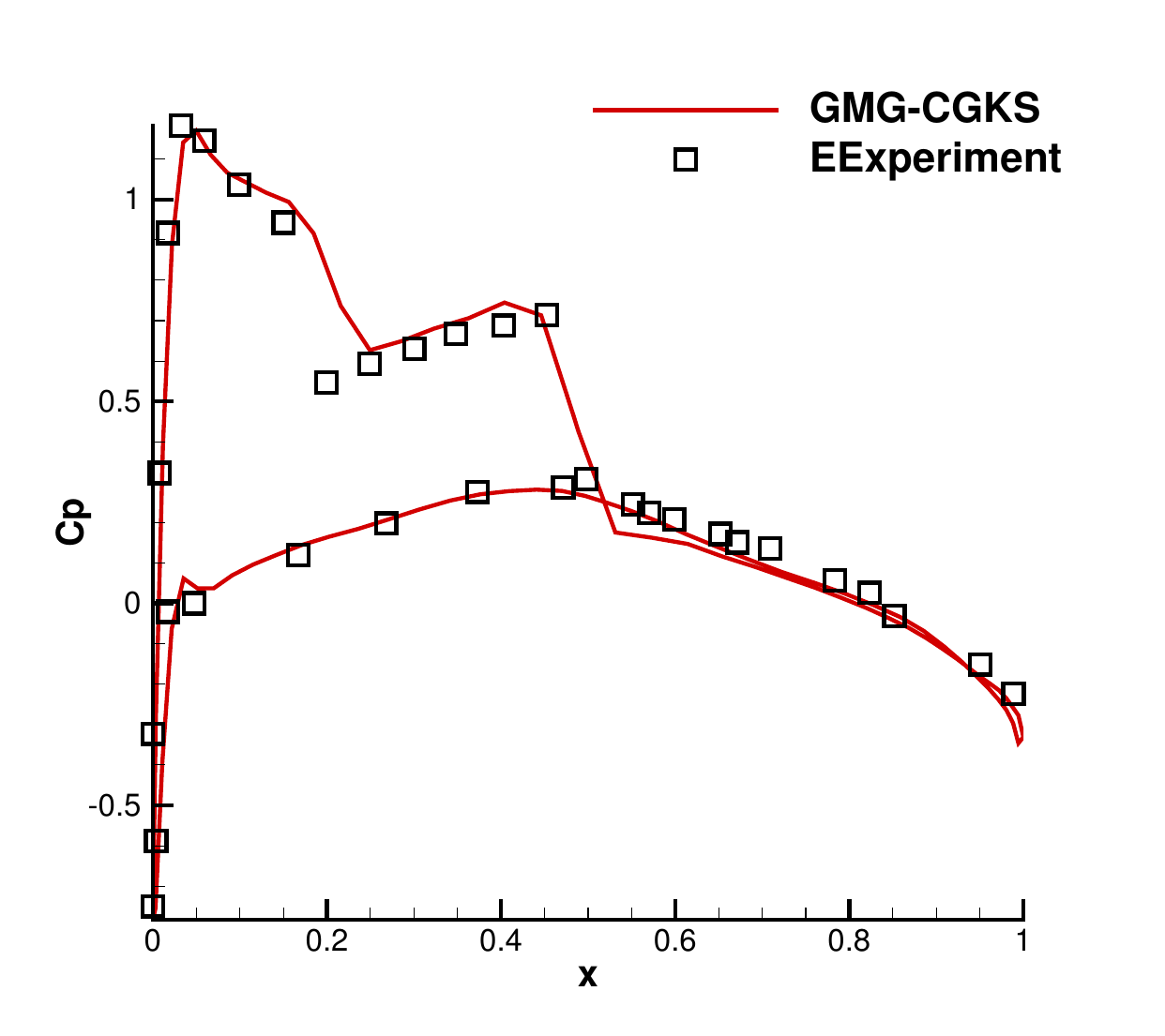}
	\includegraphics[height=0.35\textwidth]{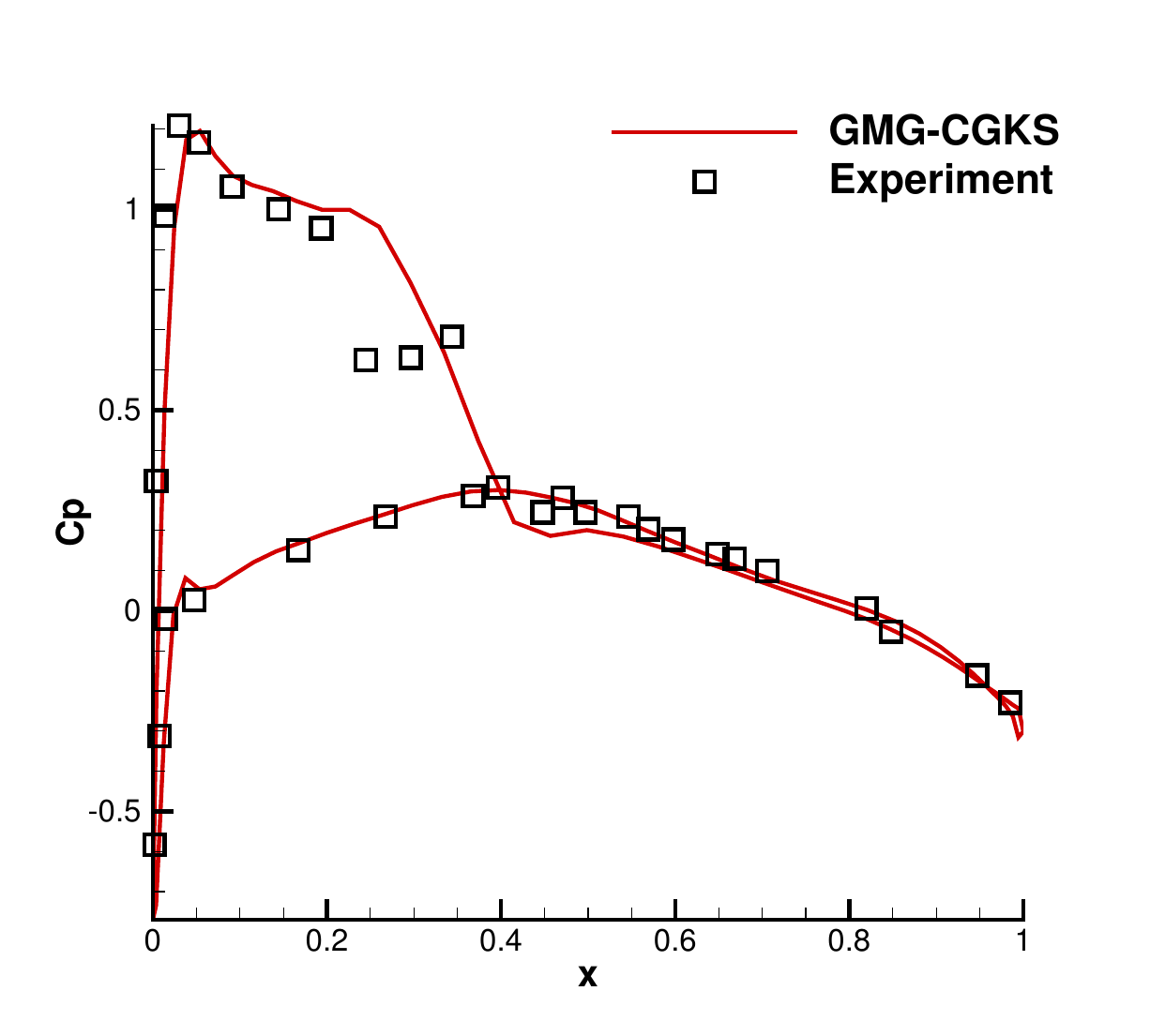}
	\includegraphics[height=0.35\textwidth]{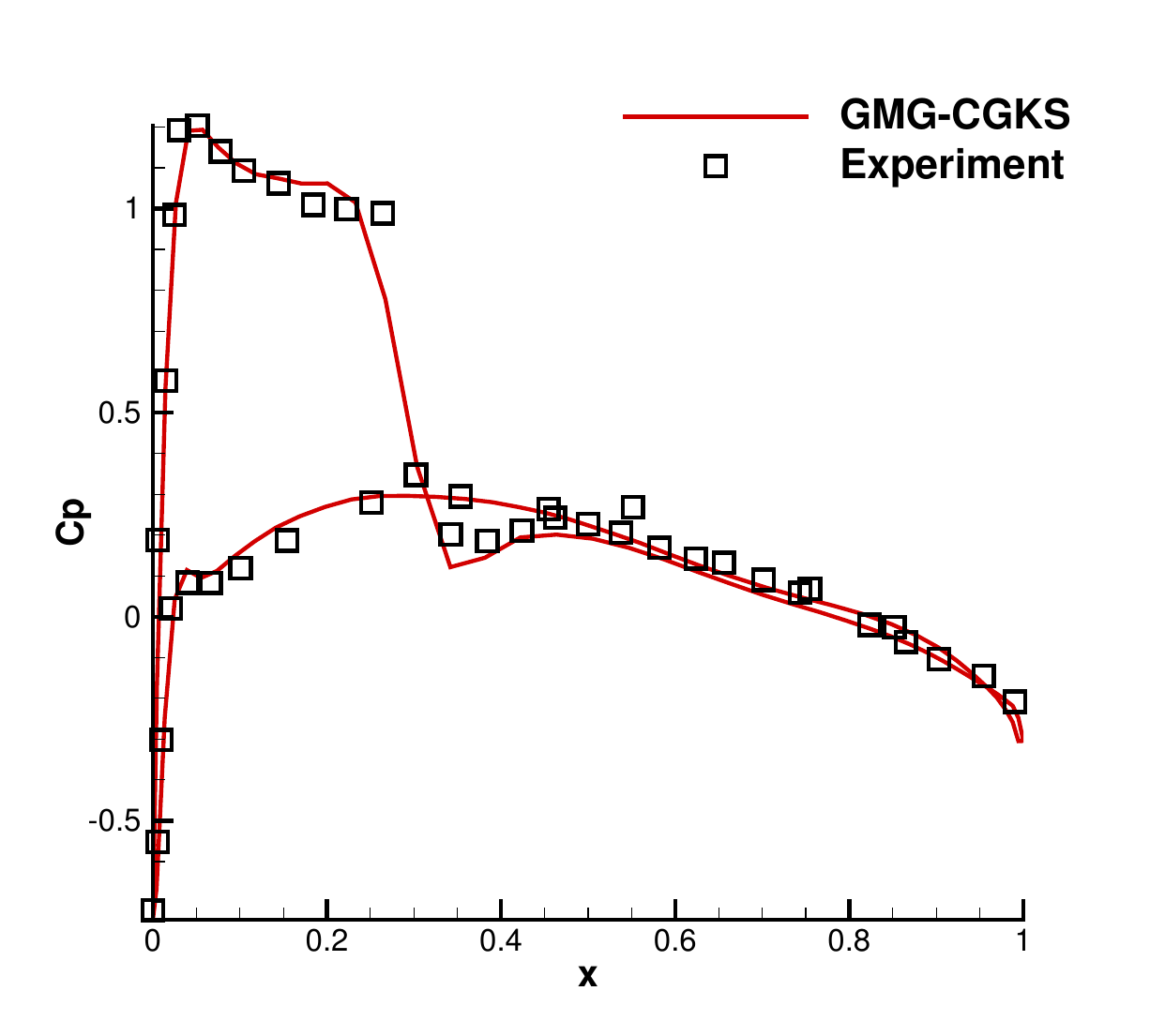}
	\includegraphics[height=0.35\textwidth]{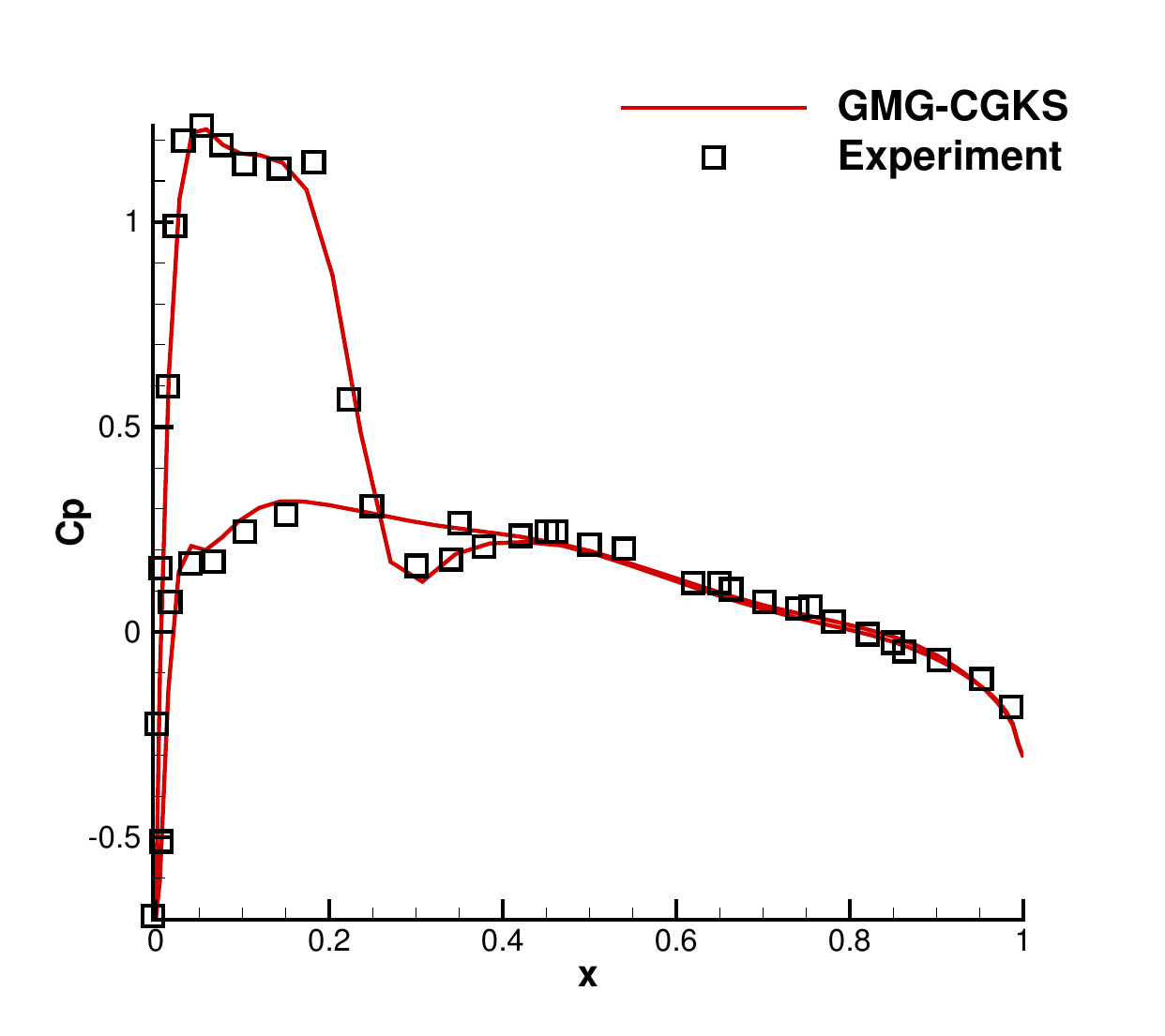}
	\caption{\label{M6 Contour-Cp}
		Transonic flow around a M6 wing. Cp distribution.}
\end{figure}

Transonic flow around an ONERA M6 \cite{eisfeld2006onera} wing is a widely used engineering case to verify the acceleration techniques used in CFD \cite{yang2023implicit}.
The flow structure of it is complicated due to the interaction of shock and wall boundary.
Moreover, three-dimensional mixed unstructured mesh is also a challenge to high-order schemes.
Thus, it is an appropriate test case to verify the accuracy and robustness of the GMG-CGKS.
The far-field Mach number is set to be 0.8395, and the angle of attack is 3.06$^{\circ}$.
The adiabatic slip wall boundary is used on the surface of the ONERA M6 wing, and the subsonic inflow boundary is set according to the local Riemann invariants.
A hybrid unstructured mesh with a near-wall size $h\approx 2e^{-3}$ is used in the computation, as shown in Fig.~\ref{M6 Contour}.
The pressure distribution on the wall surface is shown in Fig.~\ref{M6 Contour}.
The pressure contour in Fig.~\ref{M6 Contour} indicates that the memory reduction CGKS has captured the shock accurately.
The quantitative comparisons on the pressure distributions at the semi-span locations Y /B = 0.20, 0.44, 0.65, 0.80, 0.90, and 0.95 of the wing are given in Fig.~\ref{M6 Contour-Cp}. The numerical results quantitatively agree well with the experimental data.

\begin{figure}
    \centering
    \includegraphics[width=0.35\linewidth]{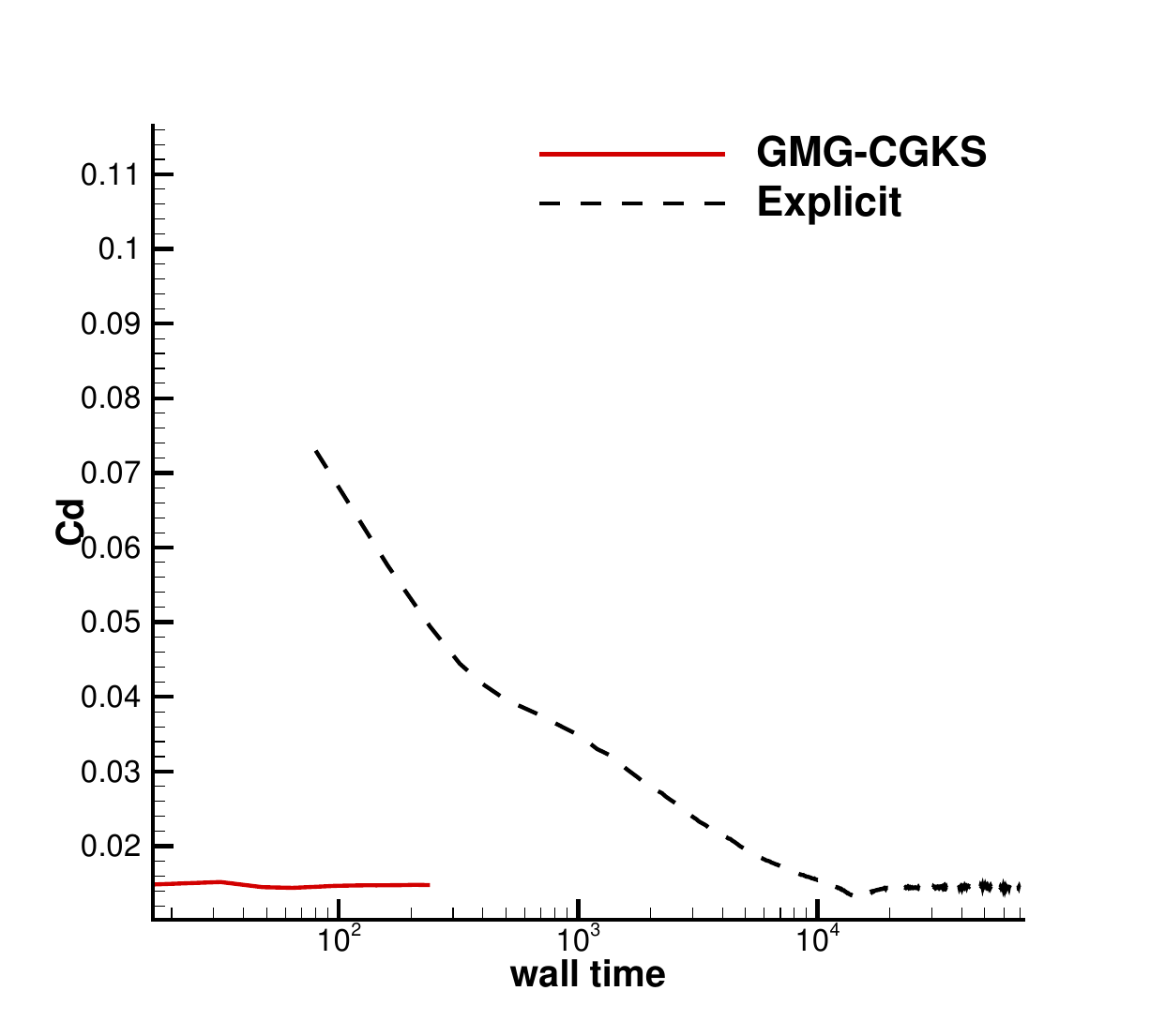}
    \includegraphics[width=0.35\linewidth]{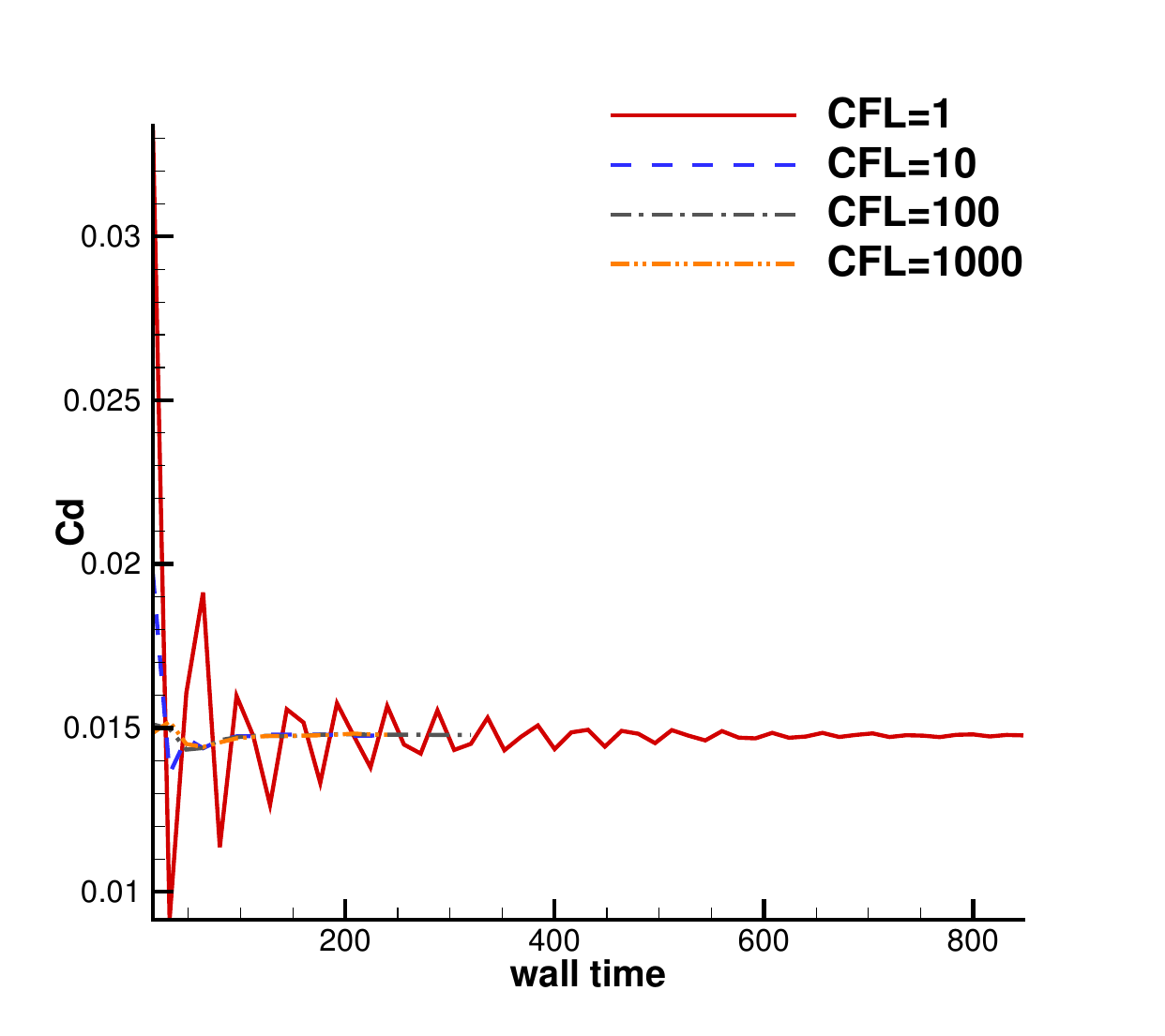}
    \caption{Transonic flow around a M6 wing. Left: The convergence history of GMG and explicit CGKS. Right: The convergence history of GMG under different CFL numbers (based on a fixed sweep number of six of LU-SGS).}
    \label{m6-GMG-1}
\end{figure}

\begin{figure}
    \centering
    \includegraphics[width=0.35\linewidth]{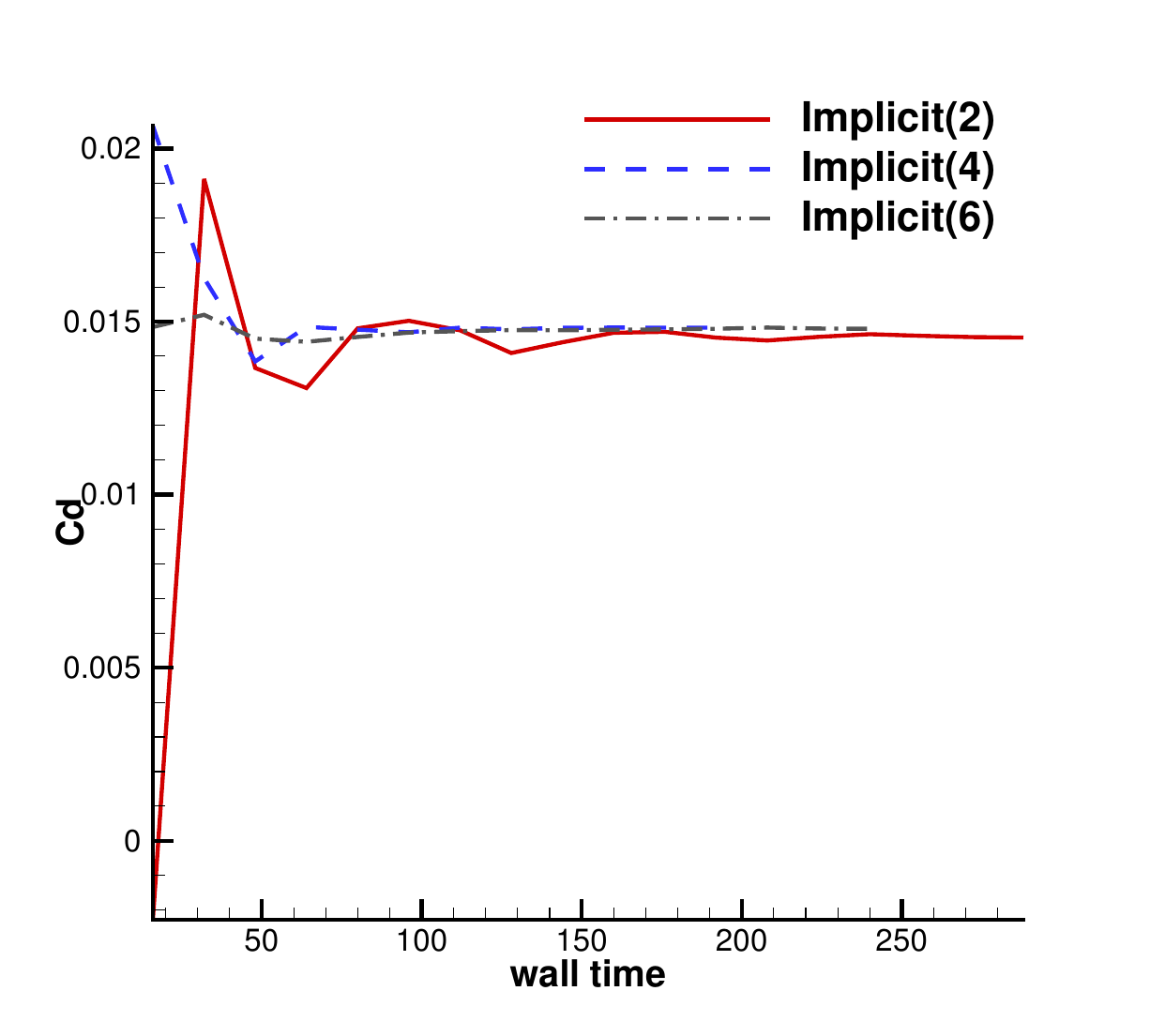}
    \caption{Transonic flow around a M6 wing. The convergence history of GMG under different LU-SGS sweep numbers (based on a fixed CFL number) }
    \label{m6-GMG-2}
\end{figure}

A series of comparisons of convergence rate for this case are shown in Fig .~\ref{m6-GMG-1} and Fig .~\ref{m6-GMG-2}, including a comparison between explicit and GMG, a comparison with a fixed CFL number while varying the LU-SGS sweep number in the coarse grid, a comparison with a fixed LU-SGS sweep number while varying the CFL number, and a comparison with and without DF-based relaxation.
Compared to the explicit scheme, GMG's convergence rate is 100 times faster, requiring only 100 seconds to converge. This is particularly notable given the GPU's already high computational speed, highlighting the significant efficiency of the GPU-accelerated multi-color LU-SGS-based GMG-CGKS.
Regarding comparing CFL numbers, the GMG method does not require a high CFL value (greater than ten is enough) to achieve a fast convergence rate.
At the same time, the comparison results indicate that the sweep count for LU-SGS should not be too low, with four to six sweeps being optimal.
Furthermore, if DF-based relaxation is not used, this case will explode, demonstrating the importance of DF-based relaxation.

\subsection{Hypersonic flow around X-38 type vehicle}
\begin{figure}
    \centering
    \includegraphics[width=0.35\linewidth]{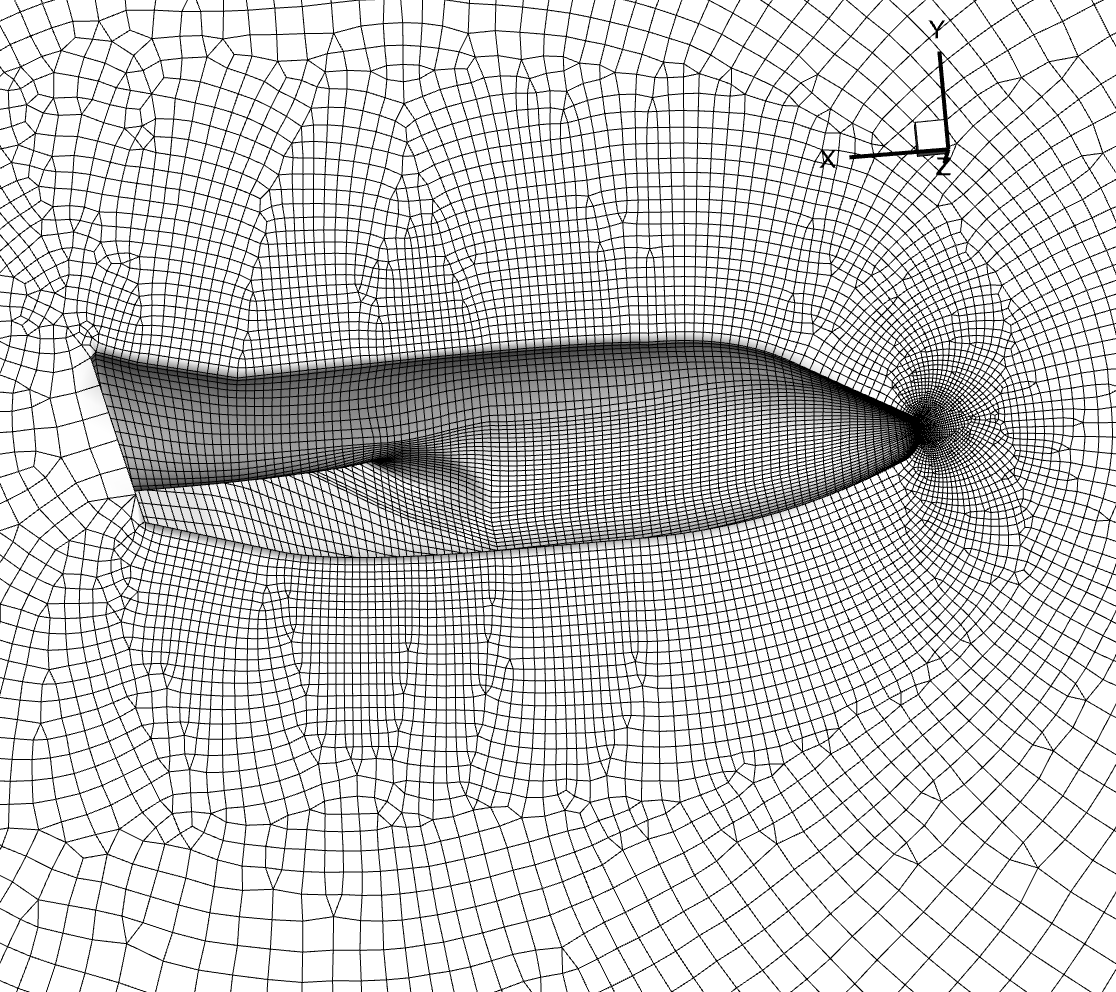}
    \caption{The sketch of meshes of the X-38 type vehicle.}
    \label{x38-mesh}
\end{figure}

\begin{figure}
    \centering
    \includegraphics[width=0.35\linewidth]{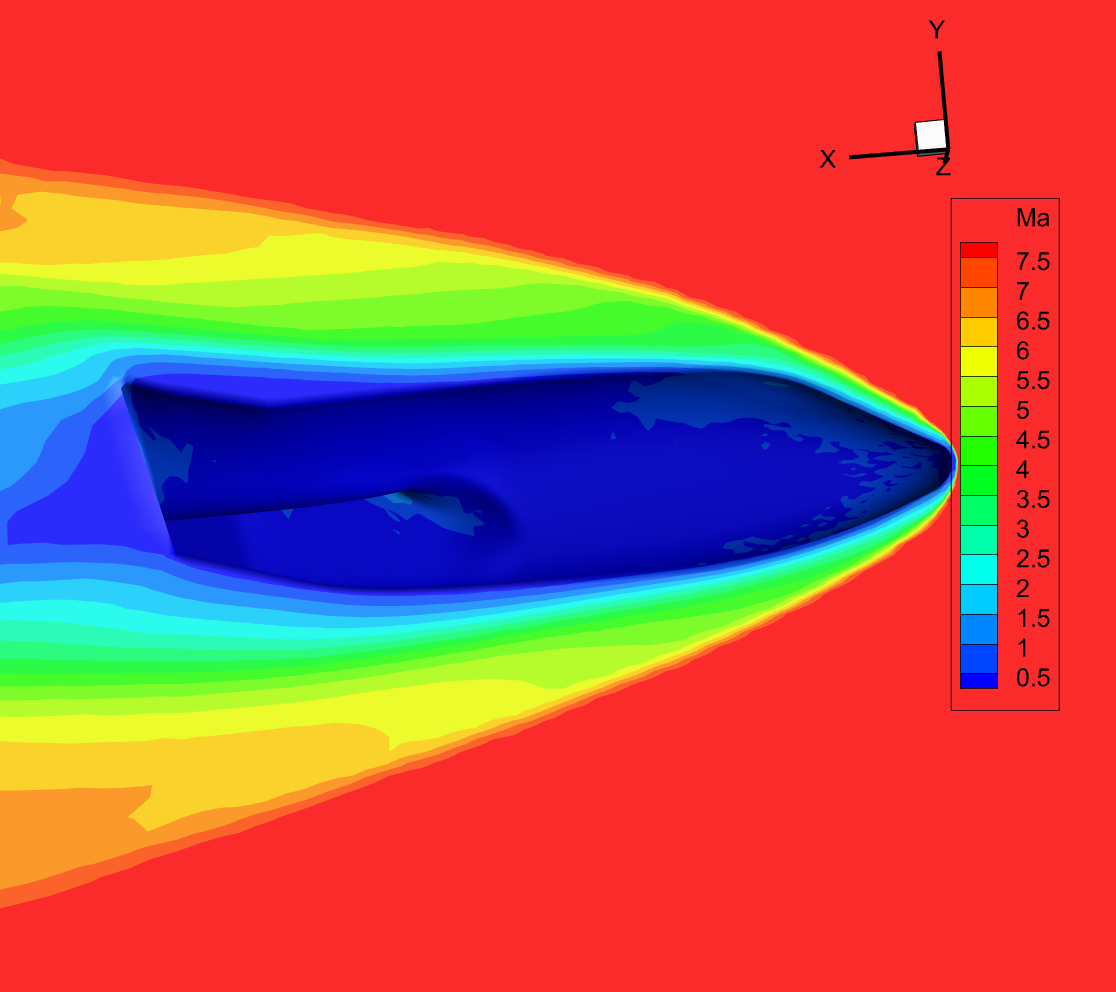}
    \caption{Mach distribution for hypersonic flow around the X-38 type vehicle.}
    \label{x38-ma}
\end{figure}

Hypersonic flow around the X-38 type vehicle is simulated to demonstrate the proposed GMG-CGKS's strong robustness and fast convergence rate.
The incoming Mach number is 8.0, and the Reynolds number is 14289, indicating it's a very tough case even for an explicit time-marching scheme.
The mixed hybrid unstructured mesh includes prim, tetrahedron, and pyramid cells with a total number equal to 560,593, as shown in Fig.~\ref{x38-mesh}.
The Mach number contour is shown in Fig.~\ref{x38-ma}.

\begin{figure}
    \centering
    \includegraphics[width=0.35\linewidth]{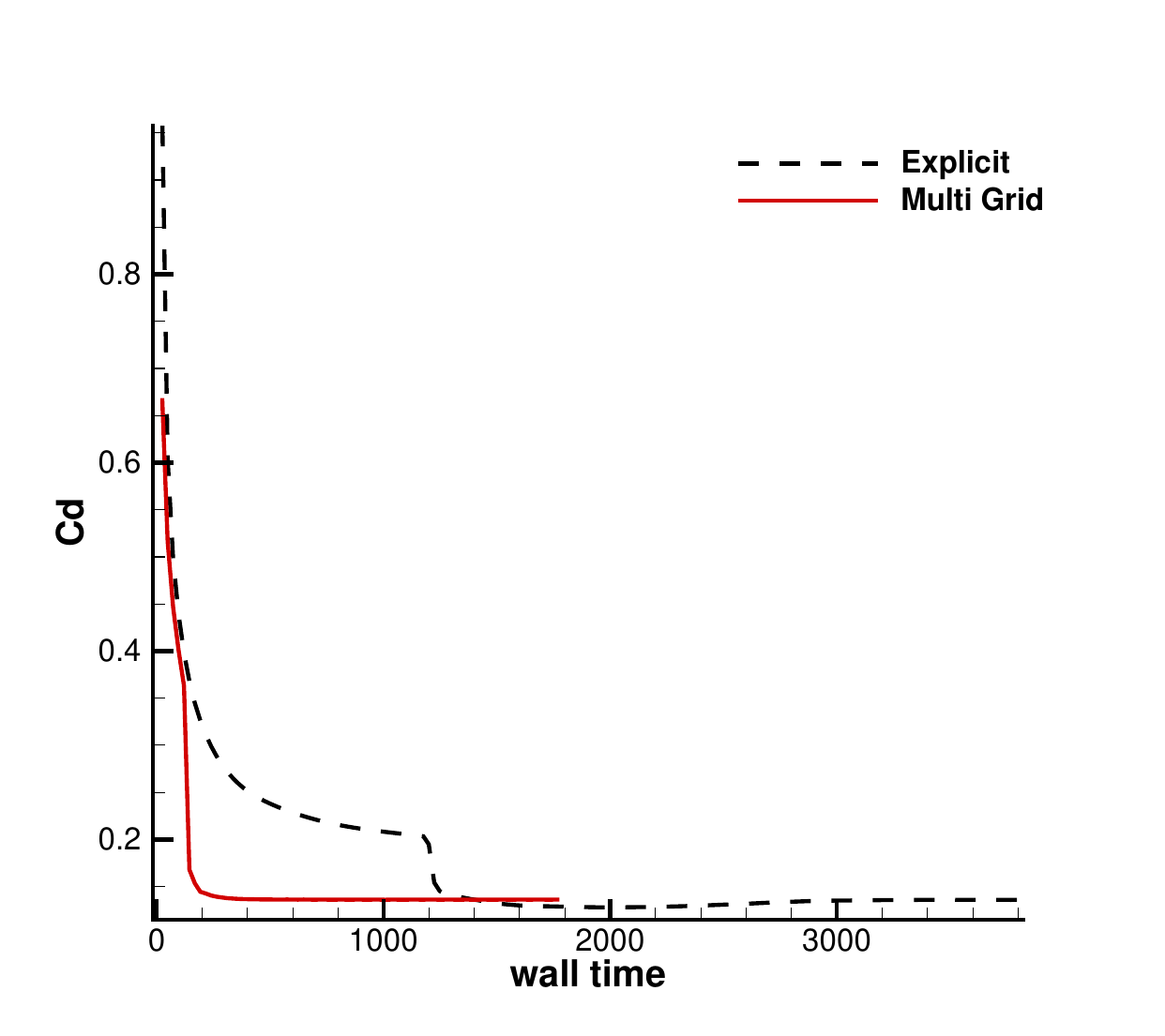}
    \caption{Drag coefficient convergence history of hypersonic flow around the X-38 type vehicle.}
    \label{x38-gmg}
\end{figure}

Drag coefficient convergence history is shown in Fig.~\ref{x38-gmg}, 16 times speedup of convergence rate has been achieved by multi-color LU-SGS based GMG on GPU.
Moreover, compared with the CPU explicit time marching scheme, an 2460 times speedup of the convergence rate is achieved.
This demonstrates that even when facing strong shocks, GMG-CGKS is very robust.
Accelerating the convergence rate is challenging for supersonic and hypersonic flow, but GMG-CGKS still has an 16-time speedup on the GPU when facing complex geometry and strong shock.

In summary, Table 5 summarizes the wall-clock times (in seconds) for each test case, comparing the CPU explicit method, GPU explicit method, and GPU multigrid method. The test cases are as follows: Case 1—subsonic flow around a cylinder; Case 2—subsonic flow around a NACA0012 airfoil; Case 3—transonic flow around a dual-NACA0012 airfoil; Case 4—subsonic viscous flow around a sphere; Case 5—transonic viscous flow around a sphere; Case 6—supersonic viscous flow around a sphere; Case 7—flow over an M6 wing; and Case 8—hypersonic flow around the X-38 vehicle. In these tests, the maximum speedup exceeded 100 times, with a minimum speedup of 13 times compared to the GPU-based explicit method. Notably, in hypersonic cases, the implicit algorithm maintained stable computations while achieving impressive acceleration ratios.

\begin{table}[!htb]\centering
		\begin{tabular}{c|c|c|c|c|c}
		\Xhline{1.5pt}
			\multicolumn{2}{c|}{Case Number} & Case 1 & Case 2 & Case 3  & Case 4   \\ \Xhline{1.5pt}
			CPU explicit  & Wall time  & 255,000 & 61,923 & 89,775  & 103,944  \\ 	\cline{1-6}
			\multirow{2}{*}{GPU explicit}  & Wall time  & 1442 & 368 & 523  & 568  \\ \cline{2-6}
			& Speedup ratio     & 177 & 168 & 171  & 183  \\ \hline
			\multirow{2}{*}{GPU multigrid}  & Wall time  & 45 & 20 & 15  & 20  \\ \cline{2-6}
& Speedup ratio     & 5667 & 3096 & 5985  & 5197  \\
\Xhline{1.5pt}
			\multicolumn{2}{c|}{Case Number}   & Case 5 & Case 6 & Case 7 & Case 8 \\ \Xhline{1.5pt}
CPU explicit  & Wall time   & 8,004,000 & 11,245 & 1,827,280 & 316,800 \\ 	\cline{1-6}
\multirow{2}{*}{GPU explicit}  & Wall time   & 46,000 & 65 & 10040 & 1924 \\ \cline{2-6}
& Speedup ratio      & 174 & 173 & 182 & 164 \\ \hline \cline{2-6}
\multirow{2}{*}{GPU multigrid}  & Wall time   & 500 & 5 & 100 & 120 \\ \cline{2-6}
& Speedup ratio     & 16,008 & 2249 & 18272 & 2640 \\ \hline \cline{2-6}
		\end{tabular}
	\caption{Summary of the speedup ratios (compared with CPU explicit wall time)}
	\end{table}

\section{Conclusion}

This paper presents a GPU-based geometric multigrid algorithm designed to accelerate steady-state computations using CGKS, offering both high convergence rates and strong robustness. An efficient, high-quality coarsening strategy is employed to generate coarse grids within the unstructured multigrid framework.
The primary innovation of this work lies in the use of multi-color LU-SGS to resolve data dependency issues on GPUs, effectively eliminating high-frequency errors on fine meshes and accommodating discontinuities in both space and time. Distinct from traditional multigrid prolongation processes, this study introduces the discontinuity feedback (DF) method—developed in previous work—as an adaptive restriction criterion. This enables the multigrid solver to seamlessly revert to explicit time-stepping in the presence of discontinuities, while eliminating the computational overhead associated with continuously selecting smooth prolongation templates. Moreover, instead of relying on a global, uniform relaxation process, the DF method is fully utilized to adaptively adjust the time-stepping for each cell. This approach accelerates convergence in smooth regions while maintaining robust explicit time-stepping near strong discontinuities.
A comprehensive set of tests on various three-dimensional unstructured grid cases—including flows ranging from nearly incompressible to hypersonic—demonstrates that the GMG-CGKS method accurately captures essential flow features and shows excellent agreement with reference solutions.

\section{Acknowledgements}
The current research is supported by the National Natural Science Foundation of China (12302378, 12172316, 92371201, and
92371107), the Funding of National Key Laboratory of Computational
Physics, the National Key R\&D Program of China (Grant
No. 2022YFA1004500), the Natural Science Basic Research Plan in Shaanxi Province of China (No.
2025SYS-SYSZD-070), and the Hong Kong Research Grant Council (16208021,
16301222, and 16208324)

\color{black}


\bibliographystyle{plain}%
\bibliography{main}

\end{document}